\documentclass{siamltex}
\usepackage{times}
\usepackage{graphics,epsfig,graphicx,color}
\usepackage{amsmath,amssymb}
\usepackage{mathrsfs,amsfonts}

\usepackage{url}
\usepackage[debug,
letterpaper=true,
colorlinks=true,
linkcolor=red,
filecolor=green,
citecolor=red,
pdfpagemode=UseNone]
{hyperref}

\graphicspath{{./FIGS/}}

\newtheorem{rem}[theorem]{Remark}
\newtheorem{ass}[theorem]{Assumption}

\DeclareMathAlphabet{\itbf}{OML}{cmm}{b}{it}
 \DeclareMathAlphabet\mathbfcal{OMS}{cmsy}{b}{n}

\renewcommand{\hat}{\widehat}
\renewcommand{\tilde}{\widetilde}
\def\RR{\mathbb{R}}
\def\bx{{{\itbf x}}}

\def\bu{{{\itbf u}}}
\def\bb{{\itbf b}}
\def\bg{{\itbf g}}
\def\be{{\itbf e}}

\def\bU{{\itbf U}}
\def\bH{{\itbf H}}
\def\bQ{{\itbf Q}}
\def\bGa{{\boldsymbol{\Gamma}}}
\def\bhGa{{\hat{\boldsymbol{\Gamma}}}}
\def\bv{{\itbf v}}
\def\bw{{\itbf w}}

\def\bphi{{\boldsymbol{\varphi}}}

\def\balpha{{\boldsymbol{\alpha}}}

\def\bY{{\itbf Y}}

\def\bV{{\itbf V}}
\def\bR{{\itbf R}}
\def\bI{{\itbf I}}
\def\bM{{\itbf M}}
\def\bS{{\itbf S}}
\def\bW{{\itbf W}}
\def\bD{{\itbf D}}
\def\cP{{\mathcal P}}
\def\cL{{\mathcal L}}
\def\cQ{{\mathcal Q}}
\def\cR{{\mathcal R}}

\def\RM{{\scalebox{0.5}[0.4]{ROM}}}

\def\cPq{{\cP(q)}}
\def\cLq{{\cL(q)}}
\def\cLo{{\cL(0)}}
\def\cPqR{{\boldsymbol{\cP}^{\RM}(q)}}

\def\cLqR{{\boldsymbol{\cL}^{\RM}(q)}}
\def\cLoR{{\boldsymbol{\cL}^{\RM}(0)}}
\def\cLqsR{{\boldsymbol{\cL}^{\RM}(q^S)}}

\def\bbR{{\bb^{\RM}}}

\def\bhu{ \hat{\bu}}
\def\bhuR{{\hat{\bu}^{\scalebox{0.5}[0.4]{ROM}}}}
\def\bvc{{\boldsymbol{\nu}}}
\def\bga{{\boldsymbol{\gamma}}}
\def\bhga{{\hat{\boldsymbol{\gamma}}}}
\def\bhphi{{\hat{\bphi}}}
\def\btphi{{\boldsymbol{\phi}}}
\def\bhtphi{{\hat{\boldsymbol{\phi}}}}
\def\bsigma{{\boldsymbol{\sigma}}}
\def\uu{{\underline{\underline{\bu}}}}
\def\uv{{\underline{\underline{\bv}}}}

\def\PY{{\boldsymbol{\Pi}^{\RM}(q)}}

\def\VY{{\boldsymbol{\mathcal{V}}}}
\def\VYd{\hat \VY}
\def\lb{\left <}
\def\rb{\right >}

\def\cT{\mathcal{T}}
\def\om{\omega}
\def\la{\lambda}

\title{Reduced Order Model Approach to Inverse Scattering} 
\author{
Liliana Borcea\footnotemark[1]
\and
Vladimir Druskin\footnotemark[2]
\and
Alexander V. Mamonov\footnotemark[3]
\and
Mikhail Zaslavsky\footnotemark[4]
\and
J\"{o}rn Zimmerling\footnotemark[5]
}

\begin{document}

\maketitle

\renewcommand{\thefootnote}{\fnsymbol{footnote}}

\footnotetext[1]{Department of Mathematics, University of Michigan,
  Ann Arbor, MI 48109-1043 (borcea@umich.edu)}
\footnotetext[2]{Mathematical Sciences, Worcester Polytechnic Institute,  Worcester, MA 01609-2280(vdruskin@wpi.edu)}
\footnotetext[3]{Department of Mathematics, University of Houston,
  Houston, TX 77004 (mamonov@math.uh.edu)}
 \footnotetext[4]{Schlumberger-Doll Research Center, 1 Hampshire St.,
  Cambridge, MA 02139-1578 (mzaslavsky@slb.com)}
\footnotetext[5]{Department of Mathematics, University of Michigan,
  Ann Arbor, MI 48109-1043 (jzimmerl@umich.edu)}  

\begin{abstract}
We study an inverse scattering problem for a  generic hyperbolic  system of equations with an unknown
coefficient called the reflectivity. The solution of the system models waves (sound, electromagnetic or elastic),
and the reflectivity models unknown scatterers embedded in a smooth and known medium.  The inverse problem is to determine
the reflectivity from the time resolved scattering matrix (the data) measured by an array of sensors. We introduce a 
novel inversion method, based on a reduced order model (ROM) of an operator called wave propagator, because it  maps  the wave from one time 
instant to the next, at interval corresponding to the discrete time sampling of the data. 
The wave propagator is  unknown in the inverse problem,  but the ROM can be computed directly from the data.  By construction, the ROM inherits key properties of the wave propagator, which facilitate the estimation of the 
reflectivity. The ROM was introduced previously and was used for two purposes: (1) to map the scattering 
matrix to that corresponding to the single scattering (Born) approximation and (2) to image i.e., obtain a qualitative estimate 
of the support of the reflectivity. Here we study further the ROM and show that it corresponds to a Galerkin
projection of the wave propagator. The Galerkin framework is useful for proving properties of the ROM that 
are used in the new inversion method which seeks a quantitative estimate of the reflectivity.
\end{abstract}

\begin{keywords}
Inverse scattering, model reduction, Galerkin approximation.
\end{keywords}

\begin{AMS}
65M32, 41A20
\end{AMS}

\section{Introduction}
\label{sect:intro}
Consider an inverse scattering problem for a hyperbolic system of equations in symmetric form
\begin{align}
\partial_t^2 \bu(t,\bx)+ L(q) L(q)^T  \bu(t,\bx) &= 0, \quad \bx \in \Omega, ~ t > 0,  \label{eq:F1} \\
\bu(0,\bx) &= \bb(\bx), \quad \bx \in \Omega, \label{eq:F2} \\
\partial_t \bu(0,\bx) &= {\bf 0}, \label{eq:F3}
\end{align}
satisfied by the wave $\bu(t,\bx)$, where $t$ denotes time and $\bx$ is the spatial variable
in the domain $\Omega \subset \RR^d$ in dimension $d \ge 1$, with piecewise smooth boundary $\partial \Omega$.
The information about the medium is in the operator $L(q)$ and its adjoint $L(q)^T$, defined on
spaces of functions satisfying some  homogeneous boundary conditions.
Both $L(q)$ and $L(q)^T$ are first order partial differential operators in the variable $\bx$, with affine dependence on the unknown
coefficient $q(\bx)$, called the reflectivity. The inverse problem is 
to determine $q(\bx)$ from data gathered by a collection (array) of sensors.  This probes the medium with incident waves, determined by the initial condition $\bb(\bx)$, and measures the backscattered waves.

\begin{figure}[t!]
\centering 
\includegraphics[width=0.25\textwidth]{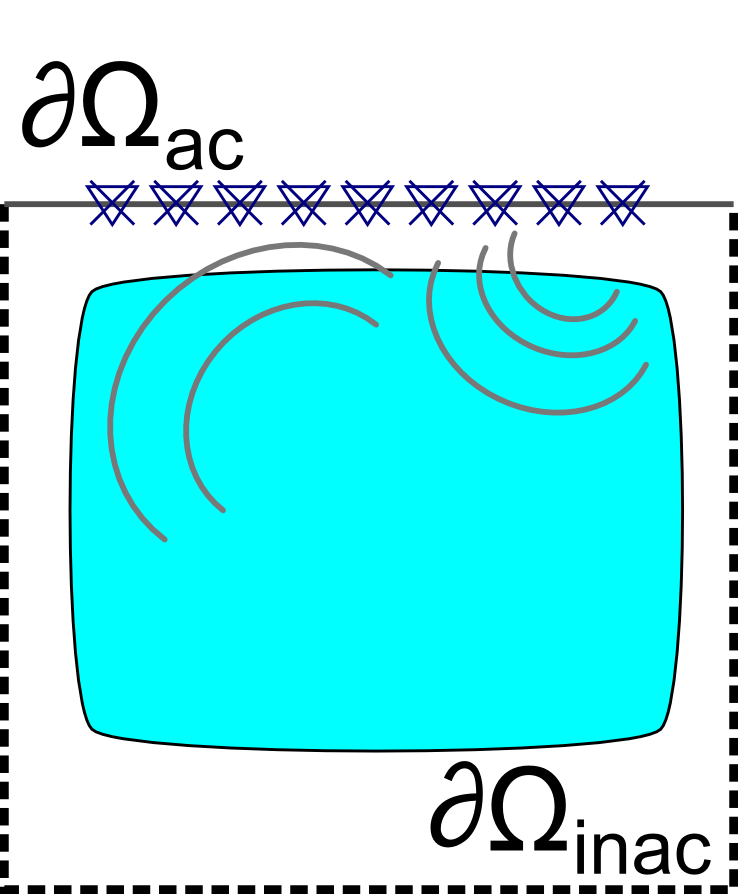} 
\vspace{-0.0in}\caption{Illustration of the setup: An array of sensors (indicated with crosses) lying near the accessible boundary 
$\partial \Omega_{\rm ac}$ probes an unknown  medium with incident waves and measures the backscattered waves.  The inaccessible 
boundary  $\partial \Omega_{\rm inac}$ is far enough from the sensors, so that it has no effect during the duration of the  measurements.}
\label{fig:setup}
\end{figure}

Problem  (\ref{eq:F1}--\ref{eq:F3}) arises in inverse scattering for sound, electromagnetic and elastic waves in isotropic media, as explained in \cite[sections 3--5]{borcea2019robust}.  In acoustics, $\bu(t,\bx)$ is  related via some transformation to the scalar valued acoustic pressure, whereas in electromagnetics and  elasticity, $\bu(t,\bx)$  is related to the vector valued electric field and displacement velocity, respectively. The medium is modeled by variable coefficients in the wave equations:  the wave speeds and wave impedances. Depending on  the  data acquisition setup, these coefficients affect in a different way the  measurements at the array.  Our definition  of the reflectivity $q(\bx)$  takes this into account, as we now explain.

We consider the setup illustrated in Figure \ref{fig:setup},  where $\Omega$ is a cube in $\RR^d$, obtained via truncation of a half space occupied by the unknown medium\footnote{One can also consider truncation of the whole space, as long as the medium is known and non-scattering on one side of the array of sensors.}. Assuming that the sensors record over the duration $t \in [0,T]$, and using that the wave speed is finite, we let the cube $\Omega$ have large enough side length, so that the measurements are not affected by the medium outside $\Omega$. The boundary $\partial \Omega = \partial \Omega_{\rm ac} \cup \partial \Omega_{\rm inac}$ is the union of the 
accessible boundary $\partial \Omega_{\rm ac}$, which is a subset of the boundary of the half space,  and the inaccessible boundary $\partial \Omega_{\rm inac}$. The name accessible means that the array of sensors
can be placed in the immediate vicinity of 
$\partial \Omega_{\rm ac}$. The inaccessible boundary is fictitious 
and has no effect on the  measurements, so the backscattered wave is due entirely to reflectors contained in $\Omega$. 
The initial condition in~\eqref{eq:F2} is a vector valued function 
\begin{equation}
\bb(\bx) = \Big(b^{(1)}(\bx), \ldots, b^{(m)}(\bx)\Big),
\label{eq:F6}
\end{equation}
where $b^{(s)}(\bx)$ is the wave emitted by one sensor\footnote{Note that typically, the source excitation is expressed 
as a time dependent force in the right hand-side of the wave equation, with homogeneous initial conditions. 
We refer to appendix \ref{ap:source} for the derivation of the initial value problem \eqref{eq:F1}--\eqref{eq:F3} 
from such a formulation. We also give there the expression of $\bb(\bx)$
which depends on the waveform emitted by the source.}. 
It is a function supported in the vicinity  of the sensor and the index $s = 1,\ldots, m$ counts the sensors and the 
polarization of the wave.

The array measures the  scattering matrix (the data),  modeled by \cite[sections 3-5]{borcea2019robust}
\begin{equation}
\bD_j = \lb \bb, \bu(j \tau,\cdot) \rb = \int_{\Omega} d \bx \, \bb(\bx)^T \bu(j \tau,\bx), \quad j = 0, \ldots, 2n-1.
\label{eq:F7}
\end{equation}
The  $s^{\rm th}$ column of this symmetric $m \times m$ matrix corresponds to the wave generated by  $b^{(s)}(\bx)$ and evaluated at all the sensors in the array, at time instant $j \tau$, where $\tau > 0$ is chosen consistent with the Nyquist sampling rate of the 
wave.

The wave $\bu(t,\bx)$ and therefore the data~\eqref{eq:F7}  depend in a complicated, nonlinear way on the coefficients (wave speed and impedance) of the wave equation. The low spatial frequency component of the wave speed determines the kinematics of the wave
\cite{symes2009seismic,bleistein2013mathematics}, since time of travel is a path integral of the slowness (the reciprocal of the velocity). The estimation of this smooth part (aka the kinematic model) is of great interest in geophysical exploration \cite{ symes2009seismic}. It is a difficult problem in the backscattering setup considered here and at high frequencies used in applications, because nearby models can 
give travel time discrepancies that exceed the short period of oscillation of the wave. Thus, unless data have low temporal frequencies,  typical least squares data fit optimization formulations \cite{tarantola1984inversion} are not amenable to solutions by Newton-type methods \cite{virieux2009overview}. Other approaches have emerged \cite{symes2008migration}, and they use redundant data sets  to separate the estimation of the kinematic model and the rough, backscattering part  of the medium, called the reflectivity. We assume that the kinematic model is known\footnote{The kinematic model (smooth wave speed) 
appears in the coefficients of the operators $L(q)$ and $L(q)^T$ (see \cite{borcea2019robust} and sections \ref{sect:Inv}--\ref{sect:Numer}). We suppress  the dependence on the known kinematic model in our notation.}
and is such that the wave front advances forward (there are no lensing effects). Then, the study in \cite{BeylkinBurridge} shows that if the 
depth of the reflectors is larger than the diameter  of the array,  backscattering is mostly due to 
relative variations of the wave impedance. This motivates our definition of $q(\bx)$ as the logarithm of the impedance \cite[sections 3-5]{borcea2019robust}. 

The estimation of the reflectivity from backscattering data i.e., inverting the mapping 
\begin{equation}
q \mapsto \{\bD_j, ~ j=0,\ldots, 2n-1\},
\label{eq:refDat}
\end{equation}
 has applications in nondestructive testing 
\cite{schmerr2016fundamentals}, ultrasound for medical diagnostics \cite{szabo2004diagnostic}, radar \cite{cheney2009fundamentals}, 
geophysical exploration \cite{symes2009seismic}, underwater sonar \cite{collins1994inverse}, and so on.  
It is a nonlinear inverse problem, even though $L(q)$ is affine in $q(\bx)$,  
as can be seen by solving (\ref{eq:F1}--\ref{eq:F3}) 
\begin{equation}
\bu(t,\bx) = \cos\Big(t \sqrt{L(q) L(q)^T}\Big) \bb(\bx), \qquad t > 0,
\label{eq:F8}
\end{equation}
and substituting the solution in the data model~\eqref{eq:F7}
\begin{equation}
\bD_j = \lb \bb, \cos\Big(j \tau \sqrt{L(q) L(q)^T}\Big)\bb \rb, \quad j = 0, \ldots, 2n-1,
\label{eq:F9}
\end{equation}
where the square root and cosine are defined as usual, using the spectral decomposition of the self-adjoint, 
nonnegative-definite operator $L(q) L(q)^T$. Basically all existing algorithms search for the reflectivity with 
a least squares data fit optimization formulation, and in many applications the mapping~\eqref{eq:refDat} 
is linearized i.e., $q(\bx)$ is estimated by the solution of the normal equation. The normal operator is not 
invertible in general, but in many setups it has the property that it preserves approximately the location of 
non-smooth features of $q(\bx)$, like jumps \cite{symes2009seismic}. Therefore, popular methods like 
reverse time migration \cite{claerbout1985imaging,biondi20063d,bleistein2013mathematics} and the related 
backprojection \cite{beylkin1985imaging,cheney2009fundamentals} use the right hand side of the normal 
equation as an image, i.e., an estimate of the support of $q(\bx)$. These imaging methods work well if the 
reflectivity $q(\bx)$ is weak, but they are qualitative. A quantitative estimate of a general reflectivity requires 
inverting, in an appropriate sense, the nonlinear map~\eqref{eq:refDat}.

We propose a method for estimating $q(\bx)$ based on a reduced order model (ROM) of the 
self-adjoint wave propagator operator
\begin{equation}
\cPq = \cos\Big(\tau \sqrt{L(q) L(q)^T}\Big).
\label{eq:Sn5}
\end{equation}
This operator is useful because it allows us to view the wave ${\bf u}(j \tau, \bx)$ as the state 
of  a discrete dynamical system, starting from $\bb(\bx)$ and evolving  with the time index $j \ge 0$.  
We can write explicitly the state $\bu(j \tau,\bx) = \cT_j(\cPq) \bb(\bx)$ using equations \eqref{eq:F8} and \eqref{eq:Sn5}, 
and substituting in the expression  (\ref{eq:F9}) of the data we obtain 
\begin{equation}
\bD_j = \lb \bb, \cos\big(j \arccos \cPq \big)\bb \rb = \lb \bb, \cT_j \big( \cPq \big)\bb \rb,
\quad j = 0, \ldots, 2n-1,
\label{eq:dataT}
\end{equation}
where $\cT_j$ are Chebyshev polynomials of the first kind \cite{rivlin1974chebyshev}.
The ROM is defined by a pair of matrices $\cPqR \in \RR^{nm \times nm}$ and $ \bbR \in \RR^{nm\times m}$, 
which are proxies of  $\cPq$ and $\bb(\bx)$, in the sense that they define a dynamical  system for the 
discrete state $\cT_j \big( \cPqR \big) \bbR$, that encodes essential features of $\bu(j \tau,\bx)$ and satisfies
\begin{equation}
\bD_j = \lb \bb, \cT_j \big( \cPq \big)\bb \rb = \bbR^T \cT_j \big( \cPqR \big) \bbR, \quad j = 0, \ldots, 2n-1.
\label{eq:ROMThm1}
\end{equation}

The matrices $\cPqR$ and $\bbR$ satisfying (\ref{eq:ROMThm1}) are calculated from the data \eqref{eq:F9} 
(i.e., the ROM is data-driven) and they capture physical aspects of the wave propagation that are needed for inversion. 
The ROM was introduced in \cite{druskin2016direct,borcea2019robust} and was used in \cite{druskin2018nonlinear} 
for imaging, and in \cite{DtB} for transforming the data~\eqref{eq:F9} to that corresponding to the single scattering 
(Born) approximation. The new results in this paper are: 

\vspace{0.05in}
\noindent 1.  We show that the ROM propagator $\cPqR$ is a Galerkin projection of the  operator~\eqref{eq:Sn5}, and use 
the Galerkin framework to prove properties of the ROM that facilitate the solution of the inverse scattering problem.

\vspace{0.05in}
\noindent 2. We use the ROM to develop a novel, quantitative inversion method for  estimating $q(\bx)$. The data are fit implicitly in our method, and the inversion is formulated as a minimization of the discrepancy 
between the data-driven ROM  and the ROM calculated for the search reflectivity.  
This optimization problem turns out to be  almost linear i.e.,  it can be solved 
in very few iterations, and it is better conditioned than the classic least squares data fit approach.

\vspace{0.05in}
The paper is organized as follows: In section \ref{sect:Galerk} we describe the ROM in the Galerkin framework and analyze its 
properties, which are then used in section \ref{sect:Inv} to introduce the new inversion method. 
We assess the performance of the method with numerical simulations in section \ref{sect:Numer}. The presentation in section 
\ref{sect:Galerk} does not depend on the expression of the operator $L(q)$ and its adjoint, 
so we work with the generic hyperbolic system (\ref{eq:F1}--\ref{eq:F3}). However, the  inversion algorithm 
requires specifying $L(q)$, so in sections  \ref{sect:Inv}--\ref{sect:Numer}  we use the operator  derived from the acoustic wave equation.  We end with a summary in section \ref{sect:Sum}.
\section{The Galerkin framework}
\label{sect:Galerk}
To introduce the  Galerkin framework, consider the approximation space 
\begin{equation}
\mathfrak{X} = \mbox{colspan}\{\bu(j\tau,\bx), ~j = 0, \ldots, n-1\},
\label{eq:F5}
\end{equation}
where 
\begin{equation}
\bu(j \tau,\bx) = \Big( u^{(1)} (j \tau, \bx), \ldots, u^{(m)}(j\tau, \bx) \Big), \qquad j \ge 0,
\label{eq:snap}
\end{equation}
are the solution snapshots, with components 
\begin{equation}
 u^{(s)}(j\tau, \bx) = \cos\Big(j \tau \sqrt{L(q) L(q)^T}\Big) b^{(s)}(\bx), \qquad s = 1, \ldots, m.
 \label{eq:snap_s}
 \end{equation}

We begin in section \ref{sect:Galerk.1} with an exact time stepping scheme satisfied by these snapshots, which 
shows the role of the propagator operator~\eqref{eq:Sn5}.
The ROM is defined from the Galerkin approximation of this time stepping scheme, 
as explained in  section \ref{sect:Galerk.2}. Note that the approximation 
space $\mathfrak{X}$ is not known in inversion, because the data~\eqref{eq:F7} correspond to the snapshots 
evaluated at  the locations of the sensors, and not inside the medium.  Nevertheless, it is possible to compute  the ROM, as explained in section \ref{sect:Galerk.3}. The properties of the ROM are stated  in section \ref{sect:Galerk.4} and are proved in appendices \ref{ap:DatFit}--\ref{ap:A}. 
In section \ref{sect:Orthog} we explain that there is a family of ROMs that share these properties, and that they are connected by special orthogonal transformations. We also give in section \ref{sect:Orthog1} an intuitive, finite differences  interpretation  of the ROM,  which is then used in section \ref{sect:Inv} to  motivate the inversion algorithm.

\subsection{The propagator and time stepping}
\label{sect:Galerk.1}
Let us introduce the notation 
\begin{equation}
\bu_j(\bx) = \bu(j \tau,\bx), \qquad j \ge 0,
\label{eq:snapj}
\end{equation}
and obtain from the definition (\ref{eq:snap}-\ref{eq:snap_s}) and the trigonometric identity satisfied by the cosine, that the 
snapshots satisfy the time stepping scheme 
\begin{align}
\bu_{j+1}(\bx) &= 2 \cPq \bu_{j}(\bx) - \bu_{j-1}(\bx), \quad j \ge 0, \label{eq:Sn2} \\
\bu_0(\bx)&= \bb(\bx), \label{eq:Sn3} \\
\bu_1(\bx)&= \bu_{-1}(\bx). \label{eq:Sn4}
\end{align}
This justifies calling $\cPq$ the wave propagator operator, because it is used  to map the wave at  consecutive time instants $(j-1)\tau$ and $j \tau$, to the wave at future time $(j+1)\tau$.

Note that~\eqref{eq:Sn2} is the three term recursion relation satisfied by the orthogonal Chebyshev polynomials 
$\cT_j$ 
and indeed, definitions~\eqref{eq:Sn5} and (\ref{eq:snap}-\ref{eq:snap_s})  give
\begin{equation}
\bu_j(\bx) = \cos\big(j \arccos \cPq\big) \bb(\bx) = \cT_j\big(\cPq\big) \bb(\bx), \quad j \ge 0.
\label{eq:Sn6}
\end{equation}

Note also that  if we subtract $2 \bu_j(\bx)$ from equation~\eqref{eq:Sn2} and divide the result by $\tau^2$ we obtain   the second order 
time stepping scheme 
\begin{align}
\frac{\bu_{j+1}(\bx)-2 \bu_j (\bx)+ \bu_{j-1}(\bx)}{\tau^2}   + \cLq \cLq^T \bu_j(\bx) &= 0, \quad j \ge 0,  ~ \bx \in \Omega, \label{eq:Sn7} \\
\bu_0(\bx)&= \bb(\bx), \label{eq:Sn8} \\
\bu_1(\bx) &= \bu_{-1}(\bx), \label{eq:Sn9}
\end{align}
where $\cLq$ is the square root of the self-adjoint, non-negative  definite operator
\begin{equation}
\frac{2}{\tau^2}\big(I - \cPq\big) = \cLq \cLq^T,
\label{eq:Sn10}
\end{equation}
and  $I$ denotes the identity.  Equations (\ref{eq:Sn7}--\ref{eq:Sn9}) are an exact time stepping scheme for the hyperbolic problem (\ref{eq:F1}--\ref{eq:F3}), with boundary conditions taken into account in the definition of  $L(q)$ and $L(q)^T$. The derivative $\partial_t^2 \bu(j \tau, \bx)$ is replaced in~\eqref{eq:Sn7}  by  second order centered differences, and the $O(\tau^2)$ error is absorbed in the operator $\cLq$ which approximates $L(q)$. That is to say,  $\cLq \cLq^T$  has the same eigenfunctions as $L(q)L(q)^T$, and the eigenvalues
\[
\frac{2}{\tau^2} \big[ 1 - \cos \big(\tau \sqrt{\theta}\big)\big] = \theta \Big[1 + O\big(\tau^2 \theta\big)\Big], 
\]
where $\theta > 0$ denotes an eigenvalue of $L(q)L(q)^T$. The larger part of the spectrum corresponds to  more oscillatory eigenfunctions\footnote{For the purpose of the explanation, we may think of $L(q)L(q)^T$ as the negative Laplacian  multiplied by $c_o^2$.}, 
so the restrictions of $\cLq \cLq^T$ and $L(q)L(q)^T$ on the space of functions that oscillate at smaller spatial  frequency $\lesssim 1/(c_o \tau)$, where $c_o$ is a 
reference wave speed, are approximately the same.

\subsection{Galerkin approximation}
\label{sect:Galerk.2}
We define the ROM using the Galerkin  approximation of  (\ref{eq:Sn2}--\ref{eq:Sn4}), under the following assumption:

\vspace{0.05in}
\begin{ass}
The solution snapshots (\ref{eq:snap}-\ref{eq:snap_s}) are linearly independent up to time $n\tau$. This implies  that the approximation space~\eqref{eq:F5} has  dimension $n m$.
\label{as.1}
\end{ass}
\vspace{0.05in}

The linear independence of the snapshots can be ensured initially by having sufficiently well separated sensors in the array (recall that the components of~\eqref{eq:F6} are approximations of the  delta  function at the sensor locations). However,  depending on the kinematic model, the waves may focus at later time in some region of the domain, or they may turn around, and the snapshots can become linearly dependent. 
We assume in the analysis that the medium is nice enough such that Assumption \ref{as.1} holds, but in the ROM construction and inversion  we can deal with a lower dimensional approximation space  using an SVD truncation. 

Let us gather  the first $n$ snapshots in the  quasimatrix
\begin{equation}
\bU(\bx) = (\bu_0(\bx), \ldots, \bu_{n-1}(\bx)),
\label{eq:ROM8}
\end{equation}
with $n m $ linearly independent columns that span the approximation space ~\eqref{eq:F5}. 
Following \cite{stewart1998afternotes},  by quasimatrix we mean a row vector valued function
defined on  $\Omega$, with the entry index viewed as a column index and $\bx \in \Omega$ 
playing the role of a row index. The quasimatrix  $\bU(\bx)$ is organized in \eqref{eq:ROM8}  in 
$n$ blocks, each with $m$ entries.

Using linear algebra terminology, we write 
\begin{equation}
\label{eq:apSpace}
\mathfrak{X} = \mbox{range} \big( \bU(\bx)\big),
\end{equation}
and we call henceforth  the components of quasimatrices like \eqref{eq:ROM8}  block columns.
The Galerkin approximation of the snapshots is 
\begin{equation}
\bu_j(\bx) \approx \bU(\bx) \bg_j, \quad j \ge 0,
\label{eq:ROM9}
\end{equation}
where $\bg_j \in \RR^{nm \times m}$ is the matrix of Galerkin coefficients. These are calculated so that when substituting~\eqref{eq:ROM9} in 
\eqref{eq:Sn2}, the residual is orthogonal to the space~\eqref{eq:apSpace}
\begin{equation}
\lb \bu_l, \bU \bg_{j+1} - 2 \cPq \bU \bg_j + \bU \bg_{j-1}  \rb = {\bf 0}, \quad \forall ~ l = 0, \ldots, n-1,
\label{eq:ROM10}
\end{equation}
where we recall the definition of $\lb \cdot, \cdot \rb$ from~\eqref{eq:F7}.

By construction, the approximation ~\eqref{eq:ROM9} is exact for $j = 0, \ldots, n-1$, so
\begin{equation}
\bg_j = \be_{j} = 
(\boldsymbol{0}_m, \ldots, \boldsymbol{0}_m, \bI_m, \boldsymbol{0}_m, \ldots, \boldsymbol{0}_m)^T, 
\quad  j = 0, \ldots, n-1,
\label{eq:ROM11}
\end{equation}
are the matrices of size $nm \times m$ with an $m \times m$ identity $\bI_m$ at block position 
$j$, and all other blocks being $m \times m$ zero matrices $\boldsymbol{0}_m$.
\footnote{Note that for convenience we count the block entries starting from $0$.} 
Using this observation and rewriting~\eqref{eq:ROM10} in matrix form, we obtain the Galerkin time stepping scheme 
\begin{align}
\bM \bg_{j+1} &= 2 \bS \bg_j - \bM\bg_{j-1}, \quad j \ge 0, \label{eq:ROM14} \\
\bg_0 &= \be_0, \label{eq:ROM15} \\
\bg_{1} &= \bg_{-1} = \be_1, \label{eq:ROM16}
\end{align}
with mass matrix 
\begin{equation}
\bM = \bU^T \bU \in \RR^{nm \times nm},
\label{eq:ROM12}
\end{equation}
and stiffness matrix 
\begin{equation}
\bS = \bU^T \cPq \bU \in \RR^{nm \times nm}.
\label{eq:ROM13}
\end{equation}
Again, we use  linear algebra notation, where  for all  $\bW(\bx) = (\bw_0(\bx), \ldots, \bw_{n-1}(\bx))$ lying in the same space as $\bU(\bx)$, we denote by 
$\bU^T \bW$   the $nm \times nm$ matrix 
with $m\times m$ blocks 
\begin{equation}
\big(\bU^T \bW \big)_{i,j} =  \lb \bu_i,\bw_j \rb, \qquad i,j = 0,\ldots, n-1.
\label{eq:ROM13p}
\end{equation}

\subsubsection{Definition of the ROM}
We conclude from Assumption \ref{as.1} and definition~\eqref{eq:ROM12} that the mass matrix is symmetric and positive 
definite, so we can take its square root using the block Cholesky factorization  \cite[Chapter 4]{golubVanLoan},
\begin{equation}
\bM = \bR^T \bR,
\label{eq:ROM23}
\end{equation}
where $\bR \in \RR^{nm \times nm}$ is block upper triangular, with $m\times m$ blocks.  This matrix is invertible, and multiplying~\eqref{eq:ROM14} 
on the left by $\bR^{-T}$ (the transpose of the inverse of $\bR$), we obtain the ROM version of the time stepping scheme (\ref{eq:Sn2}--\ref{eq:Sn4}),
\begin{align}
\bu^{\RM}_{j+1} &= 2 \cPqR\bu^{\RM}_j- \bu^{\RM}_{j-1}, \quad j \ge 0, \label{eq:ROM28} \\
\bu^{\RM}_0 &= \bbR, \label{eq:ROM29} \\
\bu^{\RM}_{1} &= \bu^{\RM}_{-1} = \bR \be_1, \label{eq:ROM30}
\end{align}
satisfied by the ROM snapshots
\begin{equation}
\bu^{\RM}_j = \bR \bg_j, \qquad j \ge 0.
\label{eq:ROM31}
\end{equation}
The ROM propagator is the $nm \times nm$ symmetric matrix 
\begin{equation}
\cPqR = \bR^{-T} \bS \bR^{-1},
\label{eq:ROMProp}
\end{equation}
and the initial snapshot is the $nm \times m$ matrix 
\begin{equation}
\bbR = \bR \bg_0 = \bR \be_0.
\label{eq:ROMb}
\end{equation}

\subsubsection{Galerkin projection}
\label{sect:GalProj}
Let us use the quasimatrix of snapshots~\eqref{eq:ROM8} and  the inverse of the Cholesky factor $\bR$ of the mass matrix to define 
the new quasimatrix 
\begin{equation}
\bV(\bx) = \bU(\bx) \bR^{-1} = (\bv_0(\bx), \ldots, \bv_{n-1}(\bx)). 
\label{eq:defV}
\end{equation}
The $n m$ columns of this quasimatrix are organized in the blocks 
\begin{equation}
\bv_j(\bx) = \bU(\bx) \bR^{-1} \be_j, \qquad j = 0, \ldots, n-1,
\label{eq:defVj}
\end{equation}
and they form an orthonormal basis of the approximation space~\eqref{eq:apSpace}. This follows from definitions~\eqref{eq:ROM12},~\eqref{eq:defVj} and 
the Cholesky factorization~\eqref{eq:ROM23},
\begin{equation}
\bV^T \bV = \bR^{-T} \bU^T \bU \bR^{-1} = \bI_{nm},
\label{eq:defVort}
\end{equation}
where $\bI_{nm}$ is the $nm \times nm$ identity matrix and we used the linear algebra notation~\eqref{eq:ROM13p}. 
Therefore, the block columns \eqref{eq:defVj} of $\bV(\bx)$  are called the orthonormal snapshots.

We now see that the ROM propagator  is the projection of  $\cPq$ on the approximation space~\eqref{eq:F5}, 
written in the basis of the orthonormal snapshots,
\begin{equation}
\cPqR \stackrel{\eqref{eq:ROMProp}}{=}  \bR^{-T} \bS \bR^{-1}\,  \stackrel{\eqref{eq:ROM13}}{ = } \bR^{-T} \bU^T \cPq \bU \bR^{-1} \, \stackrel{\eqref{eq:defV}}{=} \bV^T \cPq \bV.
\label{eq:ROM21}
\end{equation}
Moreover, 
the initial ROM snapshot~\eqref{eq:ROMb} is the projection of the initial wave,
\begin{equation}
\bbR \stackrel{\eqref{eq:ROMb}}{=}\bR \be_0 = \bV^T \bU \be_0 = \bV^T \bb,
\label{eq:ROM22}
\end{equation}
because $\bU(\bx) \be_0 = \bb(\bx)$ and equations~\eqref{eq:defV} and~\eqref{eq:defVort} give
$
\bR = \bV^T \bU. 
$

\vspace{0.1in}
\begin{rem}
\label{rem.1}
Definition~\eqref{eq:defV}  of $\bV(\bx)$ is the Gram-Schmidt orthogonalization of the snapshots.  It is a causal
construction of the orthonormal basis
\begin{equation}
\bv_{j}(\bx) \in \mbox{span}\{\bu_{0}(\bx), \ldots, \bu_{j}(\bx)\}, \quad \forall ~ j = 0, \ldots, n-1,
\label{eq:ROM19}
\end{equation}
that respects the physics of the wave propagation, captured by the ROM time stepping scheme 
(\ref{eq:ROM28}--\ref{eq:ROM30}).  Indeed,~\eqref{eq:ROM22} and the block upper triangular structure of $\bR$ give that only the first  block of $\bbR$ is non-zero. This corresponds to the initial wave $\bb(\bx)$ being supported near the array. The wavefront of $\bu(t,\bx)$ penetrates deeper inside the medium for later time, and this  is reflected in the algebraic structure of the ROM snapshots~\eqref{eq:ROM31},
where the rows of blocks are filled in sequentially, for each time step. In particular, we obtain 
from ~\eqref{eq:ROM11},~\eqref{eq:ROM31} and~\eqref{eq:defV} that 
\begin{equation}
(\bu^{\RM}_0, \ldots, \bu^{\RM}_{n-1})  = \bR = \bV^T \bU = \bV^T \big(\bu_0, \ldots, \bu_{n-1}\big).
\label{eq:ROM20}
\end{equation}
\end{rem}

The importance of this remark  will become clear in sections \ref{sect:Galerk.4} and \ref{sect:Orthog1}, 
where we show  how the causality  preserving definition of the ROM induces properties of $\cPqR$ that are useful for solving the inverse scattering problem.

\subsection{Data-driven ROM}
\label{sect:Galerk.3}

In inverse scattering we do not know the snapshots, so how can we use the definition (\ref{eq:ROMProp}--\ref{eq:ROMb}) 
of the ROM? We now explain that, in fact, the mass and stiffness matrices  can be calculated from the data~\eqref{eq:F9}.
Consequently, we can compute the ROM from  (\ref{eq:ROMProp}--\ref{eq:ROMb}) and we can also get the ROM snapshots 
via the time stepping scheme (\ref{eq:ROM28}--\ref{eq:ROM30}).  Furthermore, we can use equation~\eqref{eq:ROM31} 
to calculate the Galerkin coefficients in the approximation~\eqref{eq:ROM9}, without knowing the approximation space~\eqref{eq:F5}.


Here we recall the calculation of mass and stiffness matrices from the data introduced in \cite{DtB,druskin2016direct} 
that we repeat for the convenience of the reader. 
The blocks of the mass matrix are, by definition~\eqref{eq:ROM12} and equation~\eqref{eq:Sn6},
\begin{equation}
\bM_{i,j} = \lb \cT_i\big( \cPq \big) \bb,\cT_j\big( \cPq \big) \bb \rb = \lb \bb, \cT_i\big( \cPq \big) \cT_j\big( \cPq \big) \bb\rb, \label{eq:Mass1}
\end{equation}
where the last equality is because $\cPq$ is self-adjoint. The Chebyshev polynomials of an arbitrary argument $z$ have the multiplicative property 
\begin{equation}
\cT_i(z) \cT_j(z) = \frac{1}{2} \big[ \cT_{i+j}(z) + \cT_{|i-j|}(z) \big],
\label{eq:Mass2}
\end{equation}
so using this property in~\eqref{eq:Mass1} and recalling~\eqref{eq:dataT} we obtain 
 \begin{align}
\bM_{i,j} &= \frac{1}{2} \Big[ \lb \bb, \cT_{i+j}\big( \cPq \big) \bb\rb + \lb \bb, \cT_{|i-j|}\big( \cPq \big) \bb\rb \Big] \\&= \frac{1}{2} 
\Big( \bD_{i+j} + \bD_{|i-j|} \Big), \qquad i,j = 0, \ldots, n-1.
\label{eq:Mass3}
\end{align}

The calculation of $\bS$ is similar. Starting with its definition~\eqref{eq:ROM13} and using equation~\eqref{eq:Sn6}, we have 
\begin{equation}
\bS_{i,j} = \lb \cT_i\big( \cPq \big) \bb,\cPq \cT_j\big( \cPq \big) \bb \rb = \lb \bb, \cT_i\big( \cPq \big) \cPq \cT_j\big( \cPq \big) \bb \rb, \label{eq:Stiff1}
\end{equation}
where by the multiplicative property~\eqref{eq:Mass2},
\[
z \cT_j(z) = \cT_1(z) \cT_j(z) =\frac{1}{2} \big[ \cT_{j+1}(z) + \cT_{|j-1|}(z) \big].
\]
Substituting in ~\eqref{eq:Stiff1}, using ~\eqref{eq:Mass2} one more time, and recalling ~\eqref{eq:dataT}, we get 
\begin{align}
\bS_{i,j} = \frac{1}{4} \Big(\bD_{i+j+1} + \bD_{|i-j+1|} + \bD_{|i+j-1|} + \bD_{|i-j-1|}\Big), \quad i,j = 0, \ldots, n-1.
\label{eq:Stiff3}
\end{align}

\vspace{0.05in}
\subsection{Properties of the ROM}
\label{sect:Galerk.4}
We state here the theorems that describe the properties of the ROM.
To lighten the presentation, we write the proofs in the appendixes.

\subsubsection{Data fit}
We saw in the previous section how the ROM is calculated from the data. The next theorem, proved in appendix \ref{ap:DatFit},
states that the ROM fits exactly these data.

\vspace{0.05in}
\begin{theorem}
\label{thm.1}
The ROM snapshots~\eqref{eq:ROM31} can be written as Chebyshev polynomials of the ROM propagator, similar to 
\eqref{eq:Sn6},
\begin{equation}
\bu^{\RM}_j = \cT_j \big(\cPqR \big) \bbR, \qquad j \ge 0,
\label{eq:chebROM}
\end{equation}
and the ROM defined by (\ref{eq:ROMProp}--\ref{eq:ROMb}) satisfies the data fit relations (\ref{eq:ROMThm1}).
\end{theorem}

\vspace{0.05in} Recall from section \ref{sect:Galerk.2} that  the first $n$ snapshots are represented exactly 
in our Galerkin scheme. Therefore, it is not surprising that the ROM fits the data for the first $n$ time instants. 
The interesting point of the theorem is that the data are fit for the remaining $n$ time instants. 
Physically, this is because the measurements  at the array of duration  $(2n-1)\tau$ can only sense
the medium up to the depth traveled by waves in half the time, and all this information is contained in our 
approximation space~\eqref{eq:apSpace}. This can be seen from the following equation 
\begin{align*}
\bD_{n-1+j} &= \lb \bb, \cT_{n-1+j}(\cPq) \bb\rb \\&= 
2 \lb \bb, \cT_{n-1}(\cPq) \cT_j(\cPq) \bb\rb -  \lb \bb, \cT_{|n-1-j|}(\cPq) \bb\rb \\
&= 2  \lb \bu_{n-1}, \bu_j\rb -  \lb \bb, \bu_{|n-1-j|}\rb,
\end{align*}
obtained using the recursion relation of Chebyshev polynomials, equation \eqref{eq:Sn6} and the self-adjointness 
of $\cPq$. Indeed, if $j = 1, \ldots, n-1$, the right hand side can be calculated in terms of the waves 
$\{\bu_l(\bx)\}_{0 \le l \le n-1} $. In fact, this is the case even for $j = n$, as shown by a  more involved calculation 
given in appendix \ref{ap:DatFit}.

\subsubsection{ROM factorization}
Just as we did in section~\eqref{sect:Galerk.1}, we can subtract $2 \bu^{\RM}_j$ from equation~\eqref{eq:ROM28} and divide the 
result by $\tau^2$ to obtain the ROM equivalent of the second order time stepping scheme (\ref{eq:Sn7}--\ref{eq:Sn9}),
\begin{align}
\frac{\bu^{\RM}_{j+1}-2 \bu^{\RM}_j + \bu^{\RM}_{j-1}}{\tau^2} + \cLqR \cLqR^T \bu^{\RM}_j &= 0, \quad j  \ge 0,  \label{eq:ROM41} \\
\bu^{\RM}_0 &= \bbR, \label{eq:ROM42} \\
\bu^{\RM}_1 &= \bu^{\RM}_{-1}, \label{eq:ROM43}
\end{align}
with matrix $\cLqR$ defined by the block Cholesky factorization
\begin{equation}
\frac{2}{\tau^2}(\bI_{nm}-\cPqR) = \cLqR \cLqR^T.
\label{eq:ROM44}
\end{equation}
This is the ROM analogue of the factorization 
\begin{equation}
\frac{2}{\tau^2}(I-\cPq) = \cLq \cLq^T \approx L(q) L(q)^T,
\label{eq:approxPq}
\end{equation}
where the approximation is as discussed in section (\ref{sect:Galerk.1}).

The next theorem, proved in appendix \ref{sect:PfThm2},  gives that the ROM propagator  is a block tridiagonal invertible  matrix. 
We return to this point in section  \ref{sect:Orthog1}, where we explain that the block tridiagonal $\cLqR \cLqR^T$ can be viewed as a finite 
difference approximation of  the second order partial differential operator $L(q) L(q)^T$.  

\vspace{0.05in}
\begin{theorem}
\label{thm.2}
The ROM propagator $\cPqR$ is  symmetric, block tridiagonal and  the matrix $\bI_{nm}-\cPqR$ is invertible. Therefore,
the Cholesky factor $\cLqR$  is an invertible matrix with lower block bidiagonal structure.
\end{theorem}

\subsubsection{Galerkin-Petrov projection}
\label{sect:GP}
We now show  that $\cLqR$ is a Galerkin-Petrov projection of the operator $\cLq$ on the  subspace $\mathfrak{X}$ defined in 
\eqref{eq:F5} and the subspace
\begin{equation}
\hat{\mathfrak{X}}  = \mbox{span}\{ \bhu_{0}(\bx), \ldots, \bhu_{n-1}(\bx)\},
\label{eq:ApproxSpaceH}
\end{equation}
of the first $n$ dual snapshots denoted by the hat. 

The dual snapshots  are defined using the first order system formulation of 
the time stepping scheme (\ref{eq:Sn7}--\ref{eq:Sn9}),
\begin{align}
\frac{\bu_{j+1}(\bx)-\bu_j(\bx)}{\tau} &= -\cLq \bhu_j(\bx), \label{eq:ROM47}\\
\frac{\bhu_{j}(\bx)-\bhu_{j-1}(\bx)}{\tau} &= \cLq^T \bu_j(\bx), \quad j \ge 0, \label{eq:ROM48} \\
\bu_0(\bx) &= \bb(\bx) \label{eq:ROM49} \\
 \bhu_0(\bx) + \bhu_{-1}(\bx) &= 0. \label{eq:ROM50}
\end{align}
Indeed, it is easy to check that (\ref{eq:ROM47}--\ref{eq:ROM50}) implies (\ref{eq:Sn7}--\ref{eq:Sn8}) and the initial condition~\eqref{eq:Sn9}  follows from
\begin{align*}
\frac{\bu_{-1}(\bx)-\bu_1(\bx)}{\tau} &= -\frac{\bu_1(\bx) - \bu_0(\bx)}{\tau} - \frac{\bu_0(\bx)-\bu_{-1}(\bx)}{\tau} \\
& \stackrel{\eqref{eq:ROM47}}{=}\cLq( \bhu_0 (\bx) + \bhu_{-1}(\bx))  \stackrel{\eqref{eq:ROM50}}{=}0.
\end{align*}
The first dual snapshot is obtained from~\eqref{eq:ROM48} evaluated at $j = 0$ and~\eqref{eq:ROM50},
\begin{equation}
\bhu_0(\bx) = \hat{\bb}(\bx) = \frac{\tau}{2} \cLq^T \bb(\bx).
\label{eq:ROM51}
\end{equation}
The half time step in this equation shows that (\ref{eq:ROM47}--\ref{eq:ROM50}) is a leap-frog scheme, where the dual wave is evaluated at the time instants $(j+1/2)\tau$, and the primary wave is evaluated at the time instants $j \tau$, for $j \ge 0$. 

Similarly, we can use the first order system formulation of  the ROM time stepping scheme (\ref{eq:ROM41}--\ref{eq:ROM43}) to define the dual ROM snapshots $\bhuR_j$, 
\begin{align}
\frac{\bu^{\RM}_{j+1}-\bu^{\RM}_j}{\tau} &= -\cLqR \bhuR_j, \label{eq:ROM52}\\
\frac{\bhuR_{j}-\bhuR_{j-1}}{\tau} &= \cLqR^T \bu^{\RM}_j, \quad j \ge 0, \label{eq:ROM53} \\
\bu^{\RM}_0 &= \bbR \label{eq:ROM54} \\
\bhuR_0+ \bhuR_{-1} &= 0,\label{eq:ROM55}
\end{align}
and obtain as above that
\begin{equation}
\bhuR_0 = \hat{\bb}^{\RM} =  \frac{\tau}{2} \cLqR^T \bbR. 
\label{eq:ROM56}
\end{equation}

The orthogonalization of the dual snapshots and their use in the Galerkin-Petrov projection of $\cLq$ are in the next theorem,  
proved in appendix \ref{ap:PfThm5}.
\vspace{0.05in}
\begin{theorem}
\label{thm:5}
Denote by
$
\hat \bU(\bx) = \big(\bhu_0(\bx), \ldots, \bhu_{n-1}\big)$ the quasimatrix of the first $n$ dual snapshots, which span the space $\hat{\mathfrak{X}}$ defined in~\eqref{eq:ApproxSpaceH}. The following statements hold:

\vspace{0.05in}\begin{enumerate}
\item[(i)] 
There exists an orthonormal basis of the space $\hat{\mathfrak{X}}$, the columns of  the quasimatrix 
$\hat \bV(\bx) = \big(\hat \bv_0(\bx), \ldots, \hat \bv_{n-1}(\bx) \big),$
satisfying 
\begin{equation}
\hat \bU(\bx) = \hat \bV(\bx) \hat \bR,
\label{eq:BHU2}
\end{equation}
where $\hat \bR$ is the matrix of the first $n$ ROM dual snapshots
\begin{equation}
\hat \bR = \big(\bhuR_0,\ldots, \bhuR_{n-1}\big).
\label{eq:BHU3}
\end{equation}
This is the analogue of equations ~\eqref{eq:defV},~\eqref{eq:ROM20},  and $\hat \bR$ is block upper triangular. 
\item[(ii)] The matrix $\cLqR$ defined in the Cholesky factorization~\eqref{eq:ROM44} is the Galerkin-Petrov projection of the operator $\cLq$
defined in~\eqref{eq:Sn10}, on the spaces $\mathfrak{X}$ and $\hat{\mathfrak{X}}$, 
\begin{equation}
\cLqR = \bV^T \cLq \hat \bV.
\label{eq:ROM46}
\end{equation}
\end{enumerate}
\end{theorem}

\vspace{0.05in} Similar to Remark \ref{rem.1}, we note that definition~\eqref{eq:BHU2} of $\hat\bV(\bx)$ is the Gram-Schmidt 
orthogonalization of the dual snapshots, which gives the causal ROM dual snapshots gathered in  the block 
upper triangular matrix $\hat \bR$. 

\subsection{Orthogonal transformations}
\label{sect:Orthog}
The block Cholesky factorization~\eqref{eq:ROM23} of the mass matrix is defined up to an orthogonal transformation of the form
\begin{equation}
\bY = \mbox{diag} \left( \bY_0, \ldots, \bY_{n-1} \right),
\label{eq:Y1}
\end{equation}
 with $m \times m$ orthogonal matrices $\bY_j$ , for $j = 0, \ldots, n-1$. That is to say, the matrix
$ \bY^T \bR$ 
is also $nm \times nm$ block upper triangular, and satisfies
\begin{equation}
(\bY^T \bR)^T \bY^T \bR = \bR^T \bR = \bM.
\label{eq:Y3}
\end{equation}
Moreover,  if we replace $\bR$ with $\bY^T \bR$ in definitions (\ref{eq:ROMProp}--\ref{eq:ROMb}), we get the ROM propagator
\begin{equation}
\PY= \bY^T \cPqR \bY,
\label{eq:Y4}
\end{equation}
which  has all the properties described in section \ref{sect:Galerk.4}.  The quasimatrix of orthonormal snapshots is transformed to 
\begin{equation}
\VY(\bx) = \left(\bvc_0(\bx), \ldots, \bvc_{n-1}(\bx) \right) = \bU(\bx) (\bY^T \bR)^{-1} = \bV(\bx) \bY,
\label{eq:Y9}
\end{equation}
and for the dual snapshots we have, similarly, 
\begin{equation}
\VYd(\bx) = \left(\hat \bvc_0(\bx), \ldots, \hat \bvc_{n-1}(\bx) \right) = \hat \bV(\bx) \hat \bY.
\label{eq:Y10}
\end{equation}
The Galerkin-Petrov projection~\eqref{eq:ROM46} becomes 
\begin{equation}
\boldsymbol{\Lambda}^{\RM}(q) = \VY^T \cLq \VYd = \bY^T \cLqR \hat \bY.
\label{eq:Y11}
\end{equation}
Here $\hat \bY$ is another arbitrary orthogonal transformation of the form~\eqref{eq:Y1}, and $\boldsymbol{\Lambda}^{\RM}(q)$ is block lower bidiagonal.
In equations (\ref{eq:Y4}--\ref{eq:Y11}) we use greek letters for the transformed ROM and the orthonormal snapshots. These depend 
on  $\bY$ and $\hat \bY$, but we suppress this dependence in the notation.  

\subsection{Connection to finite differences}
\label{sect:Orthog1}
To interpret the ROM matrix  $\cLqR$ as an approximate  finite difference scheme for the operator $L(q)$, we write it here explicitly 
using the following assumption: 

\vspace{0.05in}\begin{ass}
\label{as.2}
 The iteration 
\begin{align}
\left[ \bphi_{j+1}(\bx) - \bphi_j(\bx)\right] \bga_j^{-1} &= - \cLq \bhphi_j(\bx), \label{eq:St3} \\
\left[ \bhphi_{j}(\bx) - \bhphi_{j-1}(\bx)\right] \bhga_j^{-1} &=  \cLq^T \bphi_j(\bx),  \qquad j \ge 0, \label{eq:St4}
\end{align}
with initial conditions 
\begin{equation}
\bphi_0(\bx) = \bb(x), \qquad \bhphi_{-1}(\bx) = {\bf 0},
\label{eq:St5}
\end{equation} 
and with $m\times m$ symmetric matrix coefficients 
\begin{equation}
\bhga_j = \lb \bphi_j, \bphi_j \rb^{-1}, \qquad 
\bga_j = \lb \bhphi_j, \bhphi_j \rb^{-1}, 
\label{eq:St2}
\end{equation} 
does not break down for $j = 0, \ldots, n-1$.  That is to say, the columns in each $\bphi_j(\bx)$ and $\bhphi_j(\bx)$ remain linearly independent, so the matrices~\eqref{eq:St2} are defined.\end{ass}

\vspace{0.05in}
We explain in appendix \ref{ap:A} that this assumption is basically the same as saying that the 
first-order block Lanczos procedure \eqref{eq:St3}--\eqref{eq:St2} for calculating 
orthogonal bases of the spaces $ \mathfrak{X} $ and $\hat{\mathfrak{X}}$ does not break down. 
If this is the case, we have the following result, proved in appendix \ref{ap:A}.

\vspace{0.05in}
\begin{theorem}
\label{thm.7}
Under the Assumption \ref{as.2}, there exists a choice of the square roots of the coefficients \eqref{eq:St2},
\begin{equation}
\bga_j = \bGa_j \bGa_j^T, \qquad \bhga_j = \bhGa_j \bhGa_j^T, \qquad j \ge 0,
\label{eq:St1p}
\end{equation}
which relates the orthonormal snapshots defined in \eqref{eq:defV} and \eqref{eq:BHU2} to the 
solution of the iteration (\ref{eq:St3}--\ref{eq:St5}) as follows,
\begin{equation}
\bv_j(\bx) = \bphi_j(\bx) \bhGa_j, \qquad \hat \bv_j(\bx) = \bhphi_j(\bx) \bGa_j, \qquad j = 0,\ldots, n-1.
\label{eq:St1}
\end{equation} 
Moreover, the block entries of the ROM matrix 
$\cLqR$ are defined by these square roots as 
\begin{align}
\boldsymbol{\cL}^{\RM}_{j,j}(q) &= \bhGa_j^{-1} \bGa_j^{-T}, \qquad \quad j = 0, \ldots, n-1,\label{eq:St6} \\
\boldsymbol{\cL}^{\RM}_{j+1,j}(q) &=- \bhGa_{j+1}^{-1} \bGa_j^{-T}, \qquad j = 0, \ldots, n-2 .\label{eq:St7}
\end{align}
\end{theorem}

\vspace{0.05in} 
Recall from section \ref{sect:Galerk.1} that $\cLq$ is an approximation of the first order partial differential operator $L(q)$. 
Equation~\eqref{eq:St3} shows that this operator is captured by the ROM as a  finite difference scheme, 
where each step corresponds to a time instant indexed by $j\ge 0$. The "steps" are $m \times m$ 
matrices, due to the fact that there are $m$ source excitations. As the time index  increases, the iteration (\ref{eq:St3}--\ref{eq:St5})
and definition~\eqref{eq:St1} generate orthonormal snapshots that satisfy the causality relations \eqref{eq:ROM19} and 
\eqref{eq:BHU2}.  Initially, these snapshots are in the range of $\bu_0(\bx) = \bb(\bx)$ and 
$\bhu_0(\bx) = \hat \bb(\bx)$, respectively, and are supported near the array.
At the next time instant the wave front advances  to a depth of the order $c_o \tau$, so $\bu_1(\bx)$ will have large entries around depth $c_o \tau$. Due to causality and orthogonality, $\bv_1(\bx)$ should peak around this depth. The same holds for the dual orthonormal snapshots. 
Arguing this way, we expect that the peak values of the orthonormal snapshots follow the progression of the wave front inside the medium.
This is confirmed by numerical simulations in \cite{DtB,borcea2019robust,druskin2018nonlinear}.

\vspace{0.05in}
\begin{rem}
\label{rem.3}
Every ROM matrix $\boldsymbol{\Lambda}^{\RM}(q)$ related to $\cLqR$ by \eqref{eq:Y11} has the finite differences interpretation
(\ref{eq:St3}--\ref{eq:St5}), (\ref{eq:St6}--\ref{eq:St7}), and it is the Galerkin-Petrov projection of the operator $\cLq$ in the orthonormal bases 
\[
\bvc_j(\bx) = \bv_j(\bx) \bY_j = \bphi_j(\bx) \bhGa_j \bY_j, \qquad 
\hat \bvc_j(\bx) =\hat  \bv_j(\bx) \hat \bY_j = \bhphi_j(\bx) \bGa_j \hat \bY_j, \qquad j \ge 0.
\]
This non-uniqueness  is due to the multiple choices of the square roots of $\bga_j$ and $\bhga_j$, 
\begin{align*}
\bga_j &= \bGa_j \bGa_j^T = \big(\bGa_j \hat \bY_j\big)\big(\bGa_j \hat \bY_j\big)^T, \\
\bhga_j &= \bhGa_j \bhGa_j^T = \big(\bhGa_j  \bY_j\big)\big(\bhGa_j \bY_j\big)^T, \qquad j \ge 0.
\end{align*}
The block entries of  $\boldsymbol{\Lambda}^{\RM}(q)$ are given by 
\begin{align*}
\boldsymbol{\Lambda}^{\RM}_{j,j}(q) &= \bY_j^{-1} \bhGa_j^{-1} \bGa_j^{-T} \hat \bY_j^{-T} = 
\bY_j^T \boldsymbol{\cL}^{\RM}_{j,j}(q) \hat \bY_j, \qquad \qquad  j = 0, \ldots, n-1, \\
\boldsymbol{\Lambda}^{\RM}_{j+1,j}(q) &= - \bY_{j+1}^{-1} \bhGa_{j+1}^{-1} \bGa_j^{-T} \hat \bY_j^{-T} =
\bY_{j+1}^T \boldsymbol{\cL}^{\RM}_{j+1,j}(q) \hat \bY_j , \quad \; j = 0, \ldots, n-2,
\end{align*}
which is precisely \eqref{eq:Y11} written block-wise.
\end{rem}

\vspace{0.05in}
Assumption \ref{as.2} and Theorem \ref{thm.7} are written as if we knew the operator $\cLq$, which is not the case in the inverse problem. 
Their purpose is to interpret the ROM as a finite difference scheme, which we use in the next section to motivate the new inversion method. However, matrices $\{\bga_j, \bhga_j\}_{0 \le j \le n-1}$ can be determined from the data, from the 
equation \begin{align*}
\bY^T \frac{2}{\tau^2} \Big( \bI_{nm} - \cPqR\Big) \bY= \boldsymbol{\Lambda}^{\RM}(q)\boldsymbol{\Lambda}^{\RM}(q)^T.
\end{align*}
Substituting the expression of $\boldsymbol{\Lambda}^{\RM}(q)$ described in Remark \ref{rem.3}, for a given convention of
the matrix square root,  and  equating block-wise, one obtains  an iteration which defines sequentially $\bga_j$, $\bhga_j$, starting with $\bhga_0 = \lb \bb, \bb \rb^{-1}$, and also the diagonal blocks of $\bY$.

\vspace{0.05in}
\section{Inverse scattering}
\label{sect:Inv}
We now use the ROM for solving the inverse problem.  The proposed method generalizes to all  linear waves in isotropic media, in the backscattering setup described in the introduction. Nevertheless,  to make the presentation explicit, we focus attention on inverse scattering for sound waves.

We refer to \cite{DtB,borcea2019robust} for the derivation of the hyperbolic problem (\ref{eq:F1}--\ref{eq:F3}) from the acoustic wave equation,
where $\bu(t, \bx)$ is related to the pressure and the operators $L(q)$ and $L(q)^T$ are given by
\begin{equation}
\begin{array}{rcl}
L(q) \hat\bu(t, \bx) & = & L(0) \hat\bu(t, \bx) + \dfrac{1}{2} [c(\bx) \nabla q(\bx)] \cdot \hat\bu(t, \bx),\\
L(q)^T \bu(t, \bx) & = & L(0)^T \bu(t, \bx) + \dfrac{1}{2} [c(\bx) \nabla q(\bx)] \bu(t, \bx), 
\end{array}
\label{eq:In1}
\end{equation}
with
\begin{equation}
\begin{array}{rcl}
L(0) \hat\bu(t, \bx) & = & - \sqrt{c(\bx)} \, \mbox{div}\, \big[ \sqrt{c(\bx)} \hat\bu(t, \bx) \big],\\
L(0)^T \bu(t, \bx) & = & \sqrt{c(\bx)} \, \nabla \big[ \sqrt{c(\bx)} \bu(t, \bx) \big].
\end{array}
\label{eq:In10}
\end{equation}
The dot in \eqref{eq:In1} denotes the inner product in $\mathbb{R}^d$ and the vector-valued 
function $\hat \bu(t, \bx)$ is related to the acoustic velocity as in \cite{DtB,borcea2019robust}.

Here $c(\bx)$ is the assumed smooth wave speed, the known kinematic model,  and the unknown reflectivity is defined by 
\begin{equation}
q(\bx) = \ln \sqrt{\sigma(\bx)},
\label{eq:In1p}
\end{equation}
in terms of the   acoustic impedance $\sigma(\bx)$.  As explained in the introduction, the model (\ref{eq:In1}--\ref{eq:In1p}) arises 
when separating the estimations of the kinematic model (assumed known here) and the reflectivity. It applies to the backscattering 
setup, where as shown  in \cite{BeylkinBurridge}, the reflections recorded at a small array are due mainly to relative variations of the impedance.

The main idea of our inversion method is that instead of using the conventional nonlinear least squares data fit 
minimization formulation, it is better to minimize the difference of the ROM matrices $\cLqR-\cLoR$, 
where $\cLoR$ is defined as in~\eqref{eq:ROM44}, but for the reference medium with zero reflectivity. 
We motivate this optimization formulation in section \ref{sect:Inv.1}, with a discussion based on the results 
in sections \ref{sect:Galerk.4}--\ref{sect:Orthog1}. The inversion algorithm is described in  section \ref{sect:Inv.2}.

\subsection{ROM parametrization of the reflectivity}
\label{sect:Inv.1}
There are two ways of understanding how the ROM encodes information about the unknown reflectivity~\eqref{eq:In1p}. 
The first is based on the finite difference interpretation described in 
Theorem \ref{thm.7}.  
The second is based on the Gram-Schmidt orthogonalization of the snapshots.

\subsubsection{Finite differences interpretation}
\label{sect:Inv.1.1}
We see from definition~\eqref{eq:In1} that the operator $L(q)-L(0)$ depends linearly on 
$q(\bx).$ Here we explain why the ROM version of this operator, the matrix  $\cLqR - \cLoR$,
is expected to inherit approximately this linear dependence.

Let us use Theorem \ref{thm.7} for the reference medium with zero reflectivity.  We obtain the analogue of 
the finite difference scheme (\ref{eq:St3}--\ref{eq:St5}),
\begin{align}
\big[ \bphi_{j+1}^{(0)}(\bx) - \bphi_j^{(0)}(\bx) \big] {\itbf h}_j^{-1} &= - \cLo \bhphi_j^{(0)}(\bx), \label{eq:St3o}\\
\big[ \bhphi_{j}^{(0)}(\bx) - \bhphi_{j-1}^{(0)}(\bx) \big] \hat {\itbf h}_j^{-1} &= \cLo^T \bphi_j^{(0)}(\bx),  \qquad j \ge 0, \\
\bphi_0^{(0)}(\bx) &= \bphi_0(\bx) = \bb(\bx), \\\hat \bphi_{-1}^{(0)}(\bx) &=  \hat \bphi_{-1}(\bx) = {\bf 0}, \label{eq:St5o}
\end{align}
where the superscript $(0)$ indicates that the reflectivity is zero,  and the steps
$\{{\itbf h}_j, \hat {\itbf h}_j\}_{j \ge 0}$ are the analogues of  \eqref{eq:St2},
\begin{equation}
\hat {\itbf h}_j = \lb \bphi_j^{(0)}, \bphi_j^{(0)} \rb^{-1}, \qquad 
{\itbf h}_j = \lb \bhphi_j^{(0)}, \bhphi_j^{(0)} \rb^{-1}, \qquad j \ge 0.
\label{eq:stepso}
\end{equation}
The square roots of these steps 
\begin{equation}
{\itbf h}_j = \bH_j \bH_j^T, \qquad \hat{\itbf h}_j = \hat \bH_j \hat \bH_j^T, \qquad j \ge 0,
\label{eq:SQRTo}
\end{equation}
define the block lower bidiagonal ROM matrix $\cLoR$, with entries given by the analogues
of (\ref{eq:St6}--\ref{eq:St7}),
 \begin{align}
\boldsymbol{\cL}^{\RM}_{j,j}(0) &= \hat \bH_j^{-1} \bH_j^{-T}, \qquad \quad j = 0, \ldots, n-1,\label{eq:St6o} \\
\boldsymbol{\cL}^{\RM}_{j+1,j}(0) &=- \hat \bH_{j+1}^{-1} \bH_j^{-T}, \qquad j = 0, \ldots, n-2.\label{eq:St7o}
\end{align}

We now have two exact finite differences schemes: For the operator $\cLqR$, as given in 
(\ref{eq:St3}--\ref{eq:St5}), and  for the operator $\cLoR$, as given in  (\ref{eq:St3o}--\ref{eq:St5o}). To compare the two,
let us use the transformation 
\begin{align}
\btphi_j(\bx) &= \bv_j(\bx) \hat \bH_j^{-1} = \bphi_j(\bx) \bhGa_j \hat \bH_j^{-1}, \label{eq:PhiT}\\
\bhtphi_j(\bx) &= \hat \bv_j(\bx)  \bH_j^{-1} = \bhphi_j(\bx) \bGa_j  \bH_j^{-1},  \qquad j \ge 0, \label{eq:hPhiT}
\end{align}
so that $\btphi_j(\bx)$ and $\hat \btphi_j(\bx)$ are normalized as in \eqref{eq:stepso}
\[
\lb \btphi_j,\btphi_j \rb = \lb \bphi^{(0)}_j,\bphi_j^{(0)} \rb = \hat {\itbf h}_j^{-1}, \qquad 
\lb \hat \btphi_j, \hat \btphi_j \rb = \lb \hat \bphi^{(0)}_j, \hat \bphi_j^{(0)} \rb = {\itbf h}_j^{-1}, \qquad j \ge 0.
\]
We also introduce the matrices
\begin{equation}
\bsigma_j^{\frac{1}{2}} = \bhGa_j \hat \bH_j^{-1}, \qquad 
\hat \bsigma_j^{\frac{1}{2}} = \bH_j\bGa_j^{-1}, \qquad j \ge 0,
\label{eq:defSigs}
\end{equation}
which will be interpreted below as approximations of the square root of the impedance. Substituting (\ref{eq:PhiT}--\ref{eq:hPhiT}) 
in (\ref{eq:St3}) and using definition \eqref{eq:defSigs}, we obtain 
\begin{align}
\cLq \hat \btphi_j(\bx) + \big[ \btphi_{j+1}(\bx) - \btphi_{j}(\bx)\big] {\itbf h}_j^{-1} &
=  \btphi_{j+1}(\bx) \bQ_{j}^++ \btphi_{j}(\bx)\bQ_j^-, \label{eq:ScH1}
\end{align}
where the matrices 
\begin{align}
\bQ_{j}^+ = \bsigma_{j+1}^{-\frac{1}{2}} \big[ \bsigma_{j+1}^{\frac{1}{2}} - 
\big(\hat \bsigma_j^{\frac{1}{2}}\big)^T\big] {\itbf h}_j^{-1}, \qquad 
\bQ_{j}^- = \bsigma_{j}^{-\frac{1}{2}} \big[ \big(\hat \bsigma_j^{\frac{1}{2}}\big)^T - 
\bsigma_j^{\frac{1}{2}}  \big] {\itbf h}_j^{-1}, \label{eq:Qs}
\end{align}
satisfy 
\begin{align}
\bsigma_{j+1}^{\frac{1}{2}} \bQ_j^+ + \bsigma_j^{\frac{1}{2}} \bQ_j^- &= \big[ \bsigma_{j+1}^{\frac{1}{2}} - \bsigma_j^{\frac{1}{2}} \big] {\itbf h}_j^{-1}.
\label{eq:add}
\end{align}
The second term in the left hand side in \eqref{eq:ScH1} looks like the finite differences approximation of 
$-\cLo$ in equation \eqref{eq:St3o}, although there the operator acts on a different space, spanned by the 
snapshots in the reference medium. The right hand side in \eqref{eq:ScH1} looks like a finite difference 
approximation of the operator 
\[
L(q)-L(0) \stackrel{\eqref{eq:In1}}{=}  c(\bx) \nabla q(\bx) \cdot , \qquad q(\bx) = \ln \sqrt{\sigma(\bx)}.
\]
The discretization corresponds to time stepping, so we can view ${\itbf h}_j$ and $\hat  {\itbf h}_j$ as primary and dual 
grid steps for discretization in range. These steps depend on the kinematic model $c(\bx)$ which is the same in the 
reference and the unknown medium. The matrices $\bsigma_j^{\frac{1}{2}}$ and $\hat \bsigma_j^{\frac{1}{2}}$ 
can be viewed as approximations of $\sqrt{\sigma(x)}$ on the primary grid and dual grid, respectively, 
and the matrices \eqref{eq:Qs} can be viewed as approximations of $c(\bx) \nabla q(\bx) \cdot$, 
up to some factors which add up as in \eqref{eq:add}.

It remains to study the difference of the block lower bidiagonal ROM matrices $\cLqR$ and $\cLoR$,  using the expressions (\ref{eq:St6}--\ref{eq:St7}) and (\ref{eq:St6o}--\ref{eq:St7o}) of their  entries. We obtain that 
\begin{align}
\boldsymbol{\cL}_{j,j}^{\RM}(q) - \boldsymbol{\cL}_{j,j}^{\RM}(0) &= \bhGa_j^{-1} \bGa_j^{-T} - \hat \bH_j^{-1} \bH_j^{-T} 
\stackrel{\eqref{eq:defSigs}}{=} \hat \bH_j^{-1} \bsigma_j^{-\frac{1}{2}} \big[ (\hat \bsigma_j^{\frac{1}{2}} \big)^T - 
\bsigma_j^{\frac{1}{2}} \big] \bH_j^{-T} \nonumber \\
&\hspace{-0.06in}\stackrel{\eqref{eq:Qs}}{=} \hat \bH_j^{-1} \bQ_j^{-} \bH_j, \qquad j = 0, \ldots, n-1,
\end{align}
and 
\begin{align}
\boldsymbol{\cL}_{j+1,j}^{\RM}(q) - \boldsymbol{\cL}_{j+1,j}^{\RM}(0) &
= \hat \bH_{j+1}^{-1} \bH_j^{-T} -\bhGa_{j+1}^{-1} \bGa_j^{-T} 
\stackrel{\eqref{eq:defSigs}}{=} \hat \bH_{j+1}^{-1} \bsigma_{j+1}^{-\frac{1}{2}} 
\big[ \bsigma_{j+1}^{\frac{1}{2}}  - \big(\hat \bsigma_j^{\frac{1}{2}}\big)^T \big] \bH_j^{-T} \nonumber \\
&\hspace{-0.06in} \stackrel{\eqref{eq:Qs}}{=} \hat \bH_{j+1}^{-1} \bQ_j^{+} \bH_j, \qquad j = 0, \ldots, n-2.
\end{align}
Therefore, $\cLqR-\cLoR$ is linear in the matrices $\bQ_j^\pm$ defined in \eqref{eq:Qs}, which are expected to 
approximate the gradient of the reflectivity, as explained above. 

We remark that the approximate linear dependence  of $L(q)-L(0)$ and therefore of $\cLqR-\cLoR$ on the 
gradient of the reflectivity is important in inversion, as it leads to an emphasis of the boundaries of reflectors 
and to sharp estimates of their support, as observed in the numerical results in section~\ref{sect:Numer}.

\subsubsection{Gram-Schmidt orthogonalization interpretation}
\label{sect:Inv.1.2}

Recall Remark \ref{rem.1} on the causal construction of the orthonormal snapshots, via the Gram-Schmidt orthogonalization 
\eqref{eq:defV}, and the similar result  in section \ref{sect:GP} for the orthonormal dual snapshots. We now explain that this
construction leads to projection matrices $\bV(\bx)$ and $\hat \bV(\bx)$ should be nearly independent of the unknown reflectivity. 
In light of Theorem \ref{thm:5}, this implies that $\cLqR$ has approximately the same affine dependence on $q(\bx)$ as 
$L(q) \approx \cLq$, and gives another motivation for  the inversion based on $\cLqR-\cLoR$.  

For simplicity of the argument, we assume in this section only that the kinematic model is constant $c(\bx) = c_o$. 
The extension to arbitrary $c(x)$ is straightforward in one dimension, where we can use the travel time transformation 
$x \mapsto \int_0^x ds \, c^{-1}(s)$ to eliminate the wave speed from the wave equation. In higher dimensions the 
extension is not easy and may not even be true, unless the medium is nice enough, so that the wave progresses forward 
at each time step and there are no lensing effects, as we have assumed so far.

\vspace{0.05in} 
\textbf{One dimension:}  We begin with the case $d = 1$, where the domain $\Omega$ is an interval. This is easier to 
understand because there is only one sensor $(m = 1)$ and there are no block linear algebra calculations. 

The Gram-Schmidt orthogonalization~\eqref{eq:defV} is 
\begin{equation}
\Big(u_0(x), \ldots, u_{n-1}(x)\Big) = \Big(v_0(x), \ldots, v_{n-1}(x)\Big) \bR, \qquad x \in \Omega,
\label{eq:GS1}
\end{equation}
where now $\bR$ is $n \times n$ upper triangular and we do not use bold symbols because the snapshots and $x$ are real valued. 
Let us evaluate this equation at  the locations
$x_j = c_o j \tau$ of the wavefront at the first $n$ time instants $j \tau$ of the measurements,  for $j = 0, \ldots, n-1$, and gather the results 
in the linear system
\begin{equation}
\begin{pmatrix} 
u_0(x_0) & \ldots & u_{n-1}(x_0) \\
u_0(x_1) & \ldots & u_{n-1}(x_1) \\
\vdots &  & \vdots \\
u_0(x_{n-1})  &\ldots & u_{n-1}(x_{n-1})
\end{pmatrix} = \begin{pmatrix} 
v_0(x_0) & \ldots & v_{n-1}(x_0) \\
v_0(x_1) & \ldots & v_{n-1}(x_1) \\
\vdots &  & \vdots \\
v_0(x_{n-1}) & \ldots & v_{n-1}(x_{n-1})
\end{pmatrix} \bR.
\label{eq:GS3}
\end{equation}
The first factor in the right hand side is a nearly orthogonal matrix, because 
\begin{equation}
\sum_{l = 0}^{n-1} v_i(x_l) v_j(x_l) \approx \frac{1}{c_o \tau} \int_{\Omega} dx \, v_i(x) v_j(x)  = \frac{1}{c_o \tau} \delta_{ij}, 
\qquad i,j = 0, \ldots, n-1.
\label{eq:GS4}
\end{equation}
Here the integral is approximated by a Riemann sum and the  integrand  is supported  in the interval 
$
[0,\min\{x_i,x_j\}] \subseteq [0,x_{n-1}] \subset \Omega,
$ by the causality relation~\eqref{eq:ROM19}.  

We conclude that \eqref{eq:GS3} is basically a  $QR$ factorization \cite[Section 5.2]{golubVanLoan}, 
which seeks an orthonormal basis that transforms the left hand side  to upper triangular form.  But the left hand side  is already upper triangular by construction
\[
u_j(x) = 0 ~~\mbox{for}~~ x > x_j,  \qquad j = 0, \ldots, n-1,
\]
so there is no transformation to be made, and 
\begin{equation}
\begin{pmatrix} 
v_0(x_0) & \ldots & v_{n-1}(x_0) \\
v_0(x_1) & \ldots & v_{n-1}(x_1) \\
\vdots &  & \vdots \\
v_0(x_{n-1}) & \ldots & v_{n-1}(x_{n-1})
\end{pmatrix} \approx \frac{1}{\sqrt{c_o \tau}} \bI_n,
\label{eq:GS5}
\end{equation}
up to $\pm$ sign ambiguity on the diagonal. This matrix has exactly zero entries below the diagonal, by the causality of the orthonormal snapshots, so the approximation  applies only to the upper triangular part.
The quasimatrix $\bV(x)$ of the orthonormal snapshots is an interpolation of the entries in 
\eqref{eq:GS5}, so it is  approximately independent of $q(x)$. 
The same argument applies to the quasimatrix $\hat \bV(x)$ of orthonormal dual snapshots. 

\vspace{0.05in}
\begin{rem}
\label{rem.5}
We expect from (\ref{eq:GS4}--\ref{eq:GS5}) that the approximation of the quasimatrix $\bV(x)$ and its dual analogue 
$\hat \bV(x)$ by a multiple of the identity improves when we decrease the time sampling interval $\tau$. This is the case 
up to a point, because if $\tau$ is too small, then the snapshots become linearly dependent (up to machine precision) 
and the mass matrix~\eqref{eq:ROM12} is no longer invertible.  A good strategy for choosing $\tau$ is according to the 
Nyquist criterion which takes into consideration the temporal period of oscillation of the wave.
\end{rem}

\vspace{0.05in} 
\textbf{Higher dimensions:} Here we have the block Gram-Schmidt orthogonalization
\begin{equation}
\Big(\bu_0(\bx), \ldots, \bu_{n-1}(\bx)\Big) = \Big(\bv_0(x), \ldots, \bv_{n-1}(\bx)\Big) \bR, \qquad \bx \in \Omega,
\label{eq:GS6}
\end{equation}
where $\bR$ is $nm \times nm$ block upper triangular, with $m \times m$ blocks. 

To write equation \eqref{eq:GS6}  as a block QR factorization, 
the analogue of \eqref{eq:GS3},  consider the system of coordinates 
$\bx = (x,\bx^\perp) \in \Omega,$
with origin at the center of the array, where $x \in [0,x_{\max}]$ is the depth (range) coordinate orthogonal to the array and $\bx^\perp$ is 
the cross-range in the plane of the array.  Then, we can evaluate~\eqref{eq:GS6} at points 
\[
(x_j, \bx^\perp_s), \qquad x_j = c_o j \tau, \qquad j = 0, \ldots, n-1, ~~ s = 1, \ldots, m,
\]
for some appropriate $\{\bx_s^\perp\}_{1 \le s \le m}$. Using the block notation 
\begin{equation}
\uu_j(x) = \begin{pmatrix} \bu_j(x,\bx^\perp_1) \\
\vdots \\
\bu_j(x,\bx^\perp_m) 
\end{pmatrix} \in \mathbb{R}^{m\times m}, \qquad 
\uv_j(x) = \begin{pmatrix} \bv_j(x,\bx^\perp_1) \\
\vdots \\
\bv_j(x,\bx^\perp_m) 
\end{pmatrix} \in \mathbb{R}^{m\times m}, \qquad 
\end{equation}
we get 
\begin{equation}
\begin{pmatrix} 
\uu_0(x_0) & \ldots & \uu_{n-1}(x_0) \\
\uu_0(x_1) & \ldots & \uu_{n-1}(x_1) \\
\vdots &  & \vdots \\
\uu_0(x_{n-1})  &\ldots & \uu_{n-1}(x_{n-1})
\end{pmatrix} = \begin{pmatrix} 
\uv_0(x_0) & \ldots & \uv_{n-1}(x_0) \\
\uv_0(x_1) & \ldots & \uv_{n-1}(x_1) \\
\vdots &  & \vdots \\
\uv_0(x_{n-1}) & \ldots & \uv_{n-1}(x_{n-1})
\end{pmatrix} \bR.
\label{eq:GS8}
\end{equation}
Again, by construction,
\[
\uu_j(x) = {\bf 0}, \qquad x > c_o j \tau,
\]
so the left hand side in~\eqref{eq:GS8} is block upper triangular. The products of the block columns in the right hand side 
are, similar to the one-dimensional case, 
\begin{align}
\sum_{l=0}^{n-1} \uv_i(x_l)^T \uv_j(x_l) &\approx \frac{1}{c_o \tau} \int_0^{x_{\max}} \hspace{-0.03in} dx \, \uv_i(x)^T \uv_j(x),
\label{eq:multiD}
\end{align}
where the left hand side is a Riemann sum approximation of the integral and the integrand is supported in $
[0,\min\{x_i,x_j\}] \subseteq [0,x_{n-1}] \subset \Omega,$ by the causality relation~\eqref{eq:ROM19}.  
Writing more explicitly \eqref{eq:multiD},
\begin{align}
\sum_{l=0}^{n-1} \uv_i(x_l)^T \uv_j(x_l) &\approx
 \frac{1}{c_o \tau} \int_0^{x_{\max}}  \hspace{-0.03in}dx \sum_{s=1}^m \bv_i(x,\bx_s^\perp)^T  \bv_j(x,\bx_s^\perp) \nonumber \\
 &\approx K \int_{\Omega}  d \bx \, \bv_i(\bx)^T  \bv_j(\bx) = K \bI_{m}, \label{eq:GS9}
\end{align}
where $d \bx = dx d \bx^\perp$, $K$ is a constant, and the accuracy of the last approximation depends on  the points $\{\bx_s^\perp\}_{1 \le s \le m}$ and on how the wave propagates.  Intuitively, the points $\{\bx_s^\perp\}_{1 \le s \le m}$ should be near the $m$ sensors in the array, 
and the approximation \eqref{eq:GS9} should hold  at least if $n$ is not too large, meaning that for $j = 0, \ldots, n-1$, the wave $\bu_j(\bx)$ has not spread out much in cross-range, but propagates downward like a beam. 

If the approximation~\eqref{eq:GS9} holds, then we have the analogue of the result in one dimension, where~\eqref{eq:GS8} is the 
block QR factorization of the block upper triangular matrix in the left hand side and the first factor in the right hand side is a multiple of the 
identity. The quasimatrix $\bV(\bx)$ is the interpolation of this matrix and is therefore approximately independent of $q(\bx)$. The approximate independence of the quasimatrix $\hat \bV(\bx)$ on $q(\bx)$ follows similarly. 

The numerical simulations in \cite{borcea2019robust,DtB,druskin2018nonlinear} confirm this statement, and they also show that the approximation  deteriorates for larger $n$. This is why in our inversion method we do not  rely on the  assumption that $\cLqR-\cLoR$ is linear in the reflectivity, and formulate instead a nonlinear minimization problem that is  solved iteratively. 
\subsection{Inversion method}
\label{sect:Inv.2}
The classic way of solving the inverse problem is to estimate  $q(\bx)$ using least squares data fit optimization
\begin{equation}
q^{\rm LS}(\bx) =  \mathop{\mbox{arg min}}\limits_{q^S \in \mathscr{Q}} 
\sum_{j=0}^{2n-1} \|\bD_j - \lb \bb,\cT_j\big(\cP(q^S)\big) \bb \rb \|_F^2,
\label{eq:IM0}
\end{equation}
where $\| \cdot\|_F$ denotes the Frobenius norm and we used equations~\eqref{eq:F7} and~\eqref{eq:Sn6}  
to write  the mapping of the guess reflectivity  $q^S$ to the data. The search space  is 
\begin{equation}
\mathscr{Q} = \mbox{span}\{\psi_1(\bx), \ldots, \psi_{N^S}(\bx)\},
\label{eq:Im2}
\end{equation}
for some carefully chosen basis functions $\{\psi_j(\bx)\}_{1 \le j \le N^S}$, with $\bx \in \Omega$. 
Problem~\eqref{eq:IM0} is clearly nonlinear, and depending on the space~\eqref{eq:Im2} the Jacobian of the mapping 
\begin{equation}
q^S \mapsto \left\{\lb \bb,\cT_j\big(\cP(q^S)\big) \bb \rb\right\}_{0 \le j \le 2n-1}
\label{eq:ImD}
\end{equation}
may be poorly conditioned, which means that \eqref{eq:IM0} should be regularized.
Following the geophysics literature \cite{nemeth1999least, dai2012multi} we refer to the 
inversion procedure for solving \eqref{eq:IM0} with a Gauss-Newton iteration as the least squares
reverse time migration (LS-RTM).

In contrast to the conventional approach \eqref{eq:IM0},
we estimate the reflectivity by the solution of the minimization problem 
\begin{equation}
q^\star(\bx) =  \mathop{\mbox{arg min}}\limits_{q^S \in \mathscr{Q}} \|\cLqR - \cLqsR\|_{F}^2,  
\label{eq:Im1}
\end{equation}
because as discussed in the previous section, the matrix $\cLqR-\cLoR$ is expected to be 
approximately linear in $q(\bx)$. This is confirmed by the numerical results, which show that the 
Gauss-Newton iteration \cite[Section 10.3]{nocedal2006numerical} converges in a few steps. 
For the sake of brevity we refer to such iteration for solving \eqref{eq:Im1} as ROM-GN. 
We emphasize that the construction of $\cLqR$ uses the data $\bD_j$ that depend non-linearly 
on the unknown reflectivity $q$.

\subsubsection{Parametrization and resolution}

For noisy data, regularization is needed in both the construction of the ROM (see \cite{borcea2019robust}) 
and in the inversion. In the numerical simulations we regularize the Gauss-Newton method using a truncated 
SVD approach. But regardless of the noise, the basis functions of the sample space~\eqref{eq:IM0} should 
be defined based on a resolution study, to avoid over parametrizing the unknown reflectivity. 
This ensures that we have a well conditioned Jacobian and also saves computational time by limiting
the dimension $N^S$ of the search space $\mathscr{Q}$. 

For the given excitation $\bb(\bx)$, the resolution depends on the location in $\Omega$, as we now explain. Let 
$\delta_j(\bx)$ be a non-negative function which integrates to one and has support centered at $\bx_j \in \Omega$, 
of diameter $\lambda/2$,  the Rayleigh resolution limit \cite[Chapter VIII]{bornprinciples} for imaging with waves 
at central wavelength $\lambda$. We may think of $\delta_j$ as an approximate Dirac $\delta(\bx - \bx_j)$. 
From sections \ref{sect:Galerk} and \ref{sect:Inv.1} we know that 
\begin{align*}
\Delta \boldsymbol{\cL}^{\RM}(\delta_j)  & \approx 
\bV^{{(0)}^T} \Delta \cL(\delta_j) \hat \bV^{(0)} \approx 
\bV^{{(0)}^T} \Delta L(\delta_j) \hat \bV^{(0)} \stackrel{\eqref{eq:In1}}{\approx} 
\frac{c(\bx_j)}{2}  \bV^{{(0)}^T} \nabla \delta_j(\bx) \cdot \hat \bV^{(0)}, \\
\Delta \boldsymbol{\cL}^{\RM}(\delta_j)^T  & \approx 
\hat \bV^{{(0)}^T} \Delta \cL(\delta_j)^T \bV^{(0)} \approx 
\hat \bV^{{(0)}^T} \Delta L(\delta_j)^T \bV^{(0)} \stackrel{\eqref{eq:In1}}{\approx} 
\frac{c(\bx_j)}{2} \hat \bV^{{(0)}^T} \nabla \delta_j(\bx) \bV^{(0)}.
\end{align*}
Here we denote $\Delta \boldsymbol{\cL}^{\RM}(\delta_j) = \boldsymbol{\cL}^{\RM}(\delta_j)-\cLoR$,
$\Delta \cL(\delta_j) =  \cL(\delta_j) - \cL(0)$ and $\Delta L(\delta_j) = L(\delta_j)-L(0)$, 
and similarly for the adjoints $\Delta \boldsymbol{\cL}^{\RM}(\delta_j)^T$, $\Delta \cL(\delta_j)^T$ 
and $\Delta L(\delta_j)^T$.
The quasimatrices  $\bV^{(0)}(\bx)$ and $\hat \bV^{(0)}(\bx)$ 
contain  the primary and dual orthonormal snapshots calculated in the 
known reference medium with zero reflectivity.
The dot denotes the inner product in $\mathbb{R}^d$ and is understood component-wise, i.e.,
\[
\nabla \delta_j(\bx) \cdot \hat \bV^{(0)}(\bx) = 
\big(\nabla \delta_j(\bx) \cdot \hat \bv_0^{(0)}(\bx), \ldots, \nabla \delta_j(\bx) \cdot \hat \bv_{n-1}^{(0)}(\bx) \big).
\]
Therefore, we have
\begin{align*}
\Delta \boldsymbol{\cL}^{\RM}(\delta_j) \Delta \boldsymbol{\cL}^{\RM}(\delta_j) ^T \approx 
\frac{c^2(\bx_j)}{4} \bV^{(0)^T} \nabla \delta_j(\bx) \cdot \hat \bV^{(0)} \hat \bV^{(0)^T} \nabla \delta_j(\bx) \bV^{(0)},
\end{align*}
in the ROM space, and in the physical space we get 
\begin{align}
\bV^{(0)}\Delta \boldsymbol{\cL}^{\RM}(\delta_j) \Delta \boldsymbol{\cL}^{\RM}(\delta_j) ^T\bV^{(0)^T}  \approx 
\frac{c^2(\bx_j)}{4} \mathbb{P}^{(0)} \nabla \delta_j(\bx) \cdot \hat{\mathbb{P}}^{(0)} \nabla \delta_j(\bx) \mathbb{P}^{(0)}.
\label{eq:Im5}
\end{align}
Here  $\mathbb{P}^{(0)} = \bV^{(0)} \bV^{(0)^T}$ is the  orthogonal projector on the space of the first $n$ snapshots in the reference medium, 
which takes any $\bphi(\bx)$ in the space of the snapshots and returns
\begin{equation}
\bV^{(0)} \bV^{(0)^T} \bphi(\bx) = \sum_{j=0}^{n-1}\bv^{(0)}_j(\bx)  \lb \bv^{(0)}_{j}, \bphi \rb.
\label{eq:OrtProj}
\end{equation}
Similarly,  $\hat{\mathbb{P}}^{(0)} = \hat \bV^{(0)} \hat \bV^{(0)^T}$ is the orthogonal projector on the space of the first $n$ dual snapshots in the reference medium.

Equation \eqref{eq:Im5} is the ROM approximation of the operator in
\begin{equation}
\Delta L(\delta_j) \Delta L(\delta_j)^T \bphi(\bx) = 
\frac{c^2(\bx)}{4} |\nabla \delta_j(\bx)|^2 \bphi(\bx) \approx 
\frac{c^2(\bx_j)}{4}  |\nabla \delta_j(\bx)|^2 \bphi(\bx),
\label{eq:Im6}
\end{equation}
which acts as point-wise multiplication.
Therefore, we define the resolution (point spread) function at point $\bx_j$ by 
\begin{equation}
\begin{array}{rcl}
\Psi_j(\bx) & = & \Big[ \bV^{(0)} (\bx) \Delta \boldsymbol{\cL}^{\RM}(\delta_j) 
\Delta \boldsymbol{\cL}^{\RM}(\delta_j)^T{\bV^{(0)^T}}(\bx)\Big]^{\frac{1}{2}}\\
& = & \| \bV^{(0)} (\bx) \Delta \boldsymbol{\cL}^{\RM}(\delta_j) \|_2,
\quad \bx \in \Omega,
\end{array}
\label{eq:Im7}
\end{equation}
where $\| \cdot \|_2$ is the Euclidian norm in the space of row-vectors $\mathbb{R}^{1 \times nm}$.

We will see from the display of the point spread function~\eqref{eq:Im7} in the numerical section that its support grows with the distance (range) $x_j$ of the point $\bx_j = (x_j,\bx_j^\perp)$. Moreover, the spreading is mostly in the cross-range direction, as expected from the classic resolution 
limits of imaging methods \cite[Chapter VIII]{bornprinciples}. We choose the basis $\{\psi_j(\bx)\}_{1 \le j \le N^S}$ of the search space 
\eqref{eq:Im2} as the continuous, piecewise linear (hat) functions on a mesh defined as follows: The discretization in range is determined by the 
range support of~\eqref{eq:Im7}, which is basically unchanged throughout the domain if the background wave speed does not have large variations.
Let $N_r$ be the number of range points. Then,  for any given range $x_j$,  with $j = 1,\ldots, N_r$, we discretize in cross-range at steps  determined by the support of~\eqref{eq:Im7}. This can be achieved for example  by seeking an approximate partition of unity using the point spread function~\eqref{eq:Im7} in the range direction and the cross-range direction, respectively. 
The result is a non-uniform (deformed rectangular) mesh with $N^S$ points, which we then triangularize to  define the hat functions.

\section{Numerical results}
\label{sect:Numer}

In this section we present two dimensional numerical results for  configurations of scatterers modeled by the reflectivity in figures \ref{fig:box_true} and \ref{fig:Crack_QTrue}. All lengths are normalized by $\ell$, the step size of the square mesh used to discretize the true medium in the time domain finite differences simulations for generating the synthetic data. The accessible boundary is modeled as sound hard and the inaccessible boundary as sound soft. Time is normalized by the sampling step $\tau$. 
The initial wave $\bb(\bx)$ is defined as in \cite[Equation (95)]{borcea2019robust} in terms of the pulse emitted by the sensors, 
which is a Ricker  wavelet. The central wavelength calculated at the reference wave speed $c_o = 1.8 \ell /\tau$ is $\lambda = 8.9 \ell$ and the smallest wavelength, at $5\%$ (i.e. -25dB) cut-off, is $4.5 \ell$. 

The first results, presented in section \ref{sect:Numer.1}, are with noiseless data. 
The second set of results, in section \ref{sect:Numer.2}, is for noisy data.

\subsection{Inversion with noiseless data}
\label{sect:Numer.1}
\begin{figure}[t]
\centering
\includegraphics[width=0.58\linewidth]{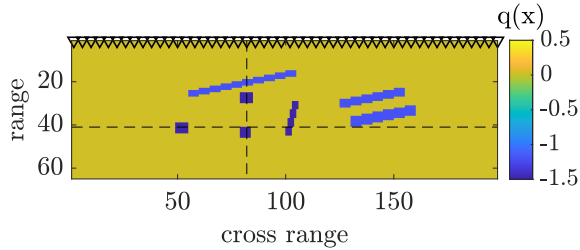}
\caption{Numerical experiment setup with sources and receivers depicted as $\nabla$. 
The dashed horizontal and vertical line are slices at which the inversion result is displayed in Fig.~\ref{fig:boxCut}. 
The range and cross range are in units of $\ell$.
The reflectivity $q(\bx)$ is dimensionless by definition~\eqref{eq:In1p} and is supported at the scatterers 
shown in different shades of blue.}
\label{fig:box_true}
\end{figure}
In the first numerical experiment we seek to estimate the reflectivity displayed in Fig.~\ref{fig:box_true}. 
The kinematic model is constant 
\[
c(\bx) = c_o = 1.8\ell/\tau,
\]
and the array has $m = 50$ sensors separated by  $4 \ell$, displayed  as triangles in the figure. The time sampling of the data is at interval $\tau$ chosen such that 
the smallest period of oscillation in the probing pulse, at $5\%$ cut-off, equals $2.5 \tau$. The data are collected at $2n = 110$ time steps, which 
leads to a data cube of dimension $110 \times 50 \times 50$. 

\begin{figure}[t]
\centering
\includegraphics[trim=0mm 4mm 8mm 7mm, clip,width=0.49\linewidth]{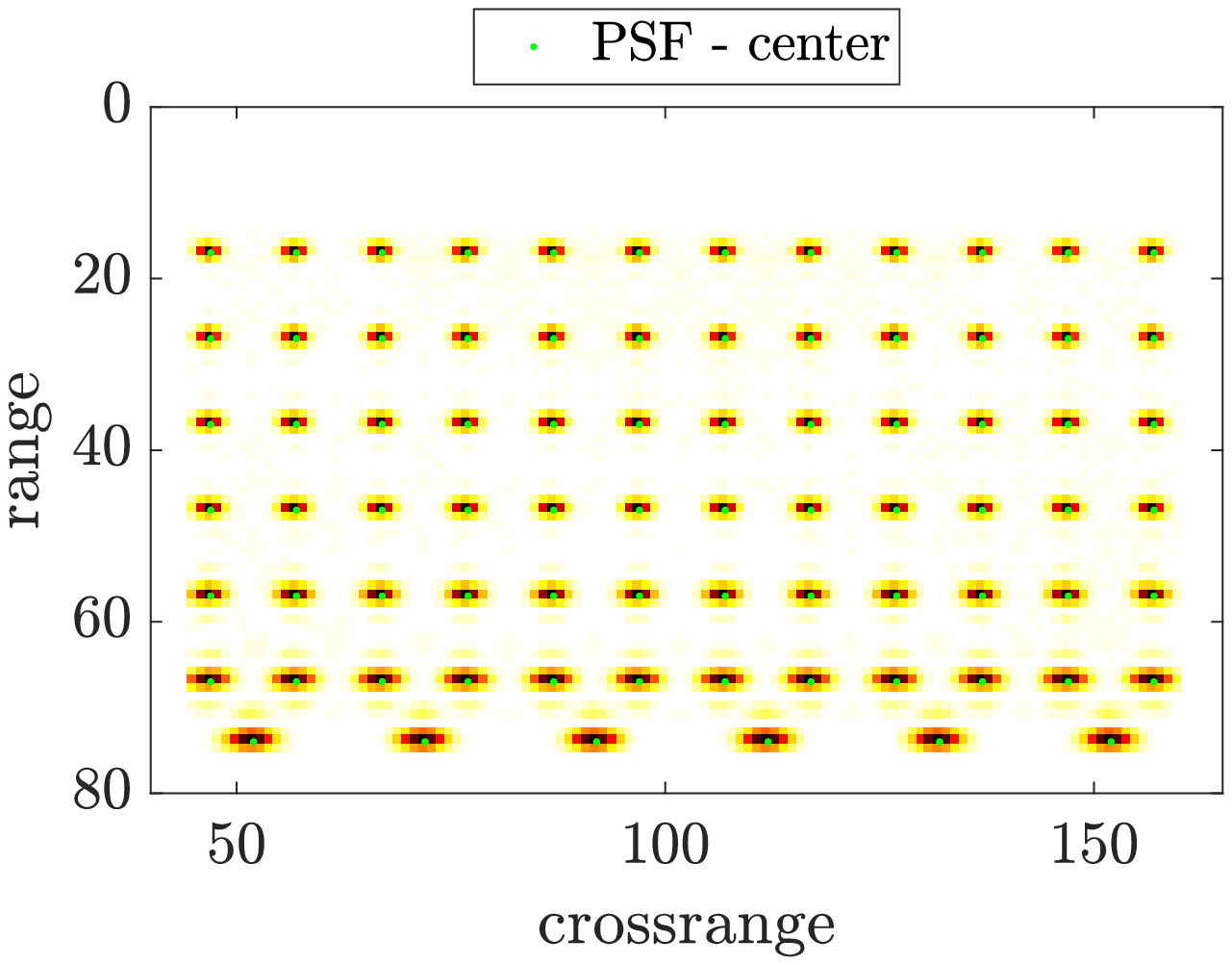}
\includegraphics[trim=0mm 4mm 8mm 7mm, clip,width=0.49\linewidth]{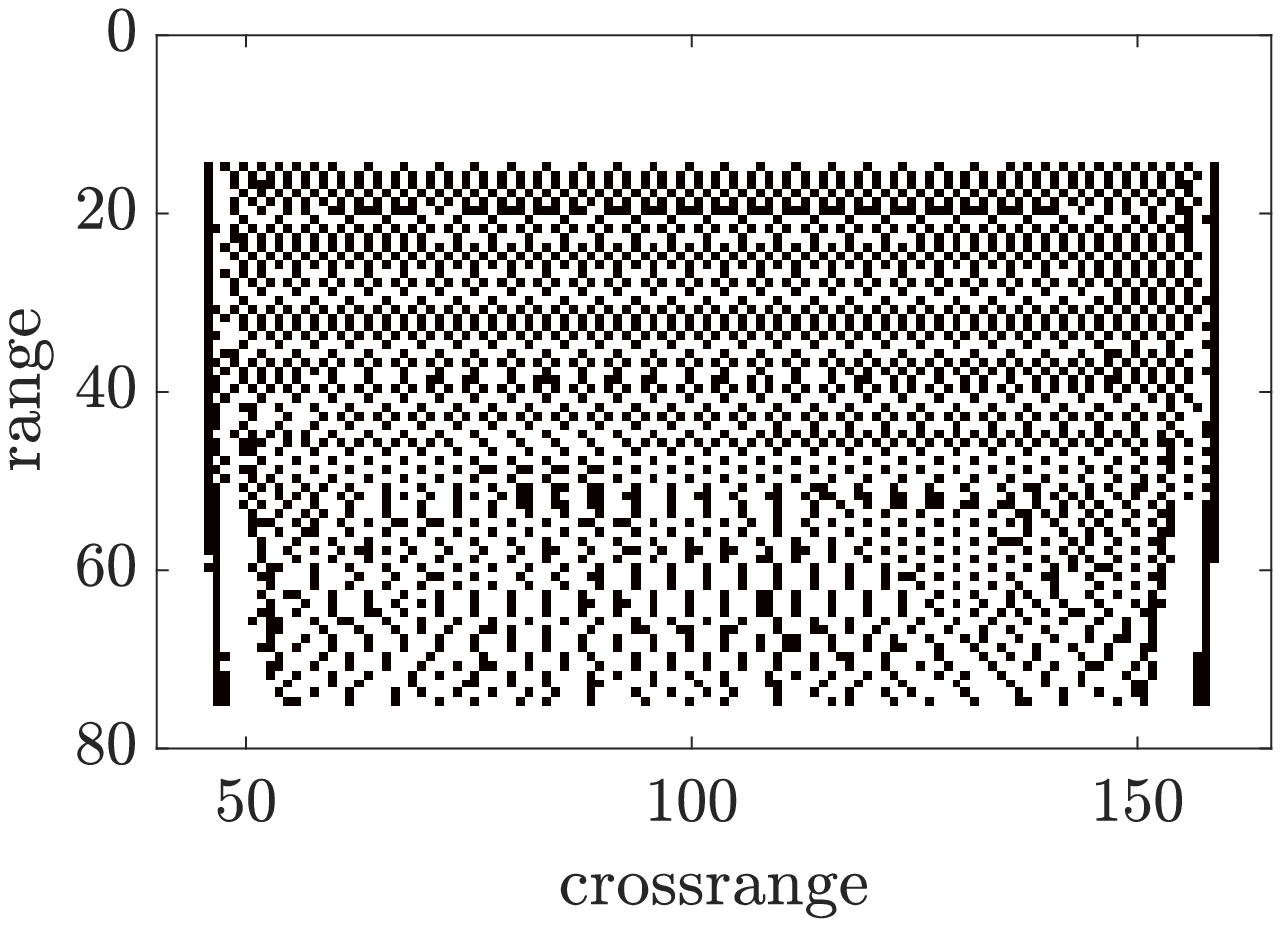}
\caption{Left: The point spread function~\eqref{eq:Im7} displayed at various 
range and cross range locations in the search region, for the setup in Fig.~\ref{fig:box_true}. 
Right: Centers of the point spread function selected from the partition of unity. 
The axes are in units of $\ell$.}
\label{fig:box_Resolution}
\end{figure}

We display in the left plot of Fig.~\ref{fig:box_Resolution} the point spread function defined in~\eqref{eq:Im7}, 
for various points in the search region,  shown with green dots. Note how its cross-range support spreads deep 
inside the medium. The parametrization 
\begin{equation}
q^S(\bx) = \sum_{j=1}^{N^S} q_j^S \psi_j(\bx),
\label{eq:parametrizedq}
\end{equation}
of the guess reflectivity is given by the continuous piecewise linear hat functions $\psi_j(\bx)$ defined on the mesh 
shown in the right plot of Fig.~\ref{fig:box_Resolution}. This mesh has the uniform spacing $c_0\tau$ in range and 
the points in the cross-range are calculated using an approximate partition of unity with the functions~\eqref{eq:Im7}. 
That is to say, at any given range $x_r$, we solved the minimization problem
\[
\min \|\balpha \|_1, \quad \mbox{such that} ~ \Big|1-\sum_{j} \alpha_j \Psi(x_r,\bx^\perp_j) \Big| \le ~\rm{tolerance},
\] 
where $\balpha$ is the vector of components $\alpha_j$, the coefficients of the point spread function at the points 
in the search cross-range interval. Due to the loss of resolution with depth, we have fewer points deep in the domain. In this example we used the tolerance of $2\%$.

\begin{figure}[t]
\centering
  \includegraphics[width=0.48 \linewidth]{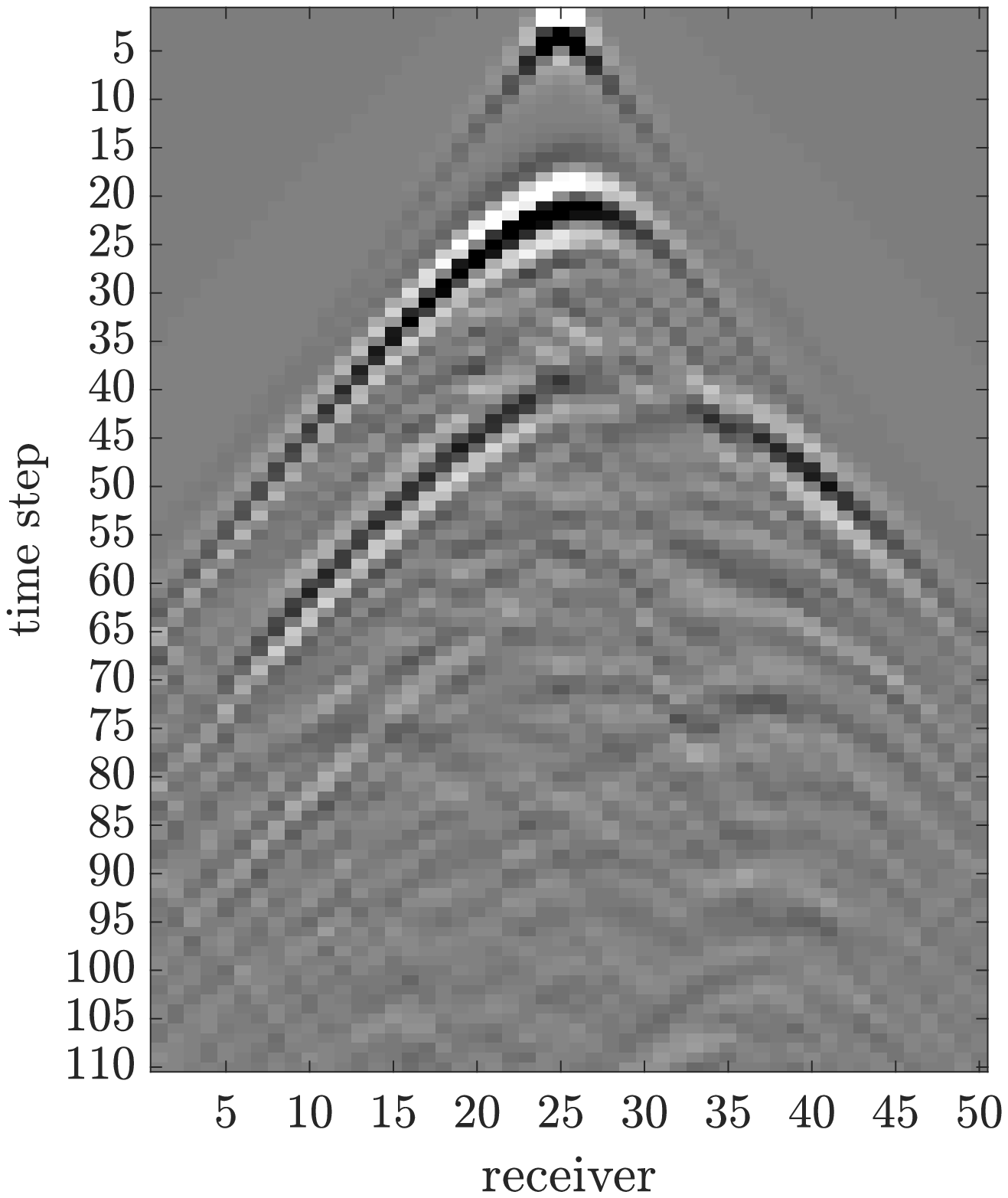}
  \includegraphics[width=0.48\linewidth]{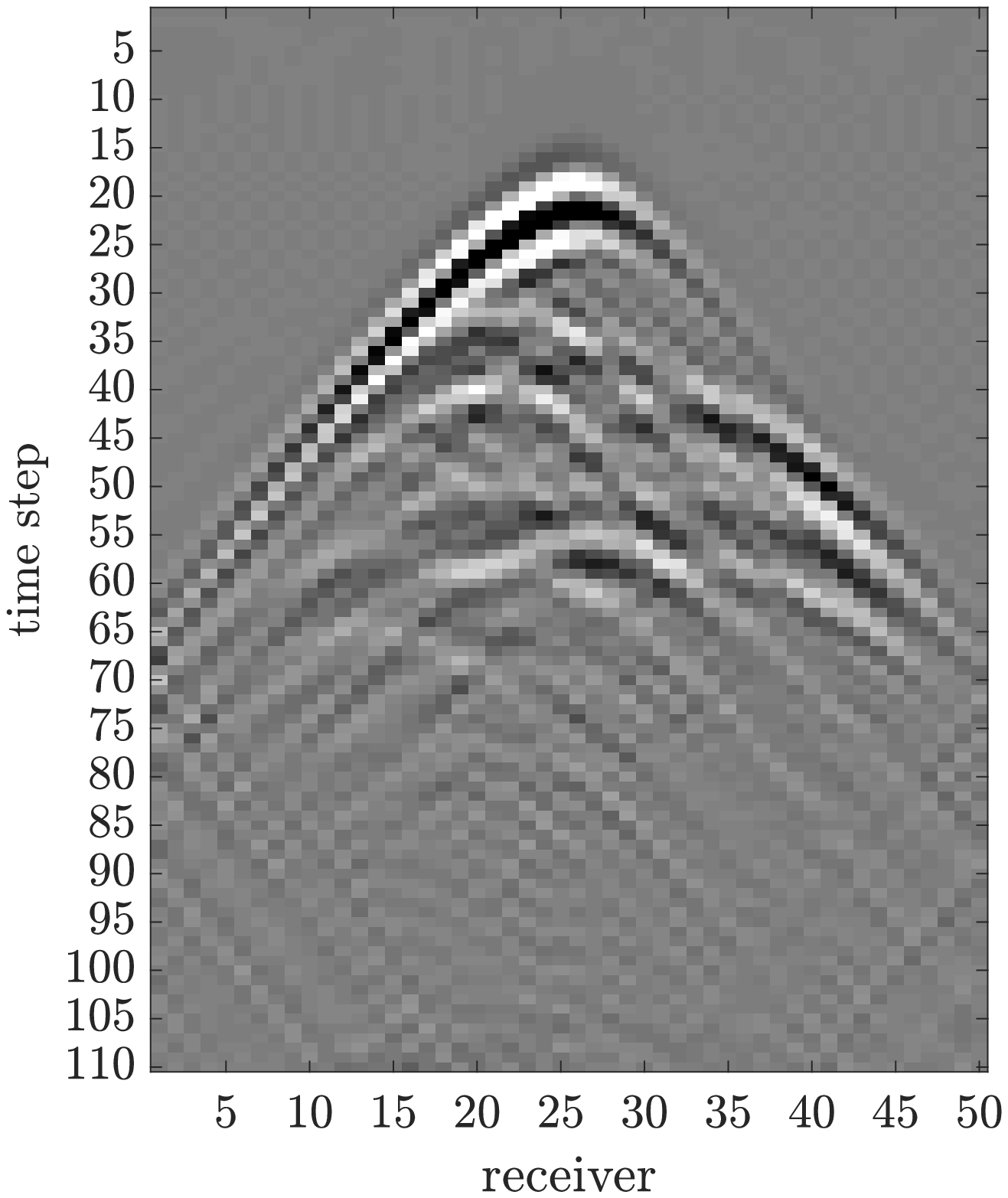}
  \caption{The $25^{th}$ column of $\bD_j$ (left plot) and of $\bD_j^{Born}$ (right plot),    for $j = 0, \ldots, 2n-1,$ as a function 
  of time in the ordinate, in units of $\tau$, and the receiver index in the abscissa. }
  \label{fig:box_data}
\end{figure}

The data obtained with the excitation from the center sensor in the array are displayed in the left plot in 
Fig.~\ref{fig:box_data}. Note that to save computational time, we made the domain $\Omega$ smaller 
than assumed in the analysis. Therefore, the fictitious boundary $\partial \Omega_{\rm inac}$ causes 
reflections that are visible at the bottom corners of the plot. In the right plot we display the 
data processed with \cite[Algorithm 1]{borcea2019robust}, which is designed to return an approximation 
of the Born (single scattering) linear data model, 
\begin{equation}
\bD_j^{Born} = 
\bD_j(0) + \bbR^T \frac{d}{d \epsilon} \cT_j \big(\boldsymbol{\cP}_{\epsilon}^{\RM}(q) \big) \big|_{\epsilon =0}\bbR, 
\qquad j = 0, \ldots, 2n-1.
\label{eq:Born1} 
\end{equation}
Here $\bD_j(0)$ are the data simulated for the reference medium with no reflectivity and the right hand side is calculated using 
\begin{align}
\boldsymbol{\cP}_{\epsilon}^{\RM}(q) & = 
\bI_{nm} - \frac{\tau^2}{2} \boldsymbol{\cL}_\epsilon^{\RM}(q) \boldsymbol{\cL}_{\epsilon}^{\RM}(q) ^T 
\approx \boldsymbol{\cP}^{\RM}(\epsilon q), \\
\boldsymbol{\cL}^{\RM}_\epsilon (q) & = 
\cLoR + \epsilon \big( \cLqR - \cLoR\big) \approx \boldsymbol{\cL}^{\RM}(\epsilon q).
\end{align}

\begin{figure}[t]
\centering
\includegraphics[trim=0mm 0mm 0mm 0mm, clip,width=0.47\linewidth]
{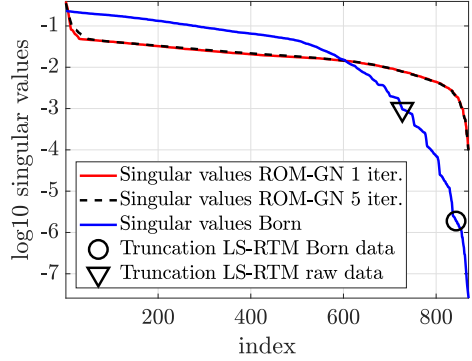}
\caption{Squared singular values of the Jacobian of the mapping 
$q^S \mapsto \{\bD^{Born}_j\}_{0 \le j \le 2n-1}$ (blue curve) and of the mapping 
$q^S \mapsto \boldsymbol{\mathcal{L}}^{\rm ROM}(q^S)$ for the first and fifth Gauss-Newton iterates 
(red and black curves, respectively). The LS-RTM formulations with Born and raw data are regularized using 
truncated SVD at the values indicated with the circle and triangle, respectively.}
\label{fig:boxSingVal}
\end{figure}

\begin{figure}[t]
\centering
\includegraphics[trim=0mm 13mm 8mm 17mm, clip,width=0.49\linewidth]{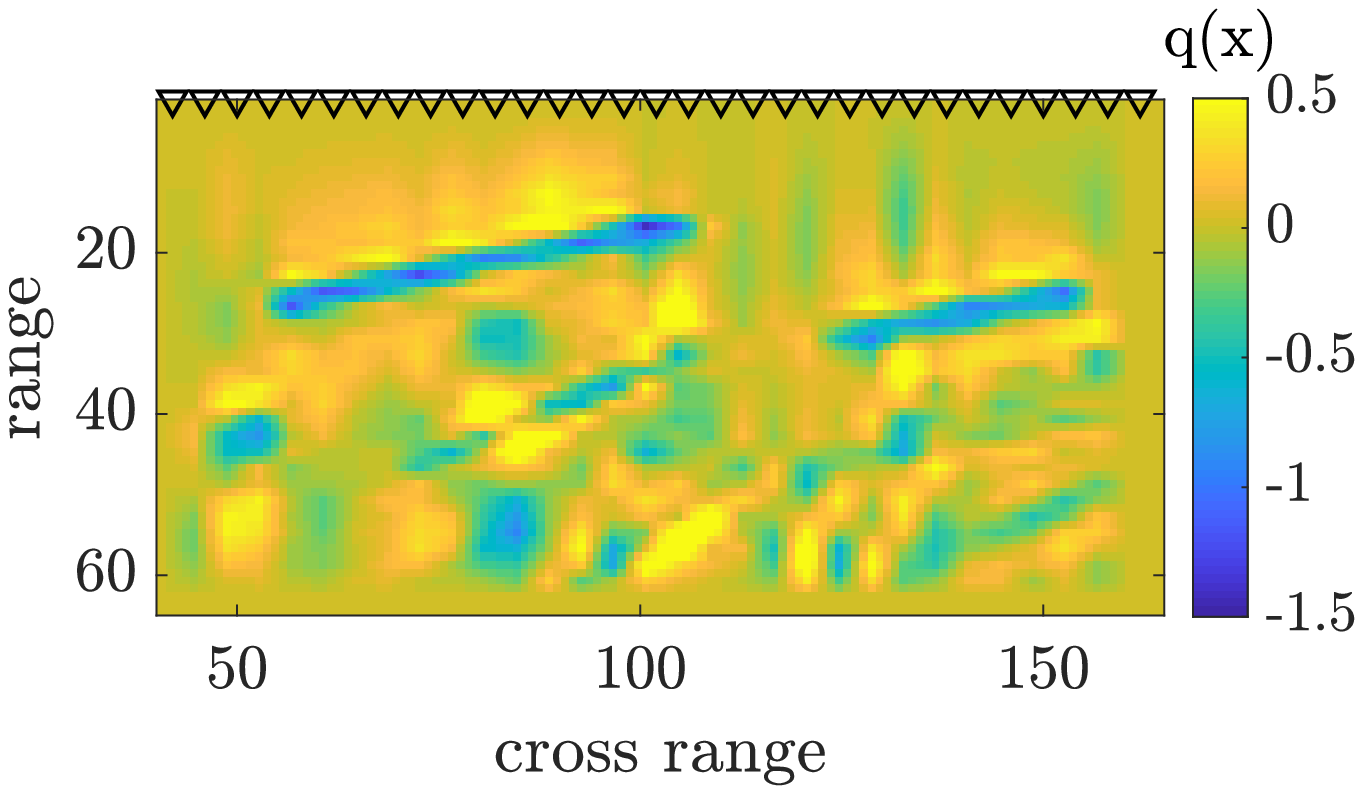}
\includegraphics[trim=0mm 13mm 8mm 17mm, clip,width=0.49\linewidth]{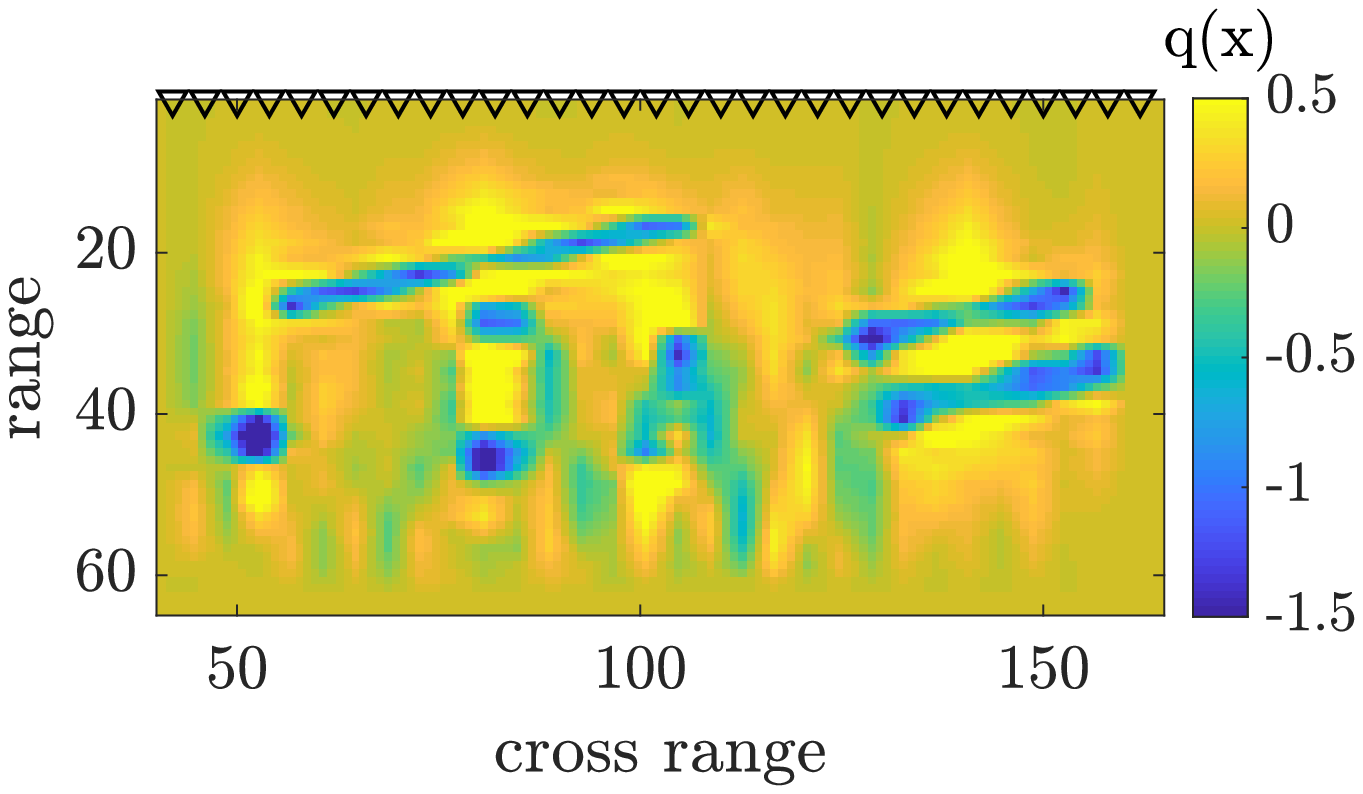}
\includegraphics[trim=0mm 13mm 8mm 17mm, clip,width=0.49\linewidth]{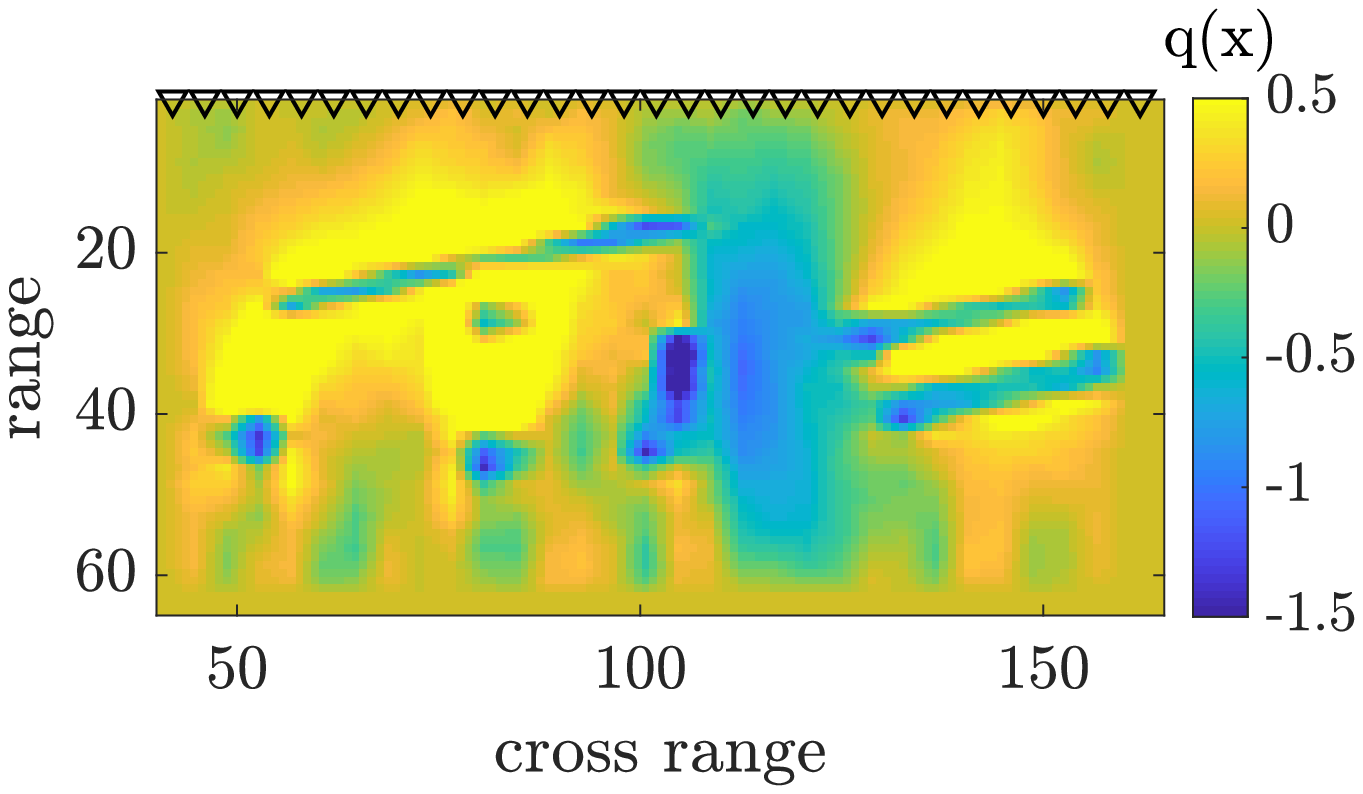}
\includegraphics[trim=0mm 13mm 8mm 17mm, clip,width=0.49\linewidth]{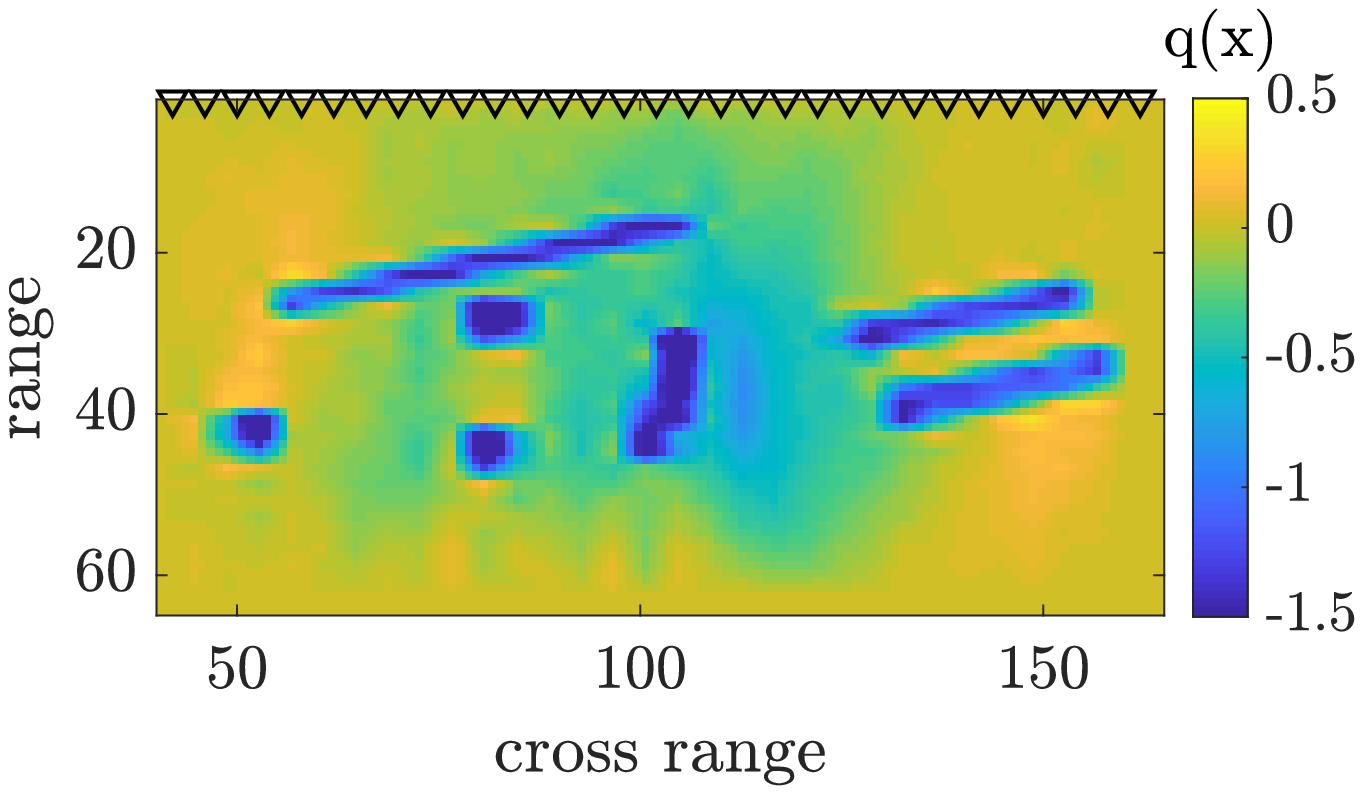}
\caption{Inversion results comparison.
Top row: LS-RTM with raw data (left plot) and the transformed (Born) data (right plot). 
Bottom row: ROM-GN after 1 iteration (left plot) and 5 iterations (right plot). 
The axes are in units of $\ell$.}
\label{fig:LSbox}
\end{figure}

\begin{figure}[ht]
\centering
\includegraphics[trim=0mm 0mm 0mm 4mm, clip,width=0.48\linewidth]
{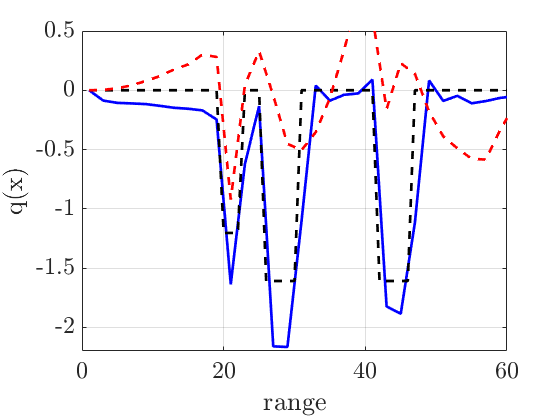}
\includegraphics[trim=0mm 0mm 0mm 4mm, clip,width=0.48\linewidth]
{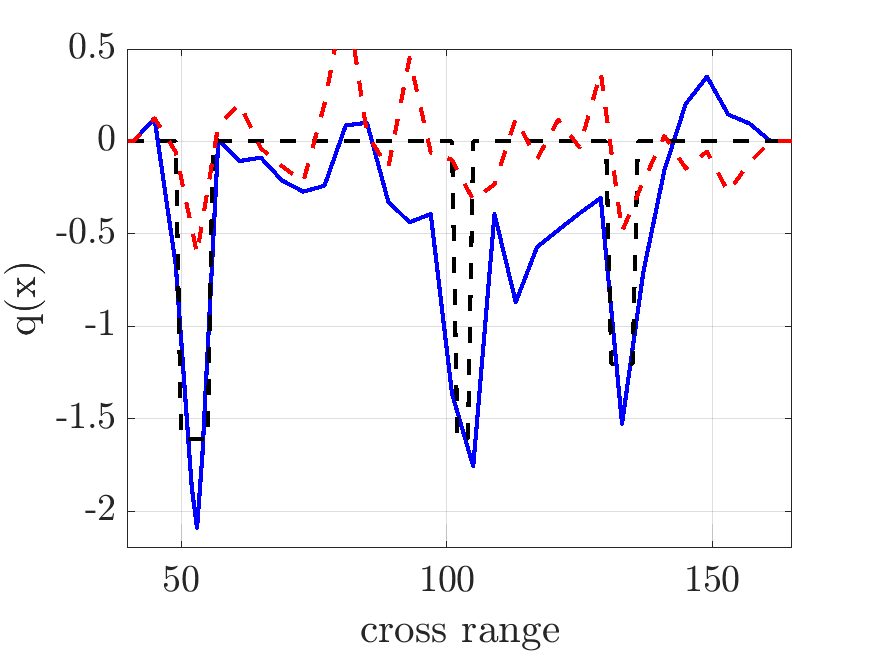}
\caption{Comparison of inversion results along the range (left plot) and cross range (right plot) slices 
taken at the lines shown in Fig.~\ref{fig:box_true}. The true reflectivity $q(x)$ is dashed black,
LS-RTM inversion result is dashed red, ROM-GN inversion result is solid blue. The abscissa is in units of $\ell$.} 
\label{fig:boxCut}
\end{figure}

We begin the comparison between the conventional LS-RTM and the proposed ROM-GN with a study 
of Jacobians of the corresponding mappings. As shown in Fig.~\ref{fig:boxSingVal}, 
the Jacobian of $q \mapsto \{\bD_j^{Born}\}_{0 \le j \le 2n-1}$ has worse conditioning compared to 
the Jacobian of $q \mapsto \cLqR$. Consequently, while LS-RTM required regularization via SVD truncation, 
we did not use regularization in ROM-GN
\footnote{However, we used the algorithm described in \cite{borcea2019robust}, based on a truncated SVD of the 
mass matrix,  for the computation of the ROM and the transformation \eqref{eq:Born1}.}.

In Fig.~\ref{fig:LSbox} we compare the inversion results for both LS-RTM and ROM-GN approaches. 
For LS-RTM (top plots) we performed a single Gauss-Newton iteration using as input both the raw data 
$\{\bD_j\}_{0 \le j \le 2n-1}$ (top left plot) and the processed (Born) data~\eqref{eq:Born1} (top right plot), 
which is intended to transform the problem into linear least squares. The image with the latter is better, 
as expected, because the multiple scattering effects have been removed approximately. Nevertheless, we observe 
image artifacts, due to the ill-conditioning of the Jacobian of the mapping $q \mapsto \{\bD_j^{Born}\}_{0 \le j \le 2n-1}$. 
In our experience, performing more Gauss-Newton iterations does not lead to an improved image, 
mostly because the transformed data \eqref{eq:Born1} are a very good approximation of the linearized (Born) data. 

The reflectivity obtained with ROM-GN is shown in the bottom two plots in Fig.~\ref{fig:LSbox} 
both after a single (bottom left plot) and five (bottom right plot) Gauss-Newton iterations, where convergence 
was achieved. Note that the shape of the scatterers is recovered well, because the operator~\eqref{eq:In1} 
depends on the gradient of the reflectivity. Thus, it is easier to get the jumps of $q(\bx)$ than its smooth part. 
However, after five iterations the magnitudes of the scatterers are also recovered very well, as clearly seen in 
the range and cross range slice plots shown in Fig.~\ref{fig:boxCut}.

The ROM-GN uses the raw data $\{\bD_j\}_{0 \le j \le 2n-1}$ and thus takes into account multiple scattering 
effects, which contain valuable information about the reflectivity $q$ that  may not be captured in 
$\{\bD_j^{Born}\}_{0 \le j \le 2n-1}$. Thus, we observe a clear advantage of our ROM-GN approach at 
recovering both the shapes and magnitudes of scatterers compared to the conventional LS-RTM.

\subsection{Inversion with noisy data}
\label{sect:Numer.2}

\begin{figure}[ht]
  \centering
\includegraphics[trim=0mm 1mm 8mm 1mm, clip, width=0.49\linewidth]{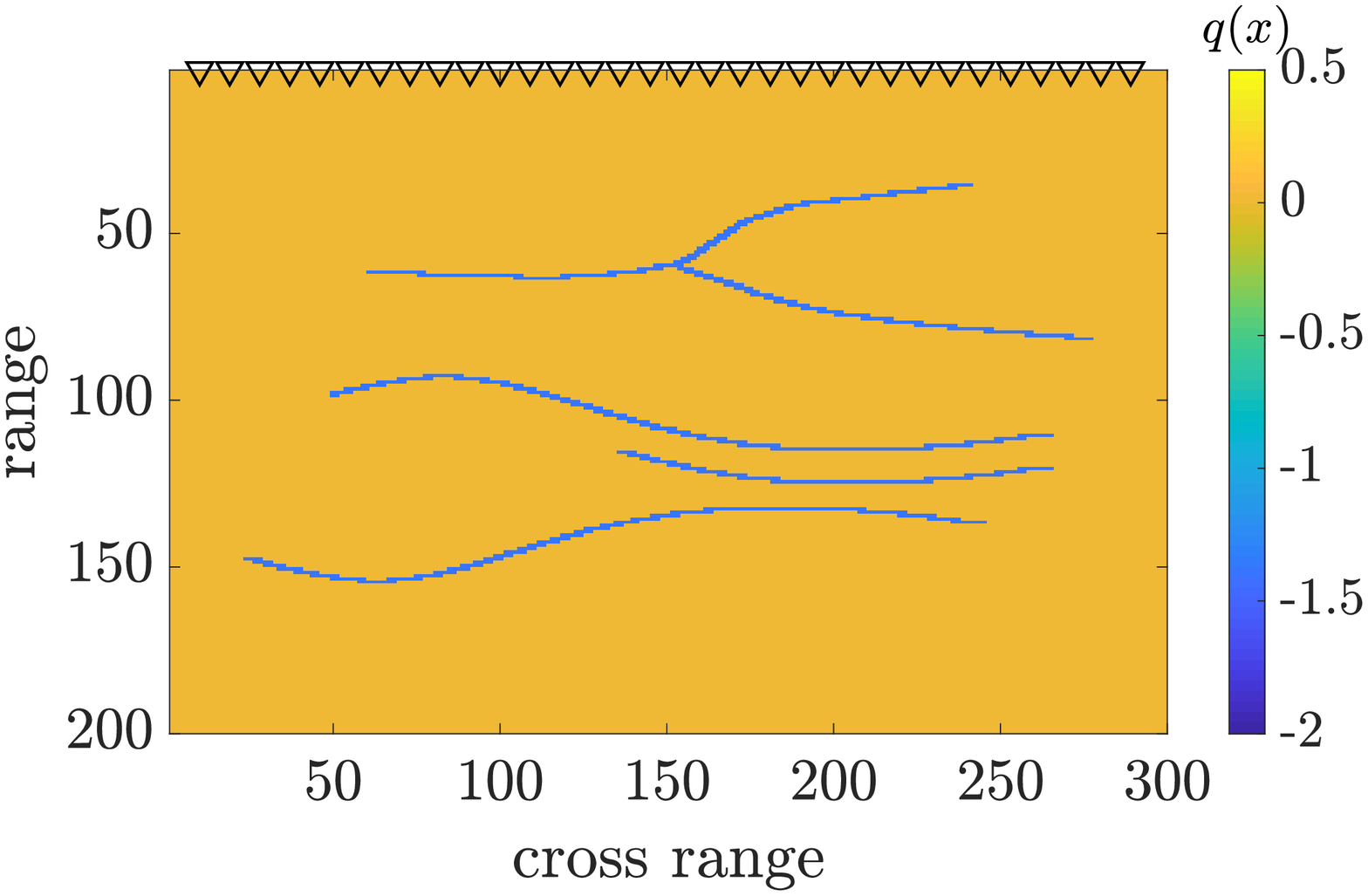}
\includegraphics[trim=0mm 1mm 8mm 1mm, clip,width=0.49\linewidth]{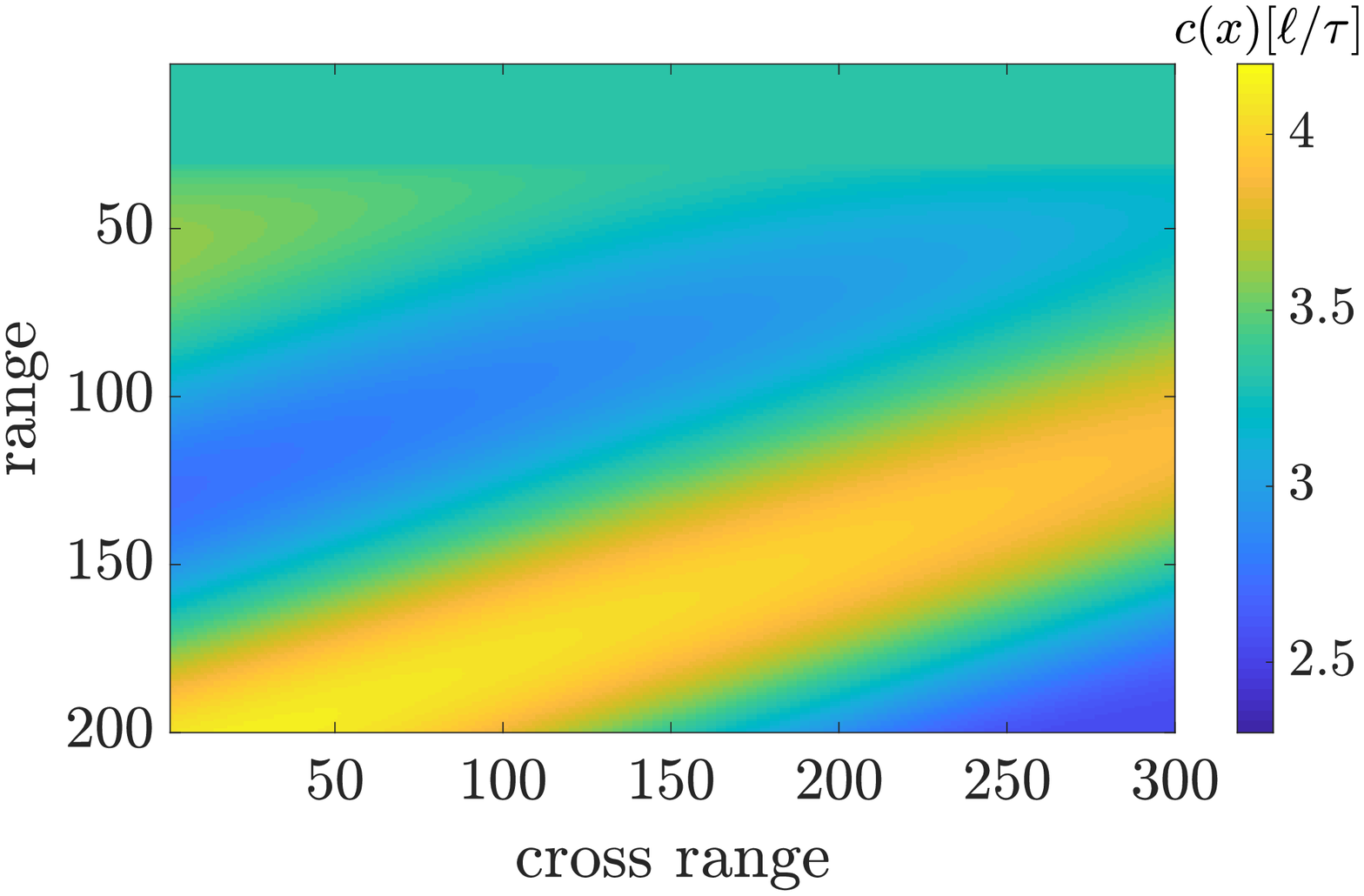}
  \caption{Left: True reflectivity $q(x)$ modeling multiple fractures (thin regions with smaller acoustic impedance). Right: Kinematic model, with 
  $c(\bx)$ displayed in units $\ell/\tau$. The array is shown on the top in the left plot. The axes are in units of $\ell$.}
  \label{fig:Crack_QTrue}
\end{figure}

\begin{figure}[ht]
  \centering
\includegraphics[width=0.48\linewidth]{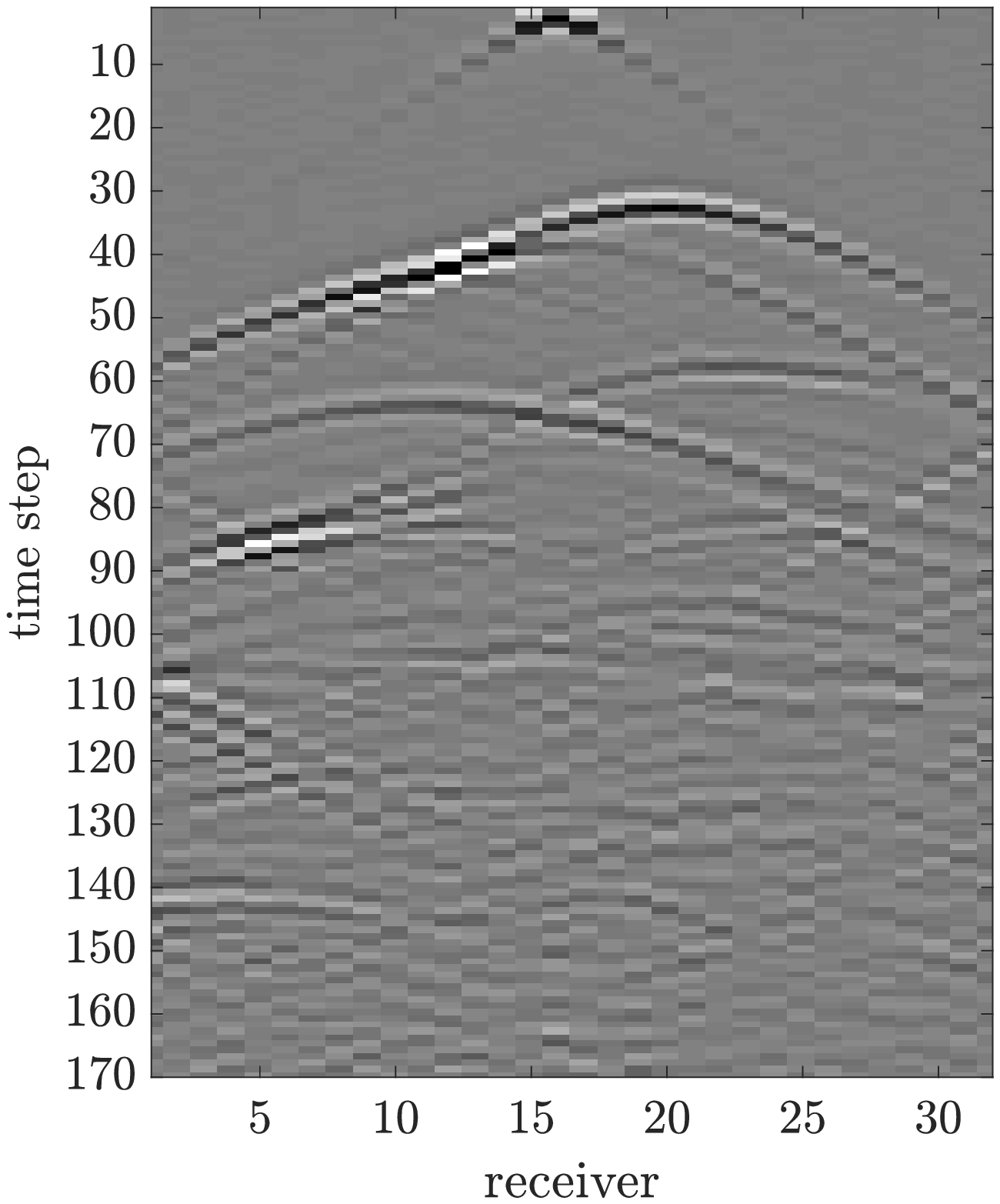}
\includegraphics[width=0.48\linewidth]{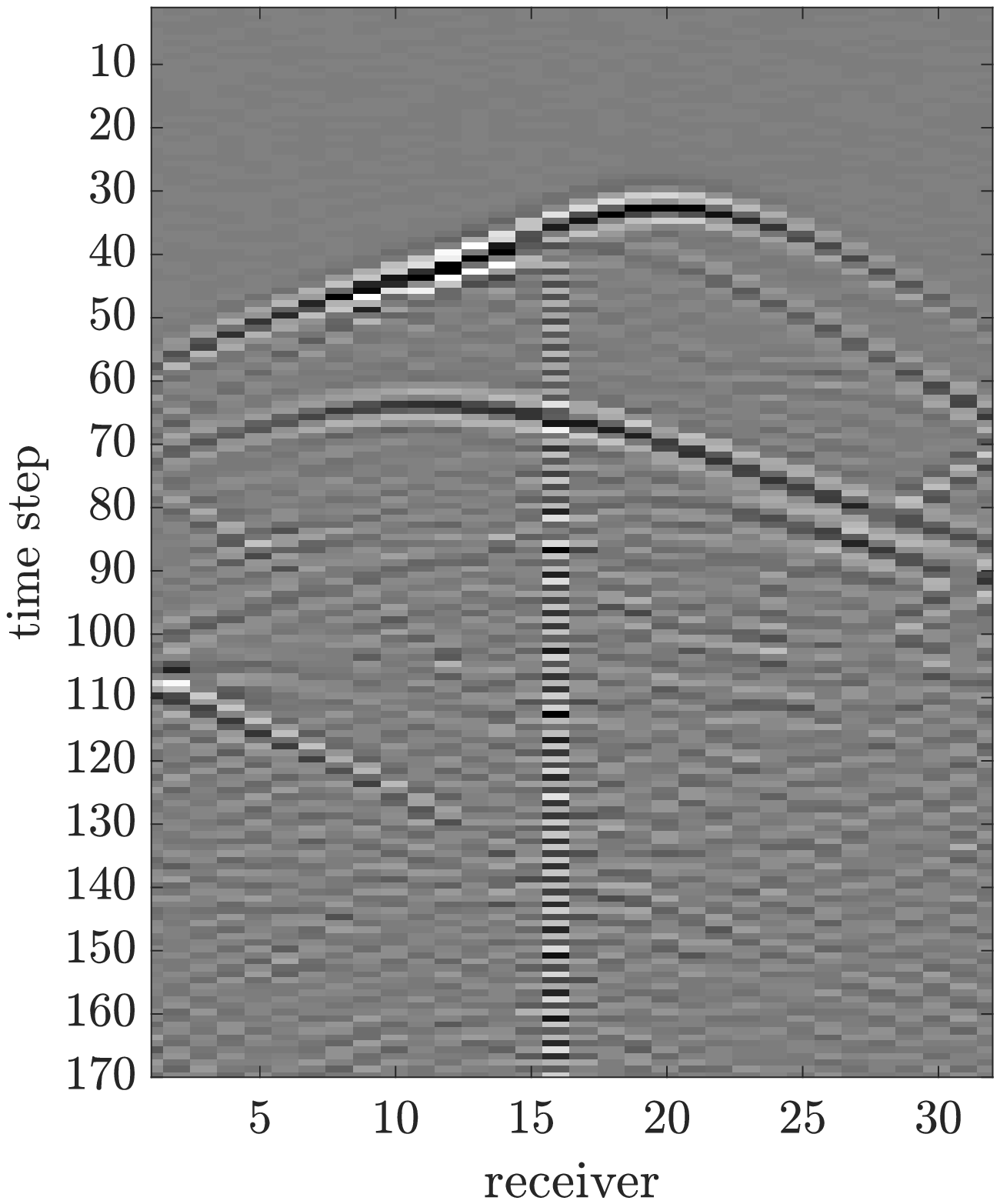}
\caption{The $16^{th}$ columns of $\bD_j$ (left plot) and of $\bD_j^{Born}$ (right plot), 
for $j = 0, \ldots, 2n-1,$ as a function of time in the ordinate, in units of $\tau$, and the receiver index in the abscissa. }  
\label{fig:DataCrack}
\end{figure}

The second numerical experiment is motivated by the application of non-destructive testing, and seeks to estimate 
multiple fractures modeled by the reflectivity displayed in Fig.~\ref{fig:Crack_QTrue}. The excitation is the same 
as in the previous experiment, except that the array has $32$ sensors separated by $8 \ell$  and the 
kinematic model is no longer constant. The data are displayed in the left plot of Fig.~\ref{fig:DataCrack} 
and are contaminated with $5\%$ additive, white Gaussian noise. They are sampled at $2n = 170$ time steps, 
at interval $\tau$ calculated so that the smallest period of oscillation in the probing pulse, at $5\%$ cut-off, 
equals $2 \tau$. The transformed data~\eqref{eq:Born1} are displayed in the right plot of Fig.~\ref{fig:DataCrack}. 
We note in particular the multiple echo around time $90 \tau$ that is suppressed after the transformation.

\begin{figure}[h]
\centering
\includegraphics[trim=0mm 7mm 8mm 10mm, clip,width=0.47\linewidth]{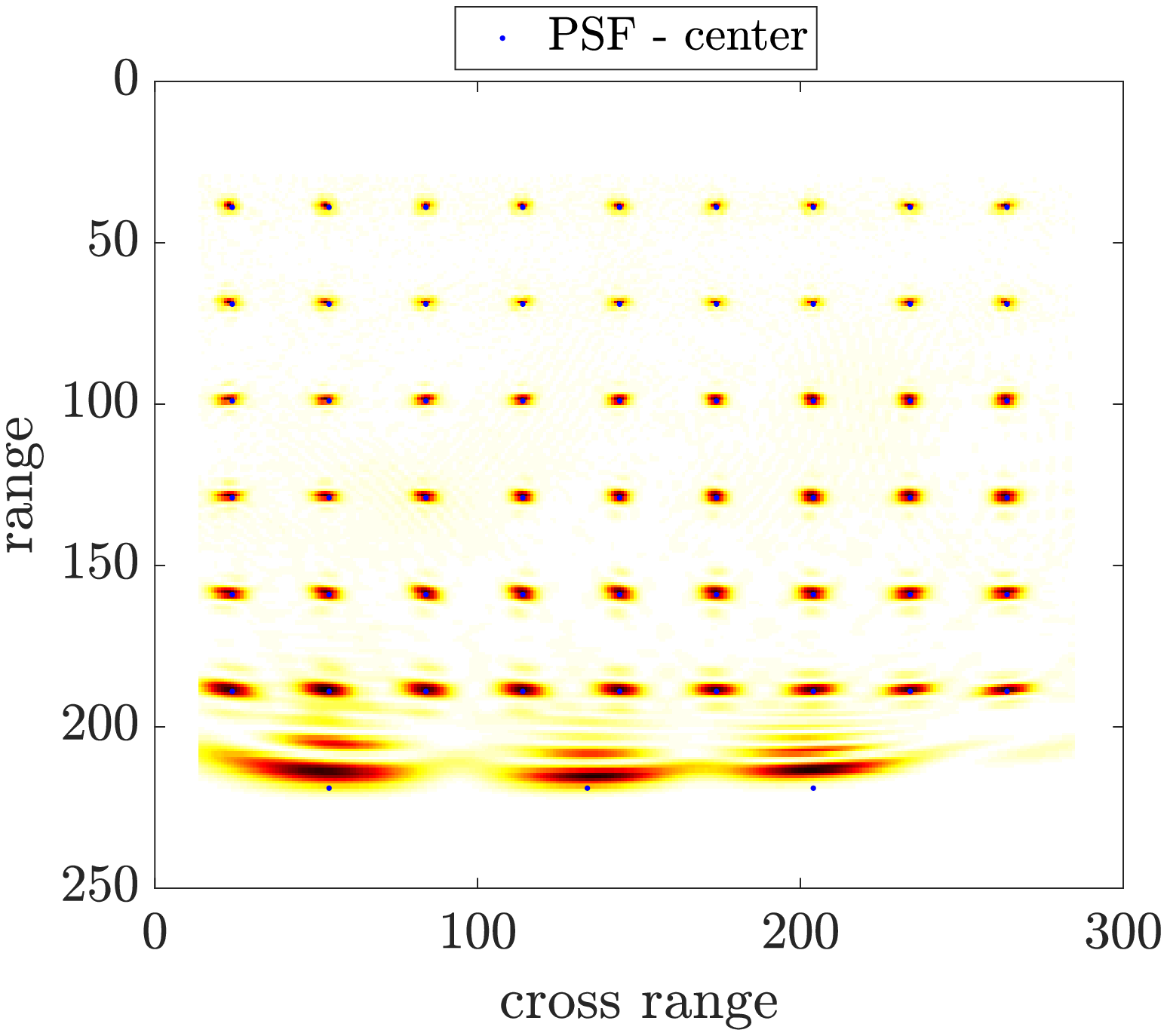}
\includegraphics[trim=0mm 8mm 8mm 10mm, clip,width=0.47\linewidth]{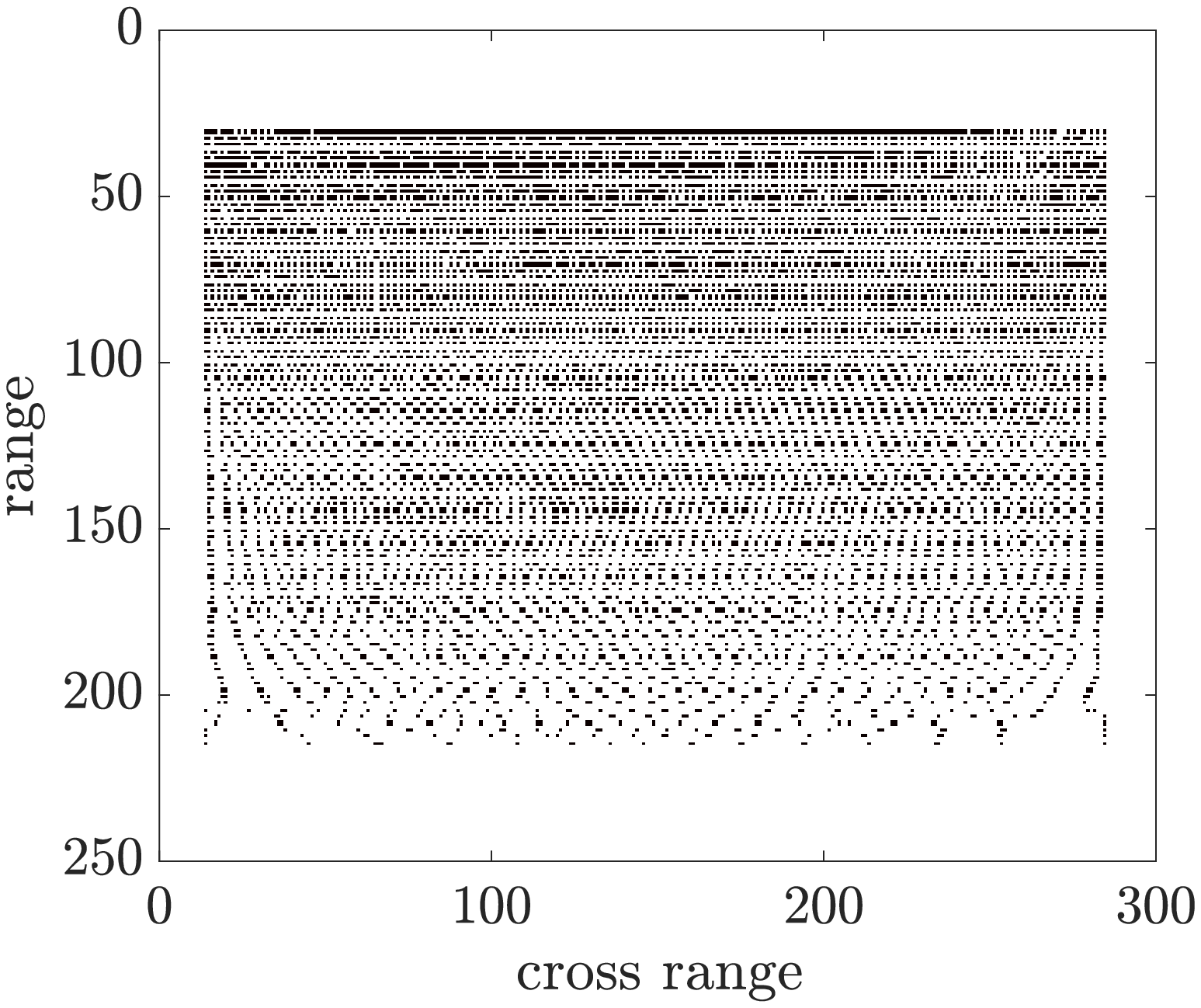}
\caption{Left: The point spread function~\eqref{eq:Im7} displayed at various range and cross range locations in the 
search region, for the setup in Fig.~\ref{fig:Crack_QTrue}. Right: Centers of the point spread function selected from the 
partition of unity. The axes are in units of $\ell$.  }
\label{fig:Crack_Resolution}
\end{figure}

We display in Fig.~\ref{fig:Crack_Resolution} the point spread function $\Psi_j(\bx)$ defined in~\eqref{eq:Im7}, 
at different locations $\bx_j$ in the search domain, indicated by the dots. Note that the spread function looks 
different than in Fig.~\ref{fig:box_Resolution} due to the variable kinematic model. The mesh calculated as 
explained in the previous section is shown in the right plot in Fig.~\ref{fig:Crack_Resolution}.

\begin{figure}[h]
  \centering
\includegraphics[trim=0mm 1mm 8.5mm 1mm, clip,width=0.48\linewidth]{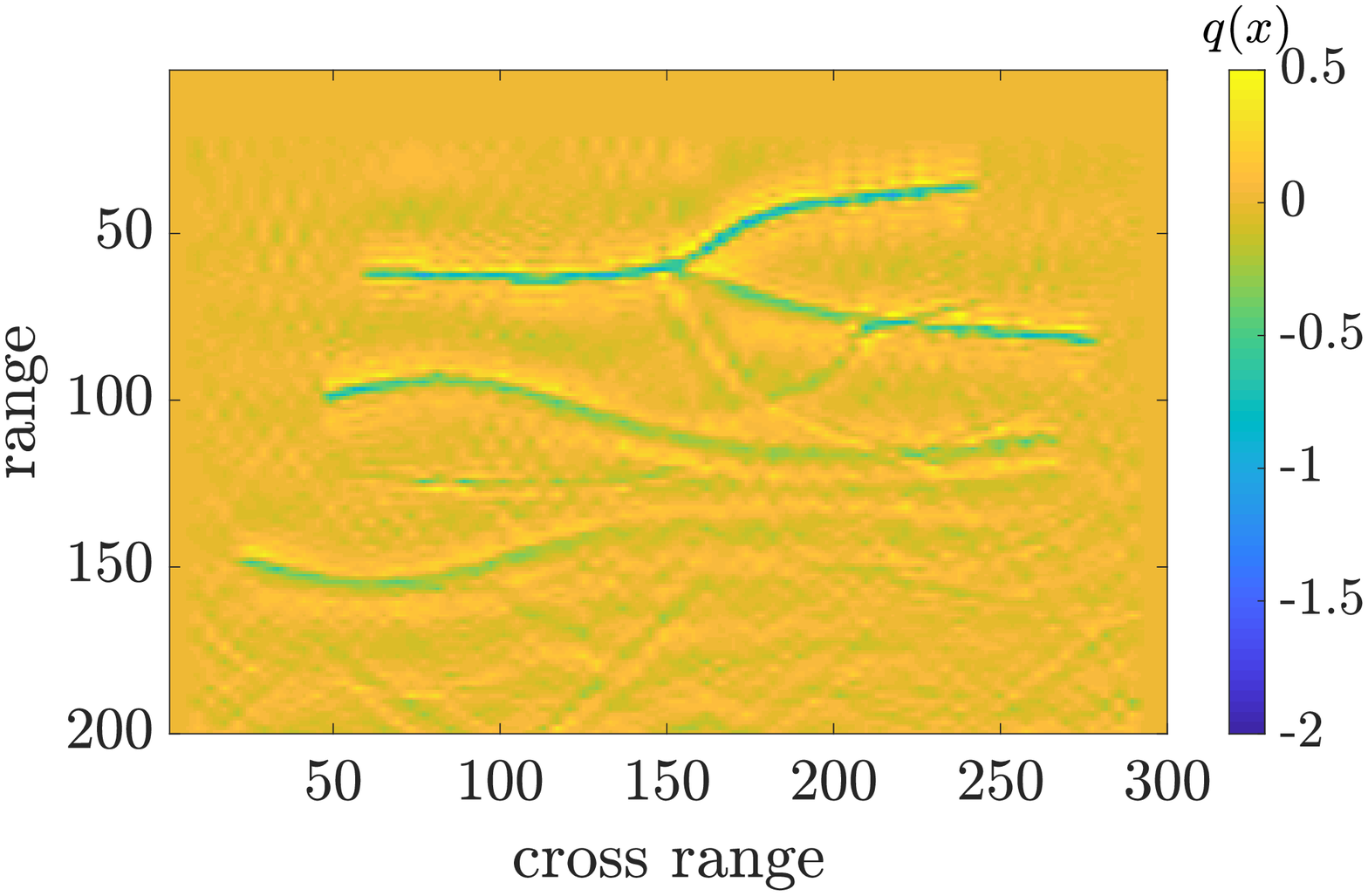}
\includegraphics[trim=0mm 1mm 8mm 1mm,  clip,width=0.48\linewidth]{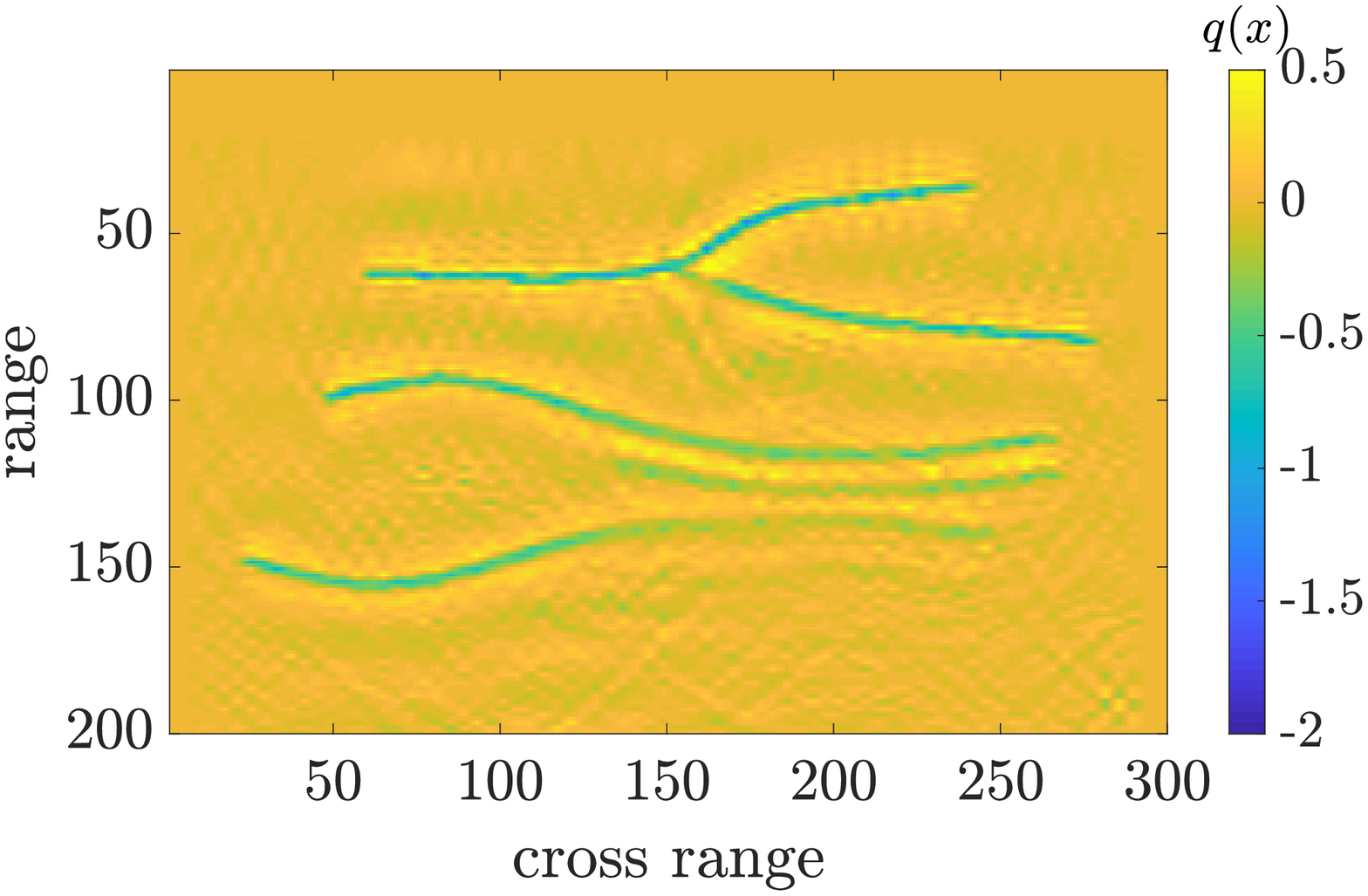} \\
\includegraphics[trim=0mm 1mm 8mm 1mm,  clip,width=0.49\linewidth]{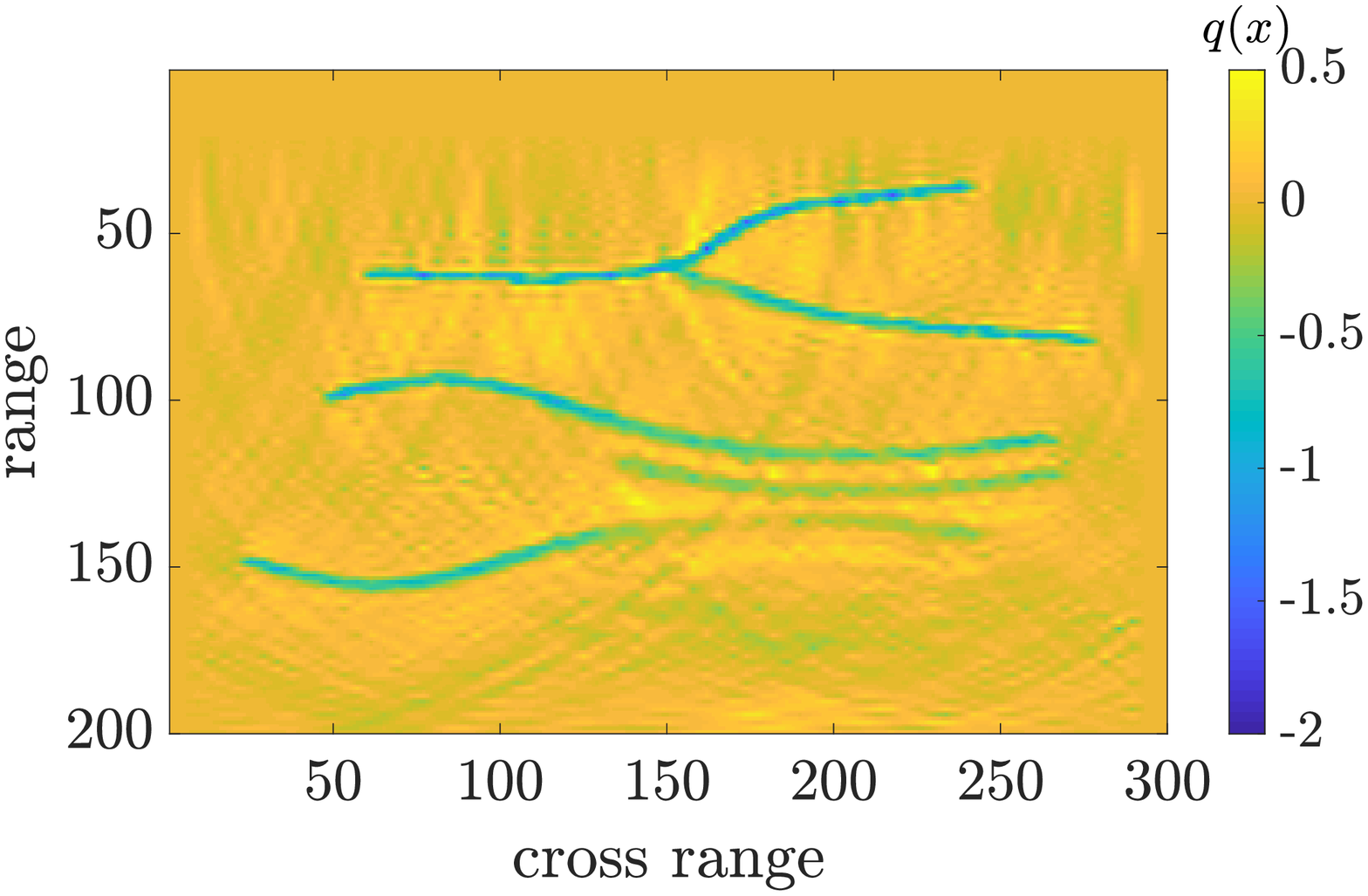}
\includegraphics[trim=0mm 1mm 8mm 1mm, clip, width=0.49\linewidth]{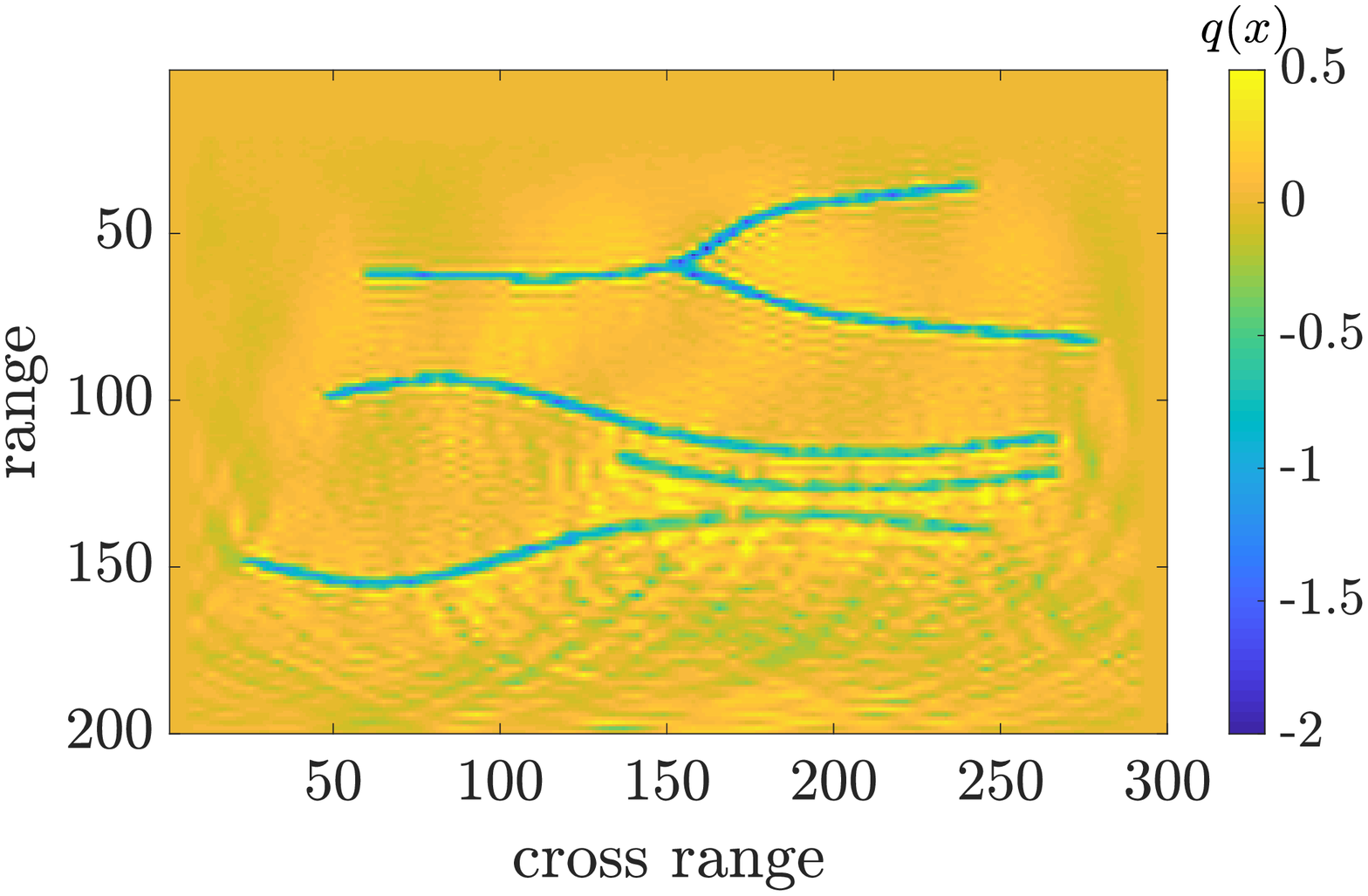}  
\caption{Inversion results comparison.
Top row: LS-RTM with raw data (left plot) and the transformed (Born) data (right plot). 
Bottom row: ROM-GN after 1 iteration (left plot) and 5 iterations (right plot). 
The axes are in units of $\ell$.}
\label{fig:Crack_LS}
\end{figure}

We compare the inversion results for LS-RTM and ROM-GN with noisy data in Fig.~\ref{fig:Crack_LS}. 
We observe that the LS-RTM inversion results are better than in the previous experiment because the 
reflectivity contrast is not as strong. 
Nevertheless, even when given the transformed (Born) data (top right plot in Fig.~\ref{fig:Crack_LS}), 
the LS-RTM does not recover the two bottom cracks very clearly. A much better inversion result is obtained 
with ROM-GN after five iterations which resolves all cracks, as shown in the bottom right plot in Fig.~\ref{fig:Crack_LS}. 
Note that both LS-RTM and ROM-GN methods were regularized with a truncated SVD of the Jacobian.


\section{Summary}
\label{sect:Sum}
We introduced a novel method for the inverse scattering problem, where the goal is to estimate reflective 
structures in a  medium from data gathered by an active array of sensors. These sensors emit waves that 
propagate through the medium and measure the backscattered returns at $2n$ time instants separated by 
an appropriately chosen  interval. The new algorithm is based on a reduced order model (ROM) of the wave 
propagator operator. This operator maps the wave from one time instant to the next, and is unknown in 
inverse scattering. However, the ROM can be calculated from the measurements at the array. We described 
the ROM for a generic hyperbolic system and showed that it corresponds to a Galerkin projection of the 
propagator operator on the space spanned by the wave at the first $n$ times instants. We analyzed the 
ROM in the Galerkin framework, and used the results to motivate the new inversion method. We described 
the implementation of the method in the context of inverse scattering for sound waves, and assessed its 
performance with numerical simulations. Compared to the conventional nonlinear least squares
data fit minimization, the new inversion method is almost unaffected by the multiple scattering
effects. It recovers robustly the locations, shapes and magnitudes of scatterers in a very small number of iterations.

\section*{Acknowledgements}
This material is based upon research supported in part by the
U.S. Office of Naval Research under award number N00014-17-1-2057 to
Borcea and Mamonov. Borcea also acknowledges support from the AFOSR
award FA9550-18-1-0131 and Mamonov acknowledges support from the National Science Foundation Grant DMS-1619821.

\vspace{0.1in}
\appendix
\textbf{Apendixes.} The next appendixes justify the wave and data model and contain the proofs of the 
ROM properties stated in Theorems \ref{thm.1}--\ref{thm.7} in section \ref{sect:Galerk.4}.

\section{The initial condition and data model}
\label{ap:source}

Typically, a wave source is modeled as a force term in the right hand side of the wave equation, and the wave field satisfies homogeneous initial conditions. In this appendix we explain how such a typical formulation can be transformed in problem  \eqref{eq:F1}--\eqref{eq:F3} and  also justify the data model \eqref{eq:F7}.

To simplify the presentation, we assume throughout the appendix that the wave field is scalar (i.e., neglect polarization), 
so $s = 1, \ldots, m$ indexes the location of the point-like sensors in the array which emit the same  pulse $f(t)$ 
supported around $t = 0$. The wave generated by the source at $\bx_s$ is denoted by  $w^{(s)}(t,\bx)$ and solves 
the wave equation 
\begin{align} 
\partial_t^2 w^{(s)}(t,\bx) + L(q) L(q)^T w^{(s)}(t,\bx) &= \partial_t f(t) \delta(\bx-\bx_s), \qquad \bx \in \Omega, \quad t \in \RR, \label{eq:S.1}\\
w^{(s)}(t,\bx)& = 0, \qquad t \ll 0, \label{eq:S.2}
\end{align}
with the same homogeneous boundary conditions  as in problem  \eqref{eq:F1}--\eqref{eq:F3}.  
We suppose that $f(t)$ is real valued, with non-negative Fourier transform\footnote{The technical condition 
\eqref{eq:S3} is needed in the derivation below but it is not a big restriction, because in imaging one usually 
convolves the received signals with the time reversed version of the emitted waveform. This is known as pulse 
compression in radar imaging \cite{curlander1991synthetic}, and it is essential because due to antenna power 
considerations, the emitted waveforms are usually long signals (chirps) $F(t)$. Using the time convolution 
$\star_t$ they are transformed into short pulses $f(t) = F(-t) \star_t F(t)$ with Fourier transform 
$\hat f(\om) = |\hat F(\om)|^2 \ge 0.$}
\begin{equation}
\hat f(\om) = \int_{-\infty}^\infty dt \, e^{i \om t} f(t) \ge 0, \qquad \forall ~\om \in \mathbb{R}.
\label{eq:S3}
\end{equation}

We can write formally the explicit expression of $w^{(s)}(t,\bx)$ using the spectral decomposition of the operator 
\begin{equation}
A := L(q)L(q)^T,
\end{equation}
which is self-adjoint and coercive.  Following \cite[Theorem 4.12]{mclean} we conclude that the eigenvalues of $A$ 
are ordered as $0 < \la_1 \le \la_2 \le \ldots$, with $\la_l \to \infty$ as $l \to \infty$, and the eigenfunctions 
$\{y_l(\bx)\}_{l \ge 1}$ form a complete orthonormal system in $L^2(\Omega)$. Therefore, we can express the wave as
\begin{equation}
w^{(s)}(t,\bx) = f(t) \star_t H(t) \sum_{l=1}^\infty  \cos \big(t \sqrt{\la_l}\big) y_l(\bx_s) y_l(\bx),
\label{eq:S4}
\end{equation}
where $H(t)$ is the Heaviside step function.

To derive the initial value problem \eqref{eq:F1}--\eqref{eq:F3}, we consider the even extension in time of this wave. 
Starting from equation \eqref{eq:S4}, using the Fourier transform formula
 \[
\int_{-\infty}^\infty dt H(t) \cos(t \sqrt{\la_l}) e^{i \om t} = 
\frac{\pi}{2} \Big[ \delta(\om - \sqrt{\la_l}) + \delta(\om + \sqrt{\la_l}) \Big] + 
\frac{i\om}{\la_l - \om^2},
\]
and the assumption that $f(t)$ is real valued, which means in light of \eqref{eq:S3} that $\hat f(\om) = \hat f(-\om)$, 
we obtain the following expression of the even time extension
\begin{align}
w_e^{(s)}(t,\bx) &= w^{(s)}(t,\bx) + w^{(s)}(-t,\bx) = 
\sum_{l=1}^\infty \hat f\big(\sqrt{\la_j}\big) \cos  \big(t \sqrt{\la_l}\big) y_l(\bx_s) y_l(\bx) \nonumber \\
&= \Big[\cos\big(t \sqrt{A}\big)\hat f\big(\sqrt{A}\big) \delta(\cdot-\bx_s)\Big](\bx),
\label{eq:S6}
\end{align}
where we use the standard definition of functions of self-adjoint operators.

The data are the matrices $\bD_j = \left( D_j^{(r,s)} \right)_{1 \le r,s \le m}$ with entries defined by this wave 
evaluated at the receivers,
\begin{align}
D_j^{(r,s)} &= w_e^{(s)}(j \tau, \bx_r) = 
\int_{\Omega} d \bx \delta(\bx-\bx_r)\Big[\cos\big(j \tau \sqrt{A }\big)\hat f\big(\sqrt{A}\big) \delta(\cdot-\bx_s)\Big](\bx),
\end{align}
for $j = 0, \ldots, 2n-1.$ 
We can rewrite them in the symmetric form \eqref{eq:F7}, in terms of  the sensor functions 
\begin{align}
b^{(s)}(\bx) &= \left[ \hat f^{\frac{1}{2}} \big(\sqrt{A}\big) \delta(\cdot-\bx_s)\right](\bx), \label{eq:S8}
\end{align}
using the commutation relations 
\[
\cos\big(t \sqrt{A}\big)\hat f\big(\sqrt{A}\big) = \hat f^{\frac{1}{2}} \big(\sqrt{A}\big)\cos\big(t \sqrt{A}\big)
 \hat f^{\frac{1}{2}} \big(\sqrt{A}\big).
\]

Note from equation \eqref{eq:S6} that at time $t = 0$,
\begin{equation}
w_e^{(s)}(0,\bx) = 2 w^{(s)}(0,\bx) = \Big[\hat f\big(\sqrt{A}\big) \delta(\cdot-\bx_s)\Big](\bx) .
\end{equation}
From equation \eqref{eq:S.1}, the homogeneous initial condition \eqref{eq:S.2}, the finite speed of propagation and the causality of the 
wave we know  that $w^{(s)}(0,\bx)$ is supported in the immediate vicinity of $\bx_s$.
The sensor function is just like it, but for a different pulse with Fourier transform $\hat f^{\frac{1}{2}}$. Therefore, it is 
supported near $\bx_s$, as stated below equation \eqref{eq:F6}. Causality also implies that $w_e^{(s)}(0,\bx)$ and therefore 
$b^{(s)}(\bx)$ are not affected by the medium outside the vicinity of $\bx_s$. Therefore, if the medium is known near the sensors, as is usually the case, the functions $b^{(s)}(\bx)$ can be calculated. This is why we treat them as known throughout the paper.

\section{Proof of Theorem \ref{thm.1}}
\label{ap:DatFit}

Equation~\eqref{eq:chebROM} follows from the time stepping scheme (\ref{eq:ROM28}--\ref{eq:ROM30}),
which is the three term recurrence relation for Chebyshev polynomials 
\begin{equation}
\begin{array}{rcl}
\cT_j(z) & = & 2 z \cT_{j-1}(z) - \cT_{j-2}(z), \quad j \geq 1, \\
\cT_0(z) & = & z^0, \\
\cT_{-1}(z) & = & \cT_{1}(z),
\end{array}
\label{eq:ttr}
\end{equation}
valid for any argument $z$.

To prove~\eqref{eq:ROMThm1}, we observe that \eqref{eq:Sn6} implies that the approximation subspace
\begin{equation}
\begin{array}{rcl}
\mathfrak{X} & = & \mbox{colspan}\{\cT_j\big(\cPq\big) \bb(\bx), ~j = 0, \ldots, n-1\} \\
& = & \mbox{colspan}\{ \bb(\bx), \cPq \bb(\bx), \ldots, \cPq^{n-1} \bb(\bx) \} 
\end{array}
\end{equation}
is a block Krylov subspace. Since 
$\mathfrak{X} = \mbox{range} \big( \bU(\bx)\big) = \mbox{range} \big( \bV(\bx)\big)$,
any polynomial $\cQ_i(\cPq) \bb(\bx)$ of degree $i \leq n-1$ is represented exactly in $\mathfrak{X}$, i.e.,
\begin{equation}
\cQ_i(\cPq) \bb(\bx) = \bV(\bx) \cQ_i(\cPqR) \bbR, \quad 0 \leq i \leq n-1.
\label{eq:poly}
\end{equation}

Any Chebyshev polynomial $\cT_j(z)$ of degree $j=0,\ldots,2n-1$ can be represented uniquely 
(via polynomial division) as 
\begin{equation}
\cT_j(z) = \cQ_i(z) \cT_n(z) + \cR_k(z),
\label{eq:qtr}
\end{equation}
for some polynomials $\cQ_i(z)$, $\cR_k(z)$ of degrees $i,k \leq n-1$. 
Setting $z = \cPq$, and using the facts that $\cQ_i(\cPq)$ is self-adjoint 
and $\cQ_i(\cPqR)$ is symmetric, we obtain
\[\begin{array}{rcl}
\bD_j & \stackrel{\eqref{eq:dataT}}{=} & \lb \bb, \cT_j(\cPq) \bb \rb \\
& \stackrel{\eqref{eq:qtr}}{=}  & \lb \bb, \cQ_i(\cPq) \cT_n(\cPq) \bb \rb + \lb \bb, \cR_k(\cPq) \bb \rb \\
& \stackrel{\eqref{eq:ttr}}{=} & \lb \cQ_i(\cPq) \bb, [2 \cPq \cT_{n-1}(\cPq) - \cT_{n-2}(\cPq)] \bb \rb  \\
& & + \lb \bb, \cR_k(\cPq) \bb \rb \\
& \stackrel{\eqref{eq:poly}}{=} & 
\lb \bV \cQ_i(\cPqR) \bbR, [2 \cPq \bV \cT_{n-1}(\cPqR) - \bV \cT_{n-2}(\cPqR)] \bbR \rb\\ 
& & + \lb \bb, \bV \cR_k(\cPqR) \bbR \rb \\
& \stackrel{\eqref{eq:defVort}}{=} & \bbR^T \cQ_i(\cPqR) [2 \cPqR \cT_{n-1}(\cPqR) - \cT_{n-2}(\cPqR)] \bbR \\ 
& & + \;\bbR^T \cR_k(\cPqR) \bbR \\
& \stackrel{\eqref{eq:ttr}}{=} & \bbR^T [\cQ_i(\cPqR) \cT_n(\cPqR) + \cR_k(\cPqR) ] \bbR \\
& \stackrel{\eqref{eq:qtr}}{=} & \bbR^T \cT_j(\cPqR) \bbR,
\end{array}\]
for all $j=0,\ldots,2n-1$.

\section{Proof of Theorem \ref{thm.2}}

\label{sect:PfThm2}
The symmetry of the ROM propagator follows immediately from equation~\eqref{eq:ROM21}, because $\cPq$ is self-adjoint. 

To prove that $\cPqR$ is block-tridiagonal, it suffices to show 
\begin{equation}
\cPqR_{j+l,j} = \lb \bv_{j+l},\cPq \bv_j \rb = 0, \qquad \forall \, l = 2, \ldots, n-j-1. 
\label{eq:ProofTridiag}
\end{equation}
From definition~\eqref{eq:defVj} of the orthogonal snapshots and the fact that the inverse $\bR^{-1}$ of the block upper triangular $\bR$ is also
block upper triangular, we get 
\begin{equation*}
\bv_j(\bx) = \sum_{i=0}^{j} \bu_i(\bx) \bR^{-1}_{i,j}, \qquad j = 0, \ldots, n-1.
\end{equation*}
We also have from the time stepping equation~\eqref{eq:Sn2} that 
\begin{align*}
\cPq \bv_j(\bx) &= \sum_{i=0}^{j} \cPq \bu_i(\bx) \bR^{-1}_{i,j} = \frac{1}{2} \Big[\sum_{i=0}^{j} \bu_{i+1}(\bx) \bR^{-1}_{i,j} + 
\sum_{i=0}^{j} \bu_{|i-1|}(\bx) \bR^{-1}_{i,j} \Big], 
\end{align*}
and therefore
\begin{align*} 
\cPqR_{j+l,j} &= \lb \bv_{j+l},\cPq \bv_j \rb = \frac{1}{2} \sum_{i=0}^{j} \Big[ \lb \bv_{j+l},\bu_{i+1} \rb + \lb \bv_{j+l}, \bu_{|i-1|} \rb 
\Big]  \bR^{-1}_{i,j} \\
&=\frac{1}{2} \sum_{i=0}^{j} \Big[ \bR_{j+l,i+1} + \bR_{j+l,|i-1|} \Big] \bR^{-1}_{i,j},
\end{align*}
where the last equality is because  $\bR = \bV^T \bU$. Since $\bR$ is block upper triangular, the right hand side in this equation is non-zero if the index $j \ge i$ satisfies $
j+l \le i+1 $ or $ j+l \le |i-1|.$
This is impossible for $l \ge 2$, so result~\eqref{eq:ProofTridiag} holds and $\cPqR$ is block tridiagonal.

It is clear from the definition~\eqref{eq:Sn5} of the propagator operator $\cPq$ that its eigenvalues must lie in the interval $[-1,1]$. Since 
$\cPqR $ is the Galerkin projection~\eqref{eq:ROM21} of the propagator, its eigenvalues also lie in $[-1,1]$. We now prove that 
\begin{equation}
\mbox{Ker} \big(I-\cPq\big) \bigcap \mathfrak{X} = \{{\bf 0}\},
\label{eq:ProofNull}
\end{equation}
which implies that 
\begin{equation*}
\bI_{nm}-\cPqR = \bV^T \big(I - \cPq \big) \bV 
\end{equation*}
is invertible.

Indeed, consider any element in $\mathfrak{X}$, written as 
\[
\sum_{j=0}^{n-1} \bu_j(\bx) \balpha_j \in \mathfrak{X},
\] 
for $m \times m$ diagonal matrices
$\balpha_j$, and suppose that it lies in the kernel of $\bI-\cPq$, 
\begin{equation}
(I - \cPq) \sum_{j=0}^{n-1} \bu_j(\bx) \balpha_j= 0.
\label{eq:LC1}
\end{equation}
We wish to show that  $\balpha_j = {\bf 0}$,  for $j = 0, \ldots, n-1$. 
Using~\eqref{eq:Sn2}--\eqref{eq:Sn4} in \eqref{eq:LC1}, we get
\[
\sum_{j=0}^{n-1} \Big[ \bu_j(\bx)  - \frac{\bu_{j+1}(\bx)+\bu_{j-1}(\bx)}{2} \Big]\balpha_j  = 0,
\]
and reordering the terms and using the initial condition~\eqref{eq:Sn9}  we have
\begin{align*}
&\bu_0(\bx)\Big(\balpha_0-\frac{\balpha_1}{2}\Big) +\bu_1(\bx) \Big(\balpha_1-\balpha_0 - \frac{\balpha_2}{2} \Big)+
 \bu_2(\bx) \Big(\balpha_2-\frac{\balpha_1 + \balpha_3}{2}\Big) + \ldots \\
&+\bu_{n-2}(\bx)\Big(\balpha_{n-2}-\frac{\balpha_{n-3} + \balpha_{n-1}}{2}\Big)
+ \bu_{n-1}(\bx)\Big(\balpha_{n-1}-\frac{\balpha_{n-2}}{2}\Big)- \bu_n(\bx)\frac{\balpha_{n-1}}{2} = 0.
\end{align*}
The wave snapshots are linearly independent up to time $n \tau$ by Assumption \ref{as.1}, so we can equate the
coefficients in this equation to 0. Starting with $\balpha_{n-1} = {\bf 0}$ and solving backward, we get that $\balpha_j = {\bf 0}$, for all $j = 0, \ldots, n-1$.
This shows that~\eqref{eq:ProofNull} holds and completes the proof of the theorem.
$\quad \Box$

\section{Proof of Theorem \ref{thm:5}}
\label{ap:PfThm5}

We obtain from definitions~\eqref{eq:Sn10},~\eqref{eq:ROM21} and~\eqref{eq:ROM44} that 
\begin{align}
\hspace{-0.15in}\frac{2}{\tau^2}(\bI_{nm}-\cPqR) &
= \bV^T \frac{2}{\tau^2}(I-\cPq)\bV = \bV^T \cLq \cLq^T \bV = \cLqR \cLqR^T,
\label{eq:ROM45}
\end{align}
where $\cLqR$ is an $nm \times nm$ block lower bidiagonal, invertible matrix by Theorem \ref{thm.2}.  We use it to define  the 
quasimatrix 
\begin{equation}
\label{eq:PG1}
\hat \bV(\bx) = \cLq^T \bV(\bx) \cLqR^{-T},
\end{equation}
and write
\begin{equation}
\cLqR = \bV^T \cLq \hat \bV.
\label{eq:PG1p}
\end{equation}
Note that we used in \eqref{eq:ROM45} the fact that $\cLqR^T = \hat \bV^T \cLq^T \bV$. This can be seen from
\begin{equation}
\begin{array}{rcl}
\big( \cLqR^T \bphi^\RM, \hat \bphi^\RM  \big) 
& = & \big( \bphi^\RM, \cLqR \hat \bphi^\RM \big) \\
& = & \big( \bphi^\RM, \bV^T \cLq \hat \bV \hat \bphi^\RM \big) \\
& = & \lb \bV \bphi^\RM, \cLq \hat \bV \hat \bphi^\RM \rb \\
& = & \lb \cLq^T \bV \bphi^\RM, \hat \bV \hat \bphi^\RM \rb \\
& = & \big( \hat \bV^T \cLq^T \bV \bphi^\RM, \hat \bphi^\RM \big),
\end{array}
\label{eq:CLQT}
\end{equation}
where $(\cdot, \cdot)$ is the inner product in $\mathbb{R}^{nm}$. Since \eqref{eq:CLQT} holds for any
$\bphi^\RM, \hat \bphi^\RM \in \mathbb{R}^{nm}$, we indeed have
\begin{equation}
\cLqR^T = \hat \bV^T \cLq^T \bV,
\label{eq:PG1h}
\end{equation}
a counterpart of \eqref{eq:PG1p}.

Returning to the quasimatrix $\hat \bV(\bx)$, we observe that it has orthonormal columns
\begin{align}
\hat \bV^T \hat \bV &= \cLqR^{-1}\bV^T \cLq \cLq^T \bV \cLqR^{-T} \hspace{0.05in} \hspace{-0.06in}\stackrel{\eqref{eq:ROM45}}{=}  \bI_{nm},
\label{eq:PG3}
\end{align}
and we now show that it satisfies the statement of the theorem. 

Recall from~\eqref{eq:ROM20} that the ROM snapshots $\bu^{\RM}_j$, for $j = 0, \ldots, n-1$, form the block upper triangular matrix $\bR$. Since $\cLqR^T$ is block 
upper bidiagonal, we get from~\eqref{eq:ROM56} that 
\[
\bhuR_0 = \frac{\tau}{2} \cLqR^T  \bbR = \frac{\tau}{2} \cLqR^T \begin{pmatrix} \bR_{0,0} \\ \vdots \\ {\bf 0} \end{pmatrix} = 
\begin{pmatrix} \hat \bR_{0,0} \\ \vdots \\ {\bf 0} \end{pmatrix},
\]
where the right hand side defines the $m \times m$ matrix $\hat \bR_{0,0}$.
 The next dual snapshot is obtained from equation~\eqref{eq:ROM53},
\[
\bhuR_1 = \bhuR_0 + \tau \cLqR^T \bu^{\RM}_1 = \begin{pmatrix} \hat \bR_{0,0} \\ \vdots \\ {\bf 0} \end{pmatrix}+ \tau \cLqR^T \begin{pmatrix} \bR_{0,1}\\
\bR_{1,1} \\ {\bf 0} \\ \vdots \\ {\bf 0} \end{pmatrix} = \begin{pmatrix} \hat \bR_{0,1}\\
\hat \bR_{1,1} \\ {\bf 0} \\ \vdots \\ {\bf 0} \end{pmatrix}
\]
and continuing this way we get ~\eqref{eq:BHU3}, with block upper triangular $\hat \bR$.

Next, we show that 
\begin{equation}
\bhu_0 \in \mbox{range}(\hat{\bV}).
\label{eq:PG4}
\end{equation}
 Indeed, using that $\hat{\bV} \hat{\bV}^T$ is the orthogonal projector  on $\mbox{range}(\hat{\bV})$, we calculate
\begin{align*}
\hat{\bV} \hat{\bV}^T  \bhu_0(\bx) &\stackrel{\eqref{eq:PG1}}{=}  \cLq^T \bV \cLqR^{-T} \cLqR^{-1} \bV^T \cLq \bhu_0(\bx) \\
&\stackrel{\eqref{eq:ROM51}}{=}\frac{\tau}{2}  \cLq^T \bV \cLqR^{-T} \cLqR^{-1} \bV^T \cLq \cLq^T \bb(\bx) \\
&\hspace{0.06in}=\frac{\tau}{2}  \cLq^T \bV \cLqR^{-T} \cLqR^{-1} \bV^T \cLq \cLq^T \bV \bV^T \bb(\bx) 
\end{align*}
where the last equality is because $\bV \bV^T$ is the orthogonal projector on the space~\eqref{eq:F5}
to which $\bb$ belongs.  The right hand side simplifies by equation~\eqref{eq:ROM45},  and ~\eqref{eq:PG4} holds because
\begin{align*}
\hat{\bV} \hat{\bV}^T  \bhu_0(\bx) &= \frac{\tau}{2}  \cLq^T \bV \cLqR^{-T} \cLqR^{-1} \cLqR \cLqR^T \bV^T \bb(\bx) \\
&\hspace{0.06in}=\frac{\tau}{2}  \cLq^T \bV  \bV^T \bb(\bx) =\frac{\tau}{2}  \cLq^T \bb(\bx) \stackrel{\eqref{eq:ROM51}}{=}\bhu_0(\bx).
\end{align*}
Furthermore, we have 
\begin{equation}
\bhu_0(\bx) = \hat{\bV}(\bx) \bhuR_0,
\label{eq:PG5}
\end{equation}
because
\begin{align*}
\hat{\bV} (\bx) \bhuR_0 &\stackrel{\eqref{eq:ROM56}}{=} \frac{\tau}{2} \hat{\bV}(\bx) \cLqR^T \bbR \\
&\stackrel{\eqref{eq:PG1p}}{=} \frac{\tau}{2} \hat{\bV} \hat{\bV}^T \cLq^T \bV \bbR \\
&\stackrel{\eqref{eq:ROM22}}{=} \frac{\tau}{2} \hat{\bV} \hat{\bV}^T \cLq^T \bV \bV^T \bb(\bx) \\
&\hspace{0.06in}=\frac{\tau}{2}   \hat{\bV} \hat{\bV}^T \cLq^T \bb(\bx) \\
&\stackrel{\eqref{eq:ROM51}}{=} \hat{\bV} \hat{\bV}^T \bhu_0(\bx) \stackrel{\eqref{eq:PG4}}{=}  \bhu_0(\bx). 
\end{align*}

Equations~\eqref{eq:ROM20},~\eqref{eq:defVort} and~\eqref{eq:ROM53} give
\begin{align*}
\bu_j(\bx) &= \bV \bu^{\RM}_j = \bV \cLqR^{-T}  \Big(\frac{\bhuR_j-\bhuR_{j-1}}{\tau} \Big), \quad j = 0, \ldots, n-1,
\end{align*}
and thefore, by~\eqref{eq:ROM48},
\begin{align*}
\frac{\bhu_j(\bx)-\bhu_{j-1}(\bx)}{\tau} &= \cLq^T \bu_j(\bx) = \cLq^T \bV(\bx) \cLqR^{-T} \Big(\frac{\bhuR_j-\bhuR_{j-1}}{\tau} \Big)
\\ &\hspace{-0.06in} \stackrel{\eqref{eq:PG1p}}{=} \hat{\bV}(\bx) \Big(\frac{\bhuR_j-\bhuR_{j-1}}{\tau}\Big), \qquad j = 0, \ldots, n-1.
\end{align*}
Starting with~\eqref{eq:PG5},  this implies  
that 
\[\bhu_j(\bx) = \hat{\bV}(\bx)  \bhuR_j, \qquad j =0, \ldots, n-1. \qquad \Box
\]

\section{Proof of Theorem \ref{thm.7}}
\label{ap:A}

The block-Lanczos iteration \cite[Chapter 4]{golubVanLoan} carried out for the skew adjoint operator 
\begin{equation}
\mathfrak{L}(q) = \begin{pmatrix} 0 & - \cLq \\ \cLq^T & 0 \end{pmatrix},
\label{eq:Y12}
\end{equation}
with a starting vector $[\bvc_0^{T}(\bx) ; {\bf 0}]^{T}$ generates the quasimatrix 
\begin{equation}
\hspace{-0.1in}\boldsymbol{\mathscr{V}}(\bx) = \begin{pmatrix} 
\bvc_0(\bx) & {\bf 0} & \bvc_1(\bx) & {\bf 0} &\ldots &\bvc_{n-1}(\bx) & {\bf 0} \\
{\bf 0} & \hat \bvc_0(\bx) & {\bf 0} & \hat \bvc_1(\bx)  &\ldots &{\bf 0} & \hat \bvc_{n-1}(\bx)
\end{pmatrix}
\end{equation}
with $2n$ block columns written in terms of  some orthonormal snapshots of the form
~\eqref{eq:Y9} and~\eqref{eq:Y10} that we wish to find. The Lanczos iteration calculates these snapshots so that  
\begin{align}
\mathfrak{L}(q) \boldsymbol{\mathscr{V}}(\bx)  = \boldsymbol{\mathscr{V}}(\bx) \tilde{\mathfrak{L}}(q) + 
\begin{pmatrix} 
 {\bf 0} & \ldots &{\bf 0} & {\itbf r}(\bx)\\
{\bf 0} & \ldots &{\bf 0} & {\bf 0},  
\end{pmatrix}
\label{eq:Y13}
\end{align}
where $ \tilde{\mathfrak{L}}(q)$ is $2nm \times 2nm$ block tridiagonal, skew-symmetric. Its diagonal consists of zero $m \times m$ blocks and the upper diagonal is, in the MATLAB notation,  
\begin{align}
\mbox{diag}(\tilde{\mathfrak{L}}(q),1) = \left( -\boldsymbol{\Lambda}^{\RM}_{0,0}(q), {\boldsymbol{\Lambda}^{\RM}_{1,0}}(q)^T,-\boldsymbol{\Lambda}^{\RM}_{1,1}(q), {\boldsymbol{\Lambda}^{\RM}_{2,1}}(q)^T, 
  \ldots, -\boldsymbol{\Lambda}^{\RM}_{n-1,n-1}(q) \right),
\label{eq:Y14}
\end{align}
where $\boldsymbol{\Lambda}^{\RM}_{i,j}(q)$ are the $m \times m$ blocks of~\eqref{eq:Y11}. 
The last term in~\eqref{eq:Y13} is the residual quasimatrix, with the single $m \times m$ non-zero block  ${\itbf r}(\bx)$. 
In this proof we will relate the $\bvc_j$ and $\hat\bvc_j$ to the fields $\bphi_{j}$ and $\hat{\bphi}_{j}$ 
from~\eqref{eq:St3} and \eqref{eq:St4}. Further, the entries of the matrix $\tilde{\mathfrak{L}}$ will be related 
to the block-finite difference coefficients $\bGa_j$ and $\hat \bGa_j$.

Equating the left and right hand sides of the Lanczos decomposition in~\eqref{eq:Y13} block column-wise, 
we obtain the following recursion scheme
\begin{align}
\cLq \hat \bvc_j(\bx) &=\bvc_{j}(\bx) \boldsymbol{\Lambda}^{\RM}_{j,j}(q) + 
\bvc_{j+1}(\bx)  \boldsymbol{\Lambda}^{\RM}_{j+1,j}(q),  \label{eq:Y15} \\
\cLq^T\bvc_j(\bx) &= \hat \bvc_{j-1}(\bx) \boldsymbol{\Lambda}^{\RM}_{j,j-1}(q)^T + 
\hat \bvc_j(\bx)  \boldsymbol{\Lambda}^{\RM}_{j,j}(q)^T,  \label{eq:Y16} 
\end{align}
for $j = 0, \ldots, n-1$, where 
\[
\hat \bvc_{-1}(\bx) = {\bf 0}, \qquad \bvc_n(\bx) \boldsymbol{\Lambda}^{\RM}_{n,n-1}(q) = -{\itbf r}(\bx).
\]
The matrices $\boldsymbol{\Lambda}^{\RM}$ follow from the normalization and orthogonality conditions for 
$\bvc_j$ and $\hat\bvc_j$.  We note that the recursion relations in~\eqref{eq:Y15} and \eqref{eq:Y16} resemble 
the recursion relations of a finite difference time-stepping scheme.
We wish to write $\bvc_j(\bx)$ and $\hat \bvc_j(\bx)$ in the form 
\begin{equation}
\bvc_j(\bx) = \bphi_j(\bx) \sqrt{\hat \bga_j}, \qquad \hat \bvc_j(\bx) = \hat \bphi_j(\bx) \sqrt{\bga_j},
\label{eq:Y18}
\end{equation}
for some arbitrary choice of the square roots
\begin{equation}
\label{eq:SQRT}
\bga_j =  \sqrt{\bga_j}  \sqrt{\bga_j}^{\,T}, \qquad \bhga_j =  \sqrt{\bhga_j}  \sqrt{\bhga_j}^T, \qquad j \ge 0.
\end{equation}
For any such choice we have,
\begin{align*}
\lb \bvc_j, \bvc_j \rb &\stackrel{ \eqref{eq:Y18}}{=} \sqrt{\bhga_j}^{\,T} \lb \bphi_j,\bphi_j\rb \sqrt{\bhga_j} \stackrel{ \eqref{eq:St2}}{=} 
\sqrt{\bhga_j}^{\,T}  \bhga_j^{-1} \sqrt{\bhga_j} \stackrel{ \eqref{eq:SQRT}}{=} \bI_m \\
\lb \hat \bvc_j, \hat \bvc_j \rb &= \sqrt{\bga_j}^{\,T} \lb \bhphi_j,\bhphi_j\rb \sqrt{\bga_j} =  \sqrt{\bga_j}^{\,T}  \bga_j^{-1} \sqrt{\bga_j} = \bI_m, 
\end{align*}
for $j =0,\ldots, n-1$, so the columns of $\bvc_j(\bx)$ and of $\hat \bvc_j(\bx)$ are orthogonal, as needed.

Substituting~\eqref{eq:St1} in (\ref{eq:St3}--\ref{eq:St4}), we get the equations 
\begin{align}
\left[ \bvc_{j+1}(\bx) \sqrt{\bhga_{j+1}}^{\, -1} - \bvc_{j} (\bx) \sqrt{\bhga_{j}}^{\, -1} \right] \bga_j^{-1}&= -\cLq \hat \bvc_j(\bx) \sqrt{\bga_j}^{-1} , 
\label{eq:Y20}\\
\left[ \hat \bvc_{j} (\bx) \sqrt{\bga_{j}}^{-1} - \hat \bvc_{j-1} (\bx) \sqrt{\bga_{j-1}}^{-1} \right] \bhga_j^{-1} &= \cLq^T \bvc_j (\bx) \sqrt{\bhga_j}^{-1} ,\label{eq:Y21}\
\end{align}
which must be consistent with (\ref{eq:St3}--\ref{eq:St4}) and the orthogonality of  $\{\bvc_j(\bx)\}_{0 \le j \le n-1}$ and
of $\{\hat \bvc_j(\bx)  \}_{0 \le j \le n-1}$.  Equations~\eqref{eq:Y15} and~\eqref{eq:Y20} are consistent if the blocks of $\boldsymbol{\Lambda}^{\RM}(q)$ satisfy 
\begin{align}
\boldsymbol{\Lambda}^{\RM}_{j,j}(q) &\stackrel{\eqref{eq:Y15}}{=} \lb \bvc_j, \cLq \hat \bvc_j\rb \stackrel{\eqref{eq:Y20}}{=}   \lb \bvc_{j}, \bvc_{j} \rb \sqrt{\bhga_{j}}^{\,-1} \bga_j^{-1}\sqrt{\bga_j} =\sqrt{ \bhga_{j}}^{-1} \sqrt{\bga_j}^{\,-T}, 
\label{eq:Ljj}
\end{align}
and 
\begin{align}
\boldsymbol{\Lambda}^{\RM}_{j+1,j}(q) &\stackrel{\eqref{eq:Y15}}{=} \lb \bvc_{j+1}, \cLq \hat \bvc_j\rb \stackrel{\eqref{eq:Y20}}{=}  -\lb \bvc_{j+1}, \bvc_{j+1} \rb \sqrt{\bhga_{j+1}}^{\,-1} \bga_j^{-1}\sqrt{\bga_j} \nonumber \\&~ = -\sqrt{\bhga_{j+1}}^{\, -1} \sqrt{\bga_j}^{\,-T}.
\label{eq:Ljp1}
\end{align}
The consistency of equations ~\eqref{eq:Y16} and~\eqref{eq:Y21},  which involve  the transposed blocks of $\boldsymbol{\Lambda}^{\RM}(q)$, follows the same way.

Next, we relate $\bv_j$ to $\bvc_{j}$ by showing that the orthonormal vectors  \eqref{eq:Y18} satisfy the 
relations~\eqref{eq:Y9} and~\eqref{eq:Y10}. To this effect, we note from definitions~\eqref{eq:ROM12},
\eqref{eq:ROM23},~\eqref{eq:defV},~\eqref{eq:St2} and~\eqref{eq:St1} that 
\begin{equation}
\label{eq:v0nu0}
\bv_0(\bx) = \bu_0(\bx) \bR_{0,0}^{-1} = \bb(\bx) \bR_{0,0}^{-1}, \qquad 
\bvc_0(\bx)  \stackrel{\eqref{eq:Y18}}{=} \bphi_0(\bx) \sqrt{\bhga_0} \stackrel{\eqref{eq:St5}}{=} \bb(\bx)\sqrt{\bhga_0}.
\end{equation}
Therefore, the columns of $\bvc_0(\bx)$ and $\bv_0(\bx)$ are an orthonormal basis of the same space $\mbox{span}\{\bu_0(\bx)\}$, 
so the two must be related by an orthogonal transformation 
$ \bY_0 \in \RR^{m \times m}$, 
\begin{equation}
\bvc_0(\bx) = \bv_0(\bx) \bY_0.
\end{equation}
We also get from equations \eqref{eq:BHU2}, \eqref{eq:ROM51} and \eqref{eq:St4} evaluated at $j = 0$ that 
\[
\hat \bv_0(\bx) = \bhu_0 \hat \bR_{0,0}^{-1}, \qquad \bhphi_0(\bx) = \cLq^T \bb(\bx) \bhga_0 = \frac{2}{\tau} \bhu_0(\bx) \bhga_0,
\]
so the columns of $\hat \bvc_0(\bx)$ and $\hat \bv_0(\bx)$ are orthonormal bases of the same space  $\mbox{span}\{\hat \bu_0(\bx)\}$,
and  must be related by an orthogonal transformation $\hat \bY_0 \in \RR^{m \times m}$, 
\begin{equation}
\hat \bvc_0(\bx) = \hat \bv_0(\bx) \hat \bY_0.
\end{equation}
Then, equations \eqref{eq:ROM19} and 
\eqref{eq:ROM47} and~\eqref{eq:St3} for $j = 0$ give that 
\[
\bv_1(\bx) \in \mbox{span} \{\bu_0(\bx),\bu_1(\bx)\}, \qquad 
\bvc_1(\bx) \in \mbox{span} \{\bu_0(\bx),\bu_1(\bx)\},
\]
so the columns of $\bvc_1(\bx)$ and $\bv_1(\bx)$ are orthonormal bases of the same space, 
the orthogonal complement of $\mbox{span}\{\bu_0(\bx)\}$ in $\mbox{span} \{\bu_0(\bx),\bu_1(\bx)\}.$ 
Therefore, they must be related by an orthogonal transformation $\bY_1 \in \RR^{m \times m}$, 
\begin{equation}
\bvc_1(\bx) = \bv_1(\bx)  \bY_1.
\end{equation}
Iterating this way we obtain the relations \eqref{eq:Y9} and \eqref{eq:Y10}. 

Using these orthogonal block diagonal transformations $\bY = \mbox{diag}\big(\bY_0, \ldots, \bY_{n-1}\big)$ and 
$\hat \bY = \mbox{diag}\big(\hat \bY_0, \ldots, \hat \bY_{n-1}\big)$, we can now define the matrices 
\begin{equation}
\bGa_j = \sqrt{\bga_j} \, \hat \bY_j^T, \qquad \bhGa_j = \sqrt{\bhga_j}\,  \bY_j^T, \qquad j \ge 0,
\label{eq:defGas}
\end{equation}
which are also square roots of $\bga_j$ and $\bhga_j$, 
\begin{align*}
\bGa_j \bGa_j^T &= \sqrt{\bga_j} \, \hat \bY_j^T \hat \bY_j \sqrt{\bga_j}^T =  
\sqrt{\bga_j} \sqrt{\bga_j}^T \stackrel{\eqref{eq:SQRT}}{=} \bga_j, \\
\bhGa_j \bhGa_j^T &= \sqrt{\bhga_j} \,  \bY_j^T \bY_j \sqrt{\bhga_j}^T =  
\sqrt{\bhga_j} \sqrt{\bhga_j}^T  \stackrel{\eqref{eq:SQRT}}{=}\bhga_j, \qquad j \ge 0.
\end{align*}
With these matrices the orthogonalized primary and dual snapshots $\bv_{j}$ and $\hat \bv_{j}$ can be shown 
to be transforms of $\bphi_j$ and $\hat \bphi_j$ from (\ref{eq:St3}--\ref{eq:St4}). We have 
\begin{align*}
\bv_j(\bx) &
\stackrel{\eqref{eq:Y9}}{=} \bvc_j(\bx) \bY_j^T 
\stackrel{\eqref{eq:Y18}}{=} \bphi_j(\bx) \sqrt{\bhga_j}\,  \bY_j^T = \bphi_j(\bx) \bhGa_j, \\
\hat \bv_j(\bx) &
\stackrel{\eqref{eq:Y10}}{=} \hat \bvc_j(\bx) \hat \bY_j^T 
\stackrel{\eqref{eq:Y18}}{=} \bhphi_j(\bx) \sqrt{\bga_j} \, \hat \bY_j^T = \bhphi_j(\bx) \bGa_j, 
\end{align*}
and the block lower bidiagonal ROM matrix $\cLqR$ follows from \eqref{eq:Y11} and (\ref{eq:Ljj}--\ref{eq:Ljp1}),
\begin{align*}
\cLqR_{j,j} &= \bY_j \boldsymbol{\Lambda}^{\RM}_{j,j}(q) \hat \bY_j^T \stackrel{\eqref{eq:Ljj}}{=} 
\bY_j \sqrt{ \bhga_{j}}^{-1} \sqrt{\bga_j}^{\,-T} \hat \bY_j^T  \stackrel{\eqref{eq:defGas}}{=} \bhGa_j^{-1} \bGa_j^{-T}, \\
\cLqR_{j+1,j} &= \bY_{j+1} \boldsymbol{\Lambda}^{\RM}_{j+1,j}(q) \hat \bY_j^T \stackrel{\eqref{eq:Ljj}}{=} 
-\bY_{j+1} \sqrt{ \bhga_{j+1}}^{-1} \sqrt{\bga_j}^{\,-T} \hat \bY_j^T  \stackrel{\eqref{eq:defGas}}{=} -\bhGa_{j+1}^{-1} \bGa_j^{-T},
\end{align*}
which allows an interpretation of the block entries of the ROM as block finite-difference coefficients.
This completes the proof of the theorem.
$\quad \Box$

\bibliography{biblio} \bibliographystyle{siam}

\end{document}